\documentclass[reqno,11pt,letterpaper]{amsart}

\usepackage{lipsum}
\usepackage{amsmath}
\usepackage{amssymb}
\usepackage{amsthm}
\usepackage{mathtools}
\usepackage{mathrsfs}
\usepackage{accents}
\usepackage{calc}
\usepackage{arydshln}
\usepackage{upgreek}
\usepackage{slashed}
\usepackage{xifthen}
\usepackage{graphicx}
\usepackage{longtable}
\usepackage[inline]{enumitem}
\usepackage{hyperref}
\usepackage{xcolor}
\usepackage{comment}

\usepackage[backend=biber, style=alphabetic, sorting=nyt]{biblatex}
\hypersetup{
  colorlinks 
}
\hypersetup{linkcolor=blue, citecolor=magenta}

\usepackage{titletoc}

\makeatletter

\def\@tocline#1#2#3#4#5#6#7{
\begingroup
  \par
    \parindent\z@ \leftskip#3 \relax \advance\leftskip\@tempdima\relax
                  \rightskip\@pnumwidth plus 4em \parfillskip-\@pnumwidth
    \ifcase #1 
       \vskip 0.6em \hskip 0em 
       \or
       \or \hskip 0em 
       \or \hskip 1em 
    \fi%
    #6
    \nobreak\relax{\leavevmode\leaders\hbox{\,.}\hfill}
    \hbox to\@pnumwidth {\@tocpagenum{#7}}
  \par
\endgroup
}

 \def\l@section{\@tocline{0}{0pt}{0pc}{}{}}

\renewcommand{\tocsection}[3]{%
  \indentlabel{\@ifnotempty{#2}{ 
    \ignorespaces\bfseries{#2. #3}}}
  \indentlabel{\@ifempty{#2}{\ignorespaces\bfseries{#3}}{}} 
    \vspace{1.5pt}}

\renewcommand{\tocsubsection}[3]{%
  \indentlabel{\@ifnotempty{#2}{
    \ignorespaces#2. #3}}
  \indentlabel{\@ifempty{#2}{\ignorespaces #3}{}}
    \vspace{1.5pt}}

\renewcommand{\tocsubsubsection}[3]{%
  \indentlabel{\@ifnotempty{#2}{
    \ignorespaces#2. #3}}
  \indentlabel{\@ifempty{#2}{\ignorespaces #3}{}}
    \vspace{1.5pt}}

\makeatother

\makeatletter
\def\@nomenstarted{0}
\newlength{\@nomenoldtabcolsep}

\newcommand{\nomenstart}
  {%
    \def\@nomenstarted{1}%
    \setlength{\@nomenoldtabcolsep}{\tabcolsep}%
    \setlength{\tabcolsep}{3.5pt}%
    \begin{longtable}{p{0.11\textwidth} p{0.86\textwidth}}
  }

\newcommand{\nomenitem}[2]{%
    \ifcase\@nomenstarted%
      \or 
      \or \\ 
    \fi%
    #1\,{\leavevmode\leaders\hbox{\,.}\hfill} & #2%
    \def\@nomenstarted{2}%
  }%
\newcommand{\nomenend}
  {\\%
      \end{longtable}%
      \setlength{\tabcolsep}{\@nomenoldtabcolsep}%
      \def\@nomenstarted{0}%
  }
\makeatother

\newcommand\numberthis{\addtocounter{equation}{1}\tag{\theequation}}

\DeclareMathSymbol{\widetildesym}{\mathord}{largesymbols}{"65}

\newcommand{\tild}
{\raise.17ex\hbox{$\scriptstyle\mathtt{\sim}$}}

\allowdisplaybreaks

\numberwithin{equation}{section}
\numberwithin{figure}{section}
\newtheorem{thm}{Theorem}[section]

\newtheorem{prop}[thm]{Proposition}
\newtheorem{lemma}[thm]{Lemma}
\newtheorem{cor}[thm]{Corollary}
\newtheorem*{thm*}{Theorem}
\newtheorem*{prop*}{Proposition}
\newtheorem*{cor*}{Corollary}
\newtheorem*{conj*}{Conjecture}

\theoremstyle{definition}
\newtheorem{definition}[thm]{Definition}

\newtheorem{example}[thm]{Example}
\AtBeginEnvironment{example}{%
  \pushQED{\qed}%
}
\AtEndEnvironment{example}{\popQED\endexample}

\newtheorem{remark}[thm]{Remark}
\AtBeginEnvironment{remark}{%
  \pushQED{\qed}%
}
\AtEndEnvironment{remark}{\popQED\endremark}

\makeatletter
\newcommand{\fakephantomsection}{%
  \Hy@MakeCurrentHref{\@currenvir.\the\Hy@linkcounter}
  \Hy@raisedlink{\hyper@anchorstart{\@currentHref}\hyper@anchorend}%
  \Hy@GlobalStepCount\Hy@linkcounter%
}
\makeatother

\DeclareMathOperator{\ind}{ind} 

\DeclareMathOperator{\End}{End}

\DeclareMathOperator{\rTr}{\prescript{R}{}Tr}
\DeclareMathOperator{\Tr}{Tr}
\DeclareMathOperator{\tr}{tr}

\newcommand{\del}{\partial}

\newcommand*{\dt}[1]{\accentset{\mbox{\large\bfseries .}}{#1}}

\DeclareSymbolFont{toneitalic}{T1}{\familydefault}{m}{it}
\DeclareMathSymbol{\dir}{\mathord}{toneitalic}{"F0}

\newcommand\mf\mathfrak
\newcommand\mc\mathcal
\newcommand\mb\mathbb

\makeatletter
\newcommand*\bigcdot{\mathpalette\bigcdot@{.5}}
\newcommand*\bigcdot@[2]{\mathbin{\vcenter{\hbox{\scalebox{#2}{$\m@th#1\bullet$}}}}}
\makeatother

\begin{document}

\title{The topology of end-periodic smoothing operators} 
\date{\today}
\author{Alex R. Taylor}

\begin{abstract}
We develop an algebra of smoothing operators adapted to the geometry of end-periodic manifolds. We extend the renormalized trace and trace defect formula of Mrowka-Ruberman-Saveliev to this algebra, and establish a Fredholm criterion and winding-number index formula for smoothing perturbations of the identity. Our main result computes the homotopy groups of the associated end-periodic smoothing group. This provides a foundation for future work on relative index theory for families of end-periodic differential operators.
\end{abstract}

\maketitle

\tableofcontents

\section{Introduction}\label{sec-0}


Our starting point is the end-periodic index theory of \cite{mrs1}, as well as the relative index theory for cusp operators developed in \cite{mr2}. We briefly recall the end-periodic setting. Let $Y$ be a closed manifold equipped with a primitive integral cohomology class, and let $\widetilde Y\to Y$ be the corresponding infinite cyclic cover. Cutting $Y$ along the dual hypersurface to this class gives the compact fundamental segment $W$ of the cover. An end-periodic manifold is obtained by attaching an infinite sequence of copies $W_{j}$ of the segment to a compact manifold $Z$:
\[
M=Z\cup W_0\cup W_1\cup W_2\cup\cdots.
\]
The geometric data on $M$ are obtained by lifting from $Y$ to $\widetilde{Y}$.



This setting gives rise to a rich Fredholm index theory, particularly for end-periodic Dirac operators \cite{tau, mrs1}. A general smoothing algebra adapted to this geometry is nevertheless missing from the framework used in these works. In \cite{mrs1}, they consider smoothing operators of the form $D^{m}e^{-tD^{2}}$ where $D$ is an end-periodic Dirac operator. For many purposes it is useful to have an intrinsic class of operators that is more expansive than heat kernels; for example, if one wants to study end-periodic differential operators not of Dirac-type, or construct Fredholm parametrices with specific mapping properties, or smoothing perturbations of end-periodic operators.


In this paper, we develop a class of smoothing operators adapted to the geometry of end-periodic manifolds. The algebra $\Psi_{\mathrm{ep}}^{-\infty}(M;E)$ consists of smoothing operators that agree, modulo a rapidly decaying residual term, with a $\mathbb{Z}$-invariant (i.e., periodic) smoothing operator along the end. It has several nice properties:
\begin{enumerate}
\item It is a Fr\'echet $*$-algebra under composition and adjoint, with a canonical topology independent of the auxiliary choices.

\item The normal operator defines a short exact sequence of Fr\'echet $*$-algebras
\[
0 \to \dt{\Psi}_{\mathrm{ep}}^{-\infty}(M;E)
\to \Psi_{\mathrm{ep}}^{-\infty}(M;E)
\xrightarrow{\;N\;}
\Psi_{\mathbb Z}^{-\infty}(\widetilde Y;\widetilde E) \to 0,
\]
whose kernel is the ideal of residual smoothing operators.

\item The Fourier-Laplace transform identifies the normal algebra with an algebra of smooth twisted $\mathbb{R}$-families of smoothing operators on the closed manifold $Y$.

\item The renormalized trace and trace defect formula introduced by \cite[Theorem 4.6]{mrs1} extend to $\Psi_{\mathrm{ep}}^{-\infty}(M;E)$ without any modification.
\end{enumerate}
Consequently the algebra provides a natural framework for the analysis of end-periodic operators and their index-theoretic invariants. Compare with similar well-known pseudodifferential calculi adapted to various geometric settings; for example the b-calculus of \cite{melrose}, the cusp calculus of \cite{mr2}, and that of \cite{aln}.

However, our goal is more limited than the construction of a full end-periodic pseudodifferential calculus: the smoothing algebra developed here is designed for the study of end-periodic \emph{differential} operators. This narrower scope permits a comparatively simple theory. We briefly describe the main aspects developed in this paper.

\medskip
\noindent\emph{The indicial family.} In the b- and cusp-calculus, the geometry at infinity is modeled using a cylinder $\mathbb{R}\times X$, where $\mathbb{R}$ acts by translation, and the normal operator is translation-invariant. Fourier transformation in the $\mathbb{R}$-direction gives a family of operators on $X$, parametrized by $\tau\in\mathbb R$. For cusp smoothing operators, the resulting \emph{indicial family} $A_{\tau}\in\mathcal S(\mathbb R;\Psi^{-\infty}(X))$ is Schwartz in $\tau$. Hence $I+A_\tau$ tends to the identity as $\tau\to\pm\infty$, giving a based loop upon one-point compactification of $\mathbb R$.

The situation is different for periodic ends. The normal operator is invariant under the discrete deck group $\mathbb Z$. The resulting indicial family depends smoothly on the Fourier-Laplace parameter $\theta \in \mathbb{R}$, but possesses no decay in that parameter. Moreover, the indicial family satisfies an extra twisting condition
\[
A_{\theta+2\pi}=CA_\theta C^{-1},\hspace{3mm}
(Cv)(p(y)) := e^{2\pi i f(y)}v(p(y)),
\]
where $f:\widetilde Y\to\mathbb R$ is a smooth function such that $f^{-1}[0,1] = W$ and $f(Ty)=f(y)+1$. Thus the indicial family is a \emph{twisted} loop, its endpoints are unitarily equivalent rather than being exactly equal. Theorem \ref{FL-transform-thm} identifies the algebra of normal operators as
\[
\Psi_{\mathbb Z}^{-\infty}(\widetilde Y;\widetilde E) \cong C^\infty_{\mathrm{FL}}
\bigl(\mathbb R;\Psi^{-\infty}(Y;E_Y)\bigr),
\]
where the right-hand side consists precisely of smooth families satisfying this twisting condition.

\medskip
\noindent\emph{Renormalized trace.} In \S\ref{sec-1}, we extend the renormalized trace of \cite{mrs1} to the end-periodic smoothing algebra and establish its trace defect formula in this setting. Compared with the cusp setting, the trace defect formula contains an additional term given by an $f$-weighted trace. When developing end-periodic versions of results for cusp smoothing operators, such as the transgression identity \cite[Proposition 5]{mr2}, this extra term must be accounted for.

\medskip
\noindent\emph{Fredholm criterion.} In \S\ref{sec-2}, for $A \in \Psi_{\mathrm{ep}}^{-\infty}(M;E)$ we establish the Fredholm criterion
\[
I+A\text{ is Fredholm} \, \Longleftrightarrow \, 
I+A_\theta\text{ is invertible for every }\theta\in\mathbb R,
\]
for end-periodic smoothing perturbations of the identity. In the Fredholm case, we establish the winding-number formula
\[
\operatorname{ind}(I+A) =-\operatorname{wind}\bigl(\det(I+A_\theta)\bigr).
\]
Thus the Fredholm index is calculated in terms of the indicial family.

\medskip
\noindent\emph{Topology of the smoothing group.} The main problem we investigate in \S\ref{sec-3} and \S\ref{sec-4} is to understand the topology of the group $G_{\mathrm{ep}}^{-\infty}(M;E)$ of invertible operators of the form $I+A$, with $A\in\Psi_{\mathrm{ep}}^{-\infty}(M;E)$. Put $\mathcal G:=G^{-\infty}(Y;E_Y)$, the ordinary smoothing group on $Y$. The end-periodic normal smoothing group is the twisted loop group $L_C\mathcal G$. The contrast with the cusp case can be summarized as follows:
\begin{itemize}
\item The ordinary free-loop group splits as $L\mathcal G\cong\mathcal G\times\Omega\mathcal G$ by evaluation and constant loops. By Bott periodicity, $\pi_j(L\mathcal G)\cong\mathbb Z$ for every $j\geq1$, and the based-loop subgroup $\Omega\mathcal G$ has homotopy group $\mathbb{Z}$ in even degrees.

\item In the cusp setting, the Schwartz decay of the indicial family means the suspended smoothing group can be viewed as a based-loop group. Hence the above argument determines its homotopy groups.
These groups are then used in the homotopy long exact sequence associated with the normal map fibration. By analyzing its connecting maps, Melrose-Rochon prove that the \emph{cusp smoothing group} is weakly contractible \cite[Proposition 6]{mr2}. 

\item In the end-periodic setting, there is no Schwartz decay, and the constant-loop splitting is unavailable. Nonetheless, we show that $\pi_j(L_C\mathcal G)\cong\pi_j(L\mathcal G)\cong\mathbb Z$ for every $j\geq1$. Using these groups in the homotopy long exact sequence determined by the normal map, we show that the even-degree connecting maps are isomorphisms, and thereby calculate the homotopy groups of $G_{\mathrm{ep}}^{-\infty}(M;E)$. 

\end{itemize}

Our main result (Theorem \ref{thm:homotopy-type-Gep}) says that evaluation of the indicial family at $\theta=0$ is a weak homotopy equivalence
\[
\varepsilon:G_{\mathrm{ep}}^{-\infty}(M;E)
\to \mathcal G, \hspace{3mm} g\mapsto g_0.
\]
Consequently, $G_{\mathrm{ep}}^{-\infty}(M;E)$ is path connected and
\[
\pi_{2k}\bigl(G_{\mathrm{ep}}^{-\infty}(M;E)\bigr)=0, \hspace{3mm} \pi_{2k-1}\bigl(G_{\mathrm{ep}}^{-\infty}(M;E)\bigr) \cong\mathbb Z,\hspace{3mm} k\geq1.
\]
The full end-periodic smoothing group is therefore not weakly contractible. However, its based subgroup
\[
G_{\mathrm{ep},*}^{-\infty}(M;E)
:=\{g\in G_{\mathrm{ep}}^{-\infty}(M;E):g_0=I\}
=\ker\varepsilon
\]
is weakly contractible (Corollary \ref{cor:based-ep-weak-contractibility}).

The study of these homotopy groups is motivated by the strategy used by Melrose-Rochon \cite{mr2} for proving their relative index theorem \cite[Theorem 2]{mr2}. Their weak contractibility theorem for the cusp smoothing group \cite[Proposition 6]{mr2} is a key input in proving that smooth families of elliptic cusp operators admit invertible smoothing perturbations \cite[Theorem 1]{mr2}, which in turn underlies their relative index theorem \cite[Theorem 2]{mr2}.

This strategy first requires an appropriate operator calculus and smoothing algebra. In this paper we develop some of these foundational elements. In subsequent work, we intend to use the weakly contractible group $G_{\mathrm{ep},*}^{-\infty}$ to construct invertible perturbations of end-periodic families, and to establish a relative index theorem for end-periodic families.

\smallskip
\noindent \textbf{Acknowledgements.} The author would like to thank Pierre Albin for many helpful conversations and valuable suggestions during the preparation of this manuscript. 

\section{End-periodic smoothing calculus}\label{sec-1}

We follow the end-periodic geometric framework of \cite[\S 2.1]{mrs1}. Let $(Y,g_{Y})$ be a smooth, closed, connected Riemannian manifold, and fix a primitive class in $H^1(Y;\mathbb Z)$. It determines an infinite cyclic covering $p:\widetilde Y\to Y$ with a chosen positive deck generator $T$. Choose a smooth map $h:Y\to S^1$ representing the primitive class and its real-valued lift $f:\widetilde Y\to\mathbb R$, satisfying
\begin{equation}\label{f-definition}
h(p(y))=e^{2\pi i f(y)},\hspace{3mm} f(Ty)=f(y)+1.
\end{equation}
We also write $y+k:=T^{k}y$ for $k\in\mathbb Z$. The differential $df$ is invariant under $T$ and therefore descends to a closed one-form on $Y$.

Let $W:=f^{-1}([0,1])\subset\widetilde Y$. Thus $W$ is a compact manifold with two boundary components. We call $W$ the \emph{fundamental segment} for the periodic end. Writing $W_j:=T^jW$, we have
\[
\widetilde Y=\bigcup_{j\in\mathbb Z}W_j, \hspace{3mm} \widetilde Y_+:=\bigcup_{j\geq0}W_j.
\]
An \emph{end-periodic manifold} $M$ is obtained by attaching $\widetilde Y_+$ to a compact manifold $Z$ along the incoming boundary of $W_0$:
\[
M=Z\cup W_0\cup W_1\cup\cdots, \hspace{3mm} M_N:=Z\cup\bigcup_{j=0}^N W_j.
\]
We equip $M$ with a Riemannian metric agreeing on the periodic end with the restriction of $p^*g_Y$. Similarly, a Hermitian bundle $E\to M$ is \emph{end-periodic} if, on the end, it is identified with $\widetilde E:=p^*E_Y$ for a Hermitian bundle $E_Y\to Y$; we require its Hermitian connection to agree there with the pullback of a connection on $E_Y$. The deck transformations act naturally on these pullback bundles and preserve their metrics and connections.

A differential operator $P:C^\infty(M;E)\to C^\infty(M;F)$ of order $m$ is \emph{end-periodic} if there is a differential operator $P(Y):C^\infty(Y;E_Y)\to C^\infty(Y;F_Y)$ of order $m$ such that, under the above identifications,
\[
P\big|_{\widetilde Y_+} =\bigl(p^*P(Y)\bigr)\big|_{\widetilde Y_+}.
\]
The periodic lift $\widetilde P:=p^*P(Y)$ is uniquely determined and is called the \emph{normal operator}, also denoted $N(P)=\widetilde P$. Using the end-periodic geometric data, we consider $L^{2}$-Sobolev spaces of sections $H^{m}(M;E)$, for $m \geq 0$. Then $P$ defines a bounded operator $P: H^m(M;E)\to L^2(M;F)$.

\subsection{Suspended algebra}\label{subsec-suspended-algebra}

An elliptic operator on a non-compact manifold need not be Fredholm.  As a general principle, in addition to symbolic invertibility, the normal operator which governs the behavior at infinity must also be invertible. This is true in the case of b- and cusp operators, and for other manifolds with nice ends; for example, the Lie manifolds studied in \cite{bcnq21}. In the cusp setting of \cite{mr2}, the normal operators live in the \emph{suspended algebra}.

We now introduce the suspended or normal smoothing algebra associated with the infinite cyclic cover $p:\widetilde Y\to Y$. By construction, any end-periodic smoothing operator will have a normal operator in this algebra. We use the pullback Hermitian connection on $\widetilde E=p^*E_Y$ to define covariant derivatives on $\widetilde E\boxtimes\widetilde E^*$.

\begin{definition}[End-periodic suspended smoothing operators]
Let
\[
\Psi^{-\infty}_{\mathbb Z}(\widetilde Y;\widetilde E)
\]
be the space of operators initially defined on
\[
P:C_{c}^{\infty}(\widetilde Y;\widetilde E)\to C^{\infty}(\widetilde Y;\widetilde E)
\]
whose Schwartz kernel $K_{P}\in
C^{\infty}(\widetilde Y\times\widetilde Y;
\widetilde E\boxtimes\widetilde E^{*})$ satisfies the following two conditions:

\begin{enumerate}
\item Diagonal periodicity:
\begin{equation}\label{eqn-sus-periodicity}
K_{P}(y+k,y'+k) = K_{P}(y,y'), \hspace{3mm} k\in\mathbb Z.
\end{equation}

\item Relative Schwartz decay along the periodic end: for any $N,r\geq 0$,
\begin{equation}\label{eqn-sus-schwartz}
\nu_{N,r}(P):= \sup_{y,y'\in \widetilde Y} \langle f(y')-f(y)\rangle^{N} \left| \nabla^{r}K_{P}(y,y')\right|<\infty
\end{equation}
where $\langle t\rangle:=(1+t^{2})^{1/2}$.
\end{enumerate}
We equip this space with the locally convex topology defined by the seminorms $\nu_{N,r}$. These are the Fr\'echet seminorms used below; we establish completeness in Lemma \ref{lem:Frechet-topology-Z-sus}.
\end{definition}

\begin{remark}
The decay in \eqref{eqn-sus-schwartz} is relative: it concerns $|f(y)-f(x)|\to\infty$, not the separate limits $f(x)\to\infty$ and $f(y)\to\infty$. A periodic kernel with Schwartz decay in both variables must vanish identically.
\end{remark}

\begin{remark}[Comparison with the cusp setting]
We refer to \cite{mr2} for the calculus of cusp pseudodifferential operators. For our purposes, the cusp smoothing algebra is most relevant. In this setting, the suspended smoothing algebra consists of operators $P \in \Psi_{\mathrm{sus}}^{-\infty}(Y)$ which are translation invariant in the sense that
\[
T_{\tau}^{*}Pu = PT_{\tau}^{*}u, \hspace{3mm} \tau \in \mb{R},
\]
where $T_{\tau}(y,t) = (y,t-\tau)$; moreover, one requires that its integral kernel lives in the space
\[
K_{P} \in C_{c}^{-\infty}(Y^{2}\times\mb{R}) + \mc{S}(Y^{2}\times\mb{R}),
\]
where $\mathcal{S}$ denotes the Schwartz space of functions rapidly decreasing, with all derivatives, at infinity.
The translation invariance implies that $P$ acts as a convolution operator in the $\mb{R}$-variable, so
\[
Pu(y,t) = \int_{Y}\int_{\mb{R}} K_{P}(y,y',t-s)u(y,'s)\,ds.
\]
Since $Y^{2}$ is compact, the only noncompact variable in the integral kernel $K_{P}(y,y',t-s)$ is the $\mathbb{R}$-variable $r = t-s$. Therefore, rapid decay at infinity means precisely that, for every $N,j\geq 0$ and every pair of multi-indices $\alpha,\beta$,
\[
\begin{split}
\sup_{\substack{r\in\mathbb R\\y,y'\in Y}}
\langle r \rangle^{N}\left|\partial_{r}^{j}\nabla_{y}^{\alpha}\nabla_{y'}^{\beta}
K_{P}(y,y',r)\right|\infty, \hspace{3mm}\langle r\rangle=(1+r^{2})^{1/2},
\end{split}
\]
away from the conormal singularity at $\{0\}\times\mathrm{Diag}(Y)$.

\medskip

In the end-periodic setting, the $\mathbb{R}$-coordinate $t$ is replaced by $f(y)$, which suggests that we should impose Schwartz decay as $|f(y)-f(y')|\to\infty$. The translation invariance is replaced by periodicity; i.e., invariance under the $\mathbb{Z}$-deck transformation on the periodic end. This is the motivation behind our definition of \emph{end-periodic suspended smoothing operators}.
\end{remark}

Let $Y$ be a closed manifold with infinite cyclic covering space $p:\widetilde Y \to Y$, positive deck generator $T$, Hermitian vector bundle $E_{Y}\to Y$, and smooth function $f$ such that
\begin{equation}\label{eqn-sus-family-height}
f(Ty)=f(y)+1.
\end{equation}
Define the unitary multiplication operator
\begin{align*}
& C: L^{2}(Y;E_{Y}) \to L^{2}(Y;E_{Y}) \\
& (Cv)(p(x)) := e^{2\pi i f(x)}v(p(x)).
\end{align*}
It is well defined since $e^{2\pi i f(Tx)}=e^{2\pi i f(x)}$. For $\theta\in\mb{R}$, the Fourier-Laplace (FL) transform is defined as
\begin{align*}
& \mathcal{F}_{\theta} : C_{c}^{\infty}(\widetilde{Y};\widetilde{E}) \to C^{\infty}(Y;E_{Y}) \\
& (\mathcal F_{\theta}u)(p(x)) := e^{i\theta f(x)} \sum_{k\in\mathbb Z} e^{ik\theta}u(x+k).
\end{align*}
The sum is finite for each $x$ and its value is independent of the lift $x$ of $p(x)$. Changing the real lift of $e^{i\theta}$ gives
\begin{equation}\label{eqn-sus-FL-vector-monodromy}
\mathcal F_{\theta+2\pi}u=C\mathcal F_\theta u.
\end{equation}

\noindent As $\theta$ varies over $\mathbb{R}$, we can regard $\mathcal{F}_{\theta}u$ as a section of an appropriate bundle over $S^{1}$. In order to incorporate the choice of real lift in an intrinsic way, we introduce the FL Hilbert bundle over $S^{1}$.

\begin{definition}[Fourier-Laplace Hilbert bundle]
\label{def:sus-FL-Hilbert-bundle}
The Fourier-Laplace Hilbert bundle over $S^1$ is the
mapping torus
\begin{equation}\label{eqn-sus-FL-gluing}
\mathcal H^{\mathrm{FL}}:= \frac{[0,2\pi]\times L^{2}(Y;E_{Y})}{(2\pi,w)\sim(0,C^{-1}w)} \to S^1.
\end{equation}
We write $[\theta,w]$ for an equivalence class lying over $e^{i\theta}$. The inner product is inherited from $L^{2}(Y;E_{Y})$. A smooth section of $\mathcal H^{\mathrm{FL}}$ is equivalently a smooth $\mathbb{R}$-family $w_{\theta}\in L^{2}(Y;E_{Y})$ satisfying
\begin{equation}\label{eqn-sus-FL-vector-equivariance}
w_{\theta+2\pi}=Cw_{\theta}, \hspace{3mm} \theta\in\mb{R}.
\end{equation}
\end{definition}

\noindent The FL transform depends on a real lift $\theta$ of $\xi=e^{i\theta}\in S^1$. Its change under $\theta\mapsto\theta+2\pi$ is precisely the relation \eqref{eqn-sus-FL-vector-equivariance}. This is why the family $\theta \mapsto \mathcal{F}_{\theta}u$ is naturally a section of $\mathcal H^{\mathrm{FL}}$, rather than an ordinary $L^{2}(Y;E_{Y})$-valued function on $S^{1}$.

Owing to the relation \eqref{eqn-sus-FL-vector-monodromy}, $\mathcal F u$ is a smooth section of the FL Hilbert bundle on $S^1$. With Haar measure $d\theta/(2\pi)$,
the Fourier series Plancherel theorem gives a unitary map
\begin{equation}\label{eqn-sus-FL-unitary}
\mathcal F:L^2(\widetilde Y;\widetilde E)
\xrightarrow{\ \cong\ } L^2(S^1;\mathcal H^{\mathrm{FL}}).
\end{equation}
For a smooth section $w_\theta$ satisfying \eqref{eqn-sus-FL-vector-equivariance}, the inverse formula is 
\begin{equation}\label{eqn-sus-FL-inverse-on-W}
 u(T^jx)=\frac1{2\pi}\int_0^{2\pi}
 e^{-ij\theta}e^{-i\theta f(x)}w_\theta(p(x))\,d\theta,
\end{equation}
evaluated for $x\in W$ and $j\in\mathbb Z$. The same formula holds in $L^2$.

One expects that conjugating the space $\Psi^{-\infty}_{\mathbb Z}(\widetilde Y;\widetilde E)$ by the FL transform should yield a natural space of operators on $Y$. Our next goal is to identify precisely what this space must be. A natural candidate is the smooth loop space $C^{\infty}(S^{1};\Psi^{-\infty}(Y;E_{Y}))$; however this is not quite correct, because conjugating $\widetilde{A} \in \Psi^{-\infty}_{\mathbb Z}(\widetilde Y;\widetilde E)$ by $\mathcal{F}$ yields a family of operators $(A_{\theta})_{\theta \in \mathbb{R}}$ satisfying
\[
A_{\theta+2\pi} = CA_{\theta}C^{-1}.
\]
Thus the family $A_{\theta}$ is not literally an $S^{1}$-family, instead changing the period by $2\pi$ yields unitarily equivalent operators. Thus, to address this issue we introduce the space
\begin{equation}\label{eqn-sus-twisted-algebra}
C^\infty_{\mathrm{FL}}\left(\mathbb{R}; \Psi^{-\infty}(Y;E_Y)\right) \subset C^{\infty}(\mathbb{R}; \Psi^{-\infty}(Y;E_Y)),
\end{equation}
consisting of smooth families $A_{\theta} \in C^{\infty}(\mathbb{R}; \Psi^{-\infty}(Y;E_Y))$ satisfying $A_{\theta+2\pi} = CA_{\theta}C^{-1}$ for all $\theta\in\mb{R}$.

\begin{example}[End-periodic Dirac operators]
Although we are primarily interested in smoothing operators, the conjugation point is illustrated by the case of an end-periodic Dirac operator. Fix an infinite cyclic covering space $p: \widetilde{Y} \to Y$ as above. Let $E_{Y} \to Y$ be an end-periodic Clifford module, let $D(Y)$ be the associated Dirac operator whose principal symbol is Clifford multiplication $\mathrm{cl}: T^{*}M \to \End(E_{Y})$, and let $\widetilde{D} = p^{*}D(Y)$ be the periodic Dirac operator acting on $\widetilde{E} \to \widetilde{Y}$. A straightforward computation shows that
\[
D_{\theta} := (\mathcal{F}\widetilde{D}\mathcal{F}^{-1})_{\theta} = D(Y) - i\theta\, \mathrm{cl}(df).
\]
Thus we obtain a family of Dirac operators $D_{\theta}$ on $Y$ parametrized by $\theta \in \mathbb{R}$, satisfying $D_{\theta+2\pi}= CD_{\theta}C^{-1}$. This is precisely the family in \cite[(9)]{mrs1}.
\end{example}

For any periodic operator $\widetilde{A} \in \Psi^{-\infty}_{\mathbb Z}(\widetilde Y;\widetilde E)$, conjugating with the FL transform yields a family of smoothing operators
\begin{equation}\label{eqn-sus-FL-intertwining}
A_{\theta} \in \Psi^{-\infty}(Y;E_{Y}), \hspace{3mm} A_\theta w_\theta := (\mathcal F\widetilde A\mathcal F^{-1}w)_\theta, \hspace{3mm} \theta\in\mathbb{R},
\end{equation}
for $w \in L^{2}(S^{1};\mathcal{H}^{\mathrm{FL}})$. The integral kernels are related by the formula
\begin{equation}\label{eqn-kernel}
K_{\widetilde A}(x,y) =\frac1{2\pi}\int_0^{2\pi}
e^{i\theta(f(y)-f(x))} K_{A_\theta}(p(x),p(y))\,d\theta, \hspace{3mm} x,y\in\widetilde Y.
\end{equation}
The family $A_{\theta}$ satisfies the twisting condition
\begin{equation}\label{eqn- sus-FL-operator-twist}
A_{\theta+2\pi}=CA_\theta C^{-1},
\end{equation}
and so defines an element $A_{(\cdot)} \in C^\infty_{\mathrm{FL}}(\mathbb{R};\Psi^{-\infty}(Y;E_Y))$. In fact, we will show in Theorem \ref{FL-transform-thm} that conjugation by $\mathcal{F}$ defines a $\ast$-algebra isomorphism 
\[
C^{\infty}_{\mathrm{FL}}(\mathbb{R};\Psi^{-\infty}(Y;E_{Y})) \xleftrightarrow{ \hspace{3mm} \cong \hspace{3mm}} \Psi_{\mathbb{Z}}^{-\infty}(\widetilde{Y};\widetilde{E}).
\]
Once these spaces are equipped with natural Fr\'echet space topologies, this is also an isomorphism of Fr\'echet spaces.

\subsection{End-periodic smoothing algebra}

Having defined the suspended smoothing algebra
$\Psi_{\mathbb Z}^{-\infty}(\widetilde Y;\widetilde E)$, we now use it to define an algebra of smoothing operators on the end-periodic manifold $M$. Roughly speaking, an end-periodic smoothing operator consists of two parts: a periodic normal operator describing its periodic behavior sufficiently far along the end, and a residual operator whose Schwartz kernel decays rapidly along the end in both variables.

\begin{definition}[Residual end-periodic smoothing operators]
Choose a smooth extension $\rho:M\to[0,\infty)$ that agrees with $f$ sufficiently far along the end. We define the space
\[
\dt{\Psi}_{\mathrm{ep}}^{-\infty}(M;E) \subset \Psi^{-\infty}(M;E)
\]
by requiring an element $R \in \dt{\Psi}_{\mathrm{ep}}^{-\infty}(M;E)$ to satisfy
\begin{equation}\label{absolute-schwartz-decay}
\sup_{x,y\in M}\langle \rho(x)+ \rho(y)\rangle^{N}
\left|\nabla_x^{\alpha}\nabla_y^{\beta}
K_R(x,y)\right|<\infty
\end{equation}
for every $N,\alpha,\beta$. We refer to this as \emph{absolute Schwartz decay}, in contrast with the relative Schwartz decay in \eqref{eqn-sus-schwartz}.
\end{definition}

\begin{definition}[End-periodic smoothing operators] \ \vspace{1mm} \\
Choose a cutoff function $\chi\in C^{\infty}(M;[0,1])$ such that:
\begin{enumerate}
\item $\chi$ is supported in the periodic end;
\item $\chi(x)=1$ sufficiently far along the periodic end.
\end{enumerate}
We define the space of end-periodic smoothing operators $\Psi_{\mathrm{ep}}^{-\infty}(M;E) \subset \Psi^{-\infty}(M;E)$ as
\[
\Psi_{\mathrm{ep}}^{-\infty}(M;E):=\bigl\{
A \in \Psi^{-\infty}(M;E) : \exists\,\widetilde A
\in\Psi_{\mathbb{Z}}^{-\infty}(\widetilde Y;\widetilde E) \text{ s.t. }A-\chi\widetilde A\chi
\in\dt{\Psi}_{\mathrm{ep}}^{-\infty}(M;E)
\bigr\}.
\]
Here 
\[
\chi\widetilde A\chi\in\Psi^{-\infty}(M;E), \hspace{3mm}K_{\chi\widetilde A\chi}(x,y)=\chi(x)K_{\widetilde A}(x,y)\chi(y).
\]
This definition is independent of the choice of cutoff function $\chi$.
\end{definition}

\begin{lemma}\label{lem:zero-normal}
Let $\widetilde A \in \Psi_{\mathbb{Z}}^{-\infty} (\widetilde Y;\widetilde E)$. If $ \chi\widetilde A\chi \in \dt{\Psi}_{\mathrm{ep}}^{-\infty}(M;E)$, then $\widetilde A=0$.
\end{lemma}
\begin{proof}
This follows from the absolute Schwartz decay \eqref{absolute-schwartz-decay}, combined with the periodicity of $\widetilde A \in \Psi_{\mathbb{Z}}^{-\infty}$.
\end{proof}

\begin{lemma}
For any $A \in \Psi_{\mathrm{ep}}^{-\infty}(M;E)$, the associated operator $\widetilde{A} \in \Psi_{\mathbb{Z}}^{-\infty}(\widetilde Y;\widetilde E)$ is uniquely defined, so we can call it the normal operator $N(A):=\widetilde{A}$.
\end{lemma}
\begin{proof}
Suppose that we have two normal operators $\widetilde A_{1}$ and $\widetilde{A}_{2}$ with
\[
A-\chi\widetilde A_1\chi \in \dt{\Psi}_{\mathrm{ep}}^{-\infty}(M;E),
\]
\[
A-\chi\widetilde A_2\chi \in \dt{\Psi}_{\mathrm{ep}}^{-\infty}(M;E).
\]
Subtracting the two relations gives
\[
\chi(\widetilde A_1-\widetilde A_2)\chi \in \dt{\Psi}_{\mathrm{ep}}^{-\infty}(M;E).
\]
Applying the preceding lemma to $\widetilde A_1-\widetilde A_2$ shows that $\widetilde A_1=\widetilde A_2$.
\end{proof}

Next, we establish some basic composition properties for the end-periodic smoothing algebra. The proof is given in Appendix \ref{appendix-1}; see Lemma \ref{schwartz-convolution-lemma}.

\begin{lemma}\label{composition-lemma}
The following composition properties hold:

\begin{enumerate}
\item $\dt{\Psi}_{\mathrm{ep}}^{-\infty}(M;E)
\circ \dt{\Psi}_{\mathrm{ep}}^{-\infty}(M;E) \subset
\dt{\Psi}_{\mathrm{ep}}^{-\infty}(M;E);$
\item if $\widetilde A \in \Psi_{\mathbb{Z}}^{-\infty} (\widetilde Y;\widetilde E)$, and $R\in\dt{\Psi}_{\mathrm{ep}}^{-\infty}(M;E)$,
then
\[
(\chi\widetilde A\chi)R, \quad R(\chi\widetilde A\chi) \in \dt{\Psi}_{\mathrm{ep}}^{-\infty}(M;E);
\]

\item if $\widetilde A,\widetilde B \in \Psi_{\mathbb{Z}}^{-\infty}(\widetilde Y;\widetilde E)$, then
\[
(\chi\widetilde A\chi) (\chi\widetilde B\chi) -
\chi(\widetilde A\widetilde B)\chi \in \dt{\Psi}_{\mathrm{ep}}^{-\infty}(M;E).
\]
\end{enumerate}
Consequently: $\Psi_{\mathrm{ep}}^{-\infty}(M;E)$ is a $\ast$-algebra under composition and adjoint, and $\dt{\Psi}_{\mathrm{ep}}^{-\infty}(M;E)$ is a two-sided $\ast$-ideal in $\Psi_{\mathrm{ep}}^{-\infty}(M;E)$.
\end{lemma}

In the proof of Lemma \ref{composition-lemma} we use the fact that $1+\rho(x)$ on $W_{j}$ is uniformly comparable to $1+j$. We will adopt a shorthand to describe this sort of estimate more concisely. For $a,b \geq 0$, we write
\begin{equation}\label{asymp-notation}
a\asymp b
\end{equation}
if there exist constants $0<c\leq C<\infty$ such that
\[
cb\leq a\leq Cb.
\]
In particular, for $\langle t\rangle=(1+t^2)^{1/2}$, the elementary inequalities
\[
\langle s+t\rangle \leq (1+|t|)\langle s\rangle,
\hspace{3mm} \langle s\rangle \leq (1+|t|)\langle s+t\rangle
\]
imply that, whenever $|t|\leq D$, one has
\begin{equation}\label{eq:uniform-weight-comparability}
(1+D)^{-1}\langle s\rangle \leq \langle s+t\rangle
\leq (1+D)\langle s\rangle.
\end{equation}
Thus $\langle s+t\rangle\asymp\langle s\rangle$ holds uniformly for bounded $t$. We will use this observation repeatedly when developing estimates in the end-periodic setting. For example, given local lifts $\sigma_i:U_i\to\widetilde Y$, the functions
\[
d_{ij}(z,z')= f(\sigma_j(z'))-f(\sigma_i(z))
\]
are uniformly bounded on the chosen coordinate patches. Consequently,
\[
\langle n+d_{ij}(z,z')\rangle \asymp\langle n\rangle, \hspace{3mm} n\in\mathbb Z,
\]
uniformly in the patch indices and spatial variables. More generally, if $x=T^jx_0$ and $y=T^\ell y_0$, with $x_0,y_0$ in a fixed compact fundamental segment $W$, then
\[
f(y)-f(x) =(\ell-j)+f(y_0)-f(x_0),
\]
and hence
\begin{equation}\label{eq:uniform-weight-comparability-2}
\langle f(y)-f(x)\rangle \asymp\langle j-\ell\rangle.
\end{equation}
Similarly, for a nonnegative extension $\rho$ agreeing with $f$ sufficiently far along the positive end,
\begin{equation}\label{eq:uniform-weight-comparability-3}
\langle\rho(x)+\rho(y)\rangle \asymp\langle j+\ell\rangle, \hspace{3mm} x\in W_j,\quad y\in W_\ell,\quad j,\ell\geq0.
\end{equation}
These estimates allow us to replace the weights in the definitions of relative and absolute Schwartz decay by the corresponding discrete weights on fundamental segments.

\begin{lemma}\label{lem:normal-op-homom}
The normal operator $N:\Psi_{\mathrm{ep}}^{-\infty}(M;E) \to \Psi_{\mathbb{Z}}^{-\infty} (\widetilde Y;\widetilde E)$ is a $\ast$-algebra homomorphism:
\[
N(AB) = N(A)N(B), \hspace{3mm} N(A^{*}) = N(A)^{*}.
\]
\end{lemma}
\begin{proof}
Write
\[
A=\chi\widetilde A\chi+R, \hspace{3mm} B=\chi\widetilde B\chi+S,
\]
where $R,S\in\dt{\Psi}_{\mathrm{ep}}^{-\infty}(M;E)$. Then
\[
\begin{split}
AB-\chi(\widetilde A\widetilde B)\chi &=
\bigl[(\chi\widetilde A\chi)(\chi\widetilde B\chi)
- \chi(\widetilde A\widetilde B)\chi\bigr]
\\
&\quad + (\chi\widetilde A\chi)S +R(\chi\widetilde B\chi) + RS.
\end{split}
\]
Every term on the right belongs to $\dt{\Psi}_{\mathrm{ep}}^{-\infty}(M;E)$ by Lemma \ref{composition-lemma}. Hence $AB\in\Psi_{\mathrm{ep}}^{-\infty}(M;E)$ with $N(AB)=\widetilde A\widetilde B$. For the adjoint property, since $\chi$ is real-valued, we have $(\chi\widetilde A\chi)^{*} = \chi\widetilde A^{*}\chi$. The residual class is invariant under adjoint, because $K_{R^{*}}(x,y)=K_R(y,x)^{*}$, and therefore
\[
A^{*}-\chi\widetilde A^{*}\chi = R^{*} \in \dt{\Psi}_{\mathrm{ep}}^{-\infty}(M;E).
\]
\end{proof}

\begin{lemma}\label{ses-lemma-1}
We have a short exact sequence of $*$-algebras
\[
0 \to \dt{\Psi}_{\mathrm{ep}}^{-\infty}(M;E) \hookrightarrow \Psi_{\mathrm{ep}}^{-\infty}(M;E) \xrightarrow{\;N\;} \Psi_{\mathbb{Z}}^{-\infty} (\widetilde Y;\widetilde E) \to 0.
\]
\end{lemma}

\begin{proof}
Every residual operator belongs naturally to the end-periodic space by choosing the suspended model to be zero. This gives us an injective inclusion $\dt{\Psi}_{\mathrm{ep}}^{-\infty}(M;E) \hookrightarrow \Psi_{\mathrm{ep}}^{-\infty}(M;E)$. The preceding lemmas in this section show that the normal map is a well-defined algebra homomorphism. If $N(A)=0$ then
\[
A = A-\chi N(A)\chi \in \dt{\Psi}_{\mathrm{ep}}^{-\infty}(M;E),
\]
hence $\ker N \subset \dt{\Psi}_{\mathrm{ep}}^{-\infty}(M;E)$. The reverse inclusion follows because Lemma \ref{lem:zero-normal} implies that a residual operator has zero normal operator. Therefore $\ker N
= \dt{\Psi}_{\mathrm{ep}}^{-\infty}(M;E)$. Finally, $N$ is surjective. Given $\widetilde A \in \Psi_{\mathbb{Z}}^{-\infty}(\widetilde Y;\widetilde E)$, the operator $\chi\widetilde A\chi$ belongs to $\Psi_{\mathrm{ep}}^{-\infty}(M;E)$ and satisfies $N(\chi\widetilde A\chi) = \widetilde A$. This proves exactness.
\end{proof}

\subsection{Fr\'echet space topologies}\label{sec:frechet-topologies}

We define natural Fr\'echet space topologies on the various function spaces and smoothing operators defined above. The most important result in this section is the Fourier-Laplace isomorphism theorem (Theorem \ref{FL-transform-thm}). Once we have Fr\'echet space topologies on the normal algebra $\Psi_{\mathbb{Z}}^{-\infty}$ and the residual algebra $\dt{\Psi}_{\mathrm{ep}}^{-\infty}$, the full smoothing algebra $\Psi_{\mathrm{ep}}^{-\infty}$ automatically obtains a natural Fr\'echet topology using the the exact normal sequence of Lemma \ref{ses-lemma-1}.

\begin{lemma}[Fr\'echet topology I]\label{lem:Frechet-topology-Z-sus}
The seminorms $\nu_{N,r}$ of \eqref{eqn-sus-schwartz} make $\Psi_{\mathbb Z}^{-\infty}(\widetilde Y;\widetilde E)$ a Fr\'echet space. The topology is independent of the function $f$, and independent of the choice of periodic Riemannian metric, Hermitian metric, and connection.
\end{lemma}

\begin{proof}
The deck group $\mathbb{Z}$ acts on $\widetilde{Y}\times\widetilde{Y}$ by $k\cdot(x,y) = (T^{k}x,T^{k}y)$, and the quotient $\mathcal{Y} := (\widetilde{Y}\times\widetilde{Y})/\mathbb{Z}$ is a smooth manifold. The function $d(x,y) = f(y)-f(x)$ descends to a proper, continuous function on $\mathcal{Y}$, and a smoothing kernel for an operator in the suspended algebra is equivalent to a smooth section of the vector bundle 
\[
\mathcal{V} := (\widetilde{E}\boxtimes\widetilde{E})/\mathbb{Z} \to \mathcal{Y}.
\]
The maps $\nu_{N,r}$ are seminorms, and they define a Hausdorff locally convex metrizable topology on the suspended algebra. It remains to verify completeness. 

Let $(K_m)_{m\geq1}$ be a Cauchy sequence with respect to the seminorms $\nu_{N,r}$. Thus, for every $N,r\geq0$ and every $\varepsilon>0$, there exists $m_0$ such that
\begin{equation}\label{eq:Z-sus-Cauchy-estimate}
\sup_{z\in\mathcal Y} \, \langle\delta(z)\rangle^N
|\nabla^r(K_m-K_n)(z)| <\varepsilon
\end{equation}
whenever $m,n\geq m_0$. Take a relatively compact coordinate neighborhood $U \subset \mathcal{Y}$ trivializing $\mathcal{V}$, and write the connections in a local frame. A straightforward local computation shows that one can find a global smooth section $K \in C^{\infty}(\mathcal{Y},\mathcal{V})$ such that $\nabla^{r}K_{m} \to \nabla^{r}K$ uniformly on compact subsets of $\mathcal{Y}$, for every $r \geq 0$. We must prove convergence in the weighted seminorms, which is stronger than convergence in $C^\infty$ on compact subsets.

Fix $N,r\geq0$ and $\varepsilon>0$.
Choose $m_0$ so that
\eqref{eq:Z-sus-Cauchy-estimate} holds for
all $m,n\geq m_0$. For fixed $m\geq m_0$, let $n\to\infty$.
The local smooth convergence established above implies pointwise convergence of $\nabla^rK_n$ to $\nabla^rK$. Therefore,
\[
\langle\delta(z)\rangle^N |\nabla^r(K_m-K)(z)| \leq\varepsilon
\]
for every $z\in\mathcal Y$. Taking the supremum over $z$ gives
\[
\nu_{N,r}(K_m-K)\leq\varepsilon, \hspace{3mm} m\geq m_0.
\]
Hence $\nu_{N,r}(K_m-K)\to$ for every $N,r\geq0$. Now fixing $m\geq m_0$, we have
\[
\nu_{N,r}(K) \leq \nu_{N,r}(K_m)+\nu_{N,r}(K-K_m) <\infty.
\]
Thus the limiting section satisfies the required relative Schwartz estimates. Lifting $K$ to $\widetilde Y\times\widetilde Y$ produces a smooth diagonally periodic kernel with the estimates \eqref{eqn-sus-schwartz}. Therefore the space is complete, and $\Psi_{\mathbb Z}^{-\infty}(\widetilde Y;\widetilde E)$ is a Fr\'echet space.

If \(f_0,f_1\) are two functions on $\widetilde{Y}$ satisfying \eqref{f-definition}, then \(f_1-f_0\) is deck-invariant and thus bounded, since it descends to the compact quotient $Y$. The estimate \eqref{eq:uniform-weight-comparability-2} therefore yields equivalence of the corresponding seminorms. Likewise, any two choices of periodic metric or connection are uniformly equivalent on $\widetilde{Y}$, because they descend to $Y$. The resulting kernel seminorms are therefore equivalent.
\end{proof}

For a smoothing operator $S$ on the closed manifold $S$, let
\begin{equation}\label{eq:FL-smoothing-seminorms}
p_r(S):=\max_{0\le j\le r}\sup_{(z,z')\in Y^2}
|\nabla^{j}K_{S}(z,z')|.
\end{equation}
These are the standard Fr\'echet seminorms on
$\Psi^{-\infty}(Y;E_Y)$. For a twisted family $A_{(\cdot)}$ as in \eqref{eqn-sus-twisted-algebra}, define
\begin{equation}\label{eq:FL-twisted-seminorms}
P_{r,q}(A):=\sup_{0\le\theta\le2\pi}p_r(\partial_\theta^q A_\theta), \hspace{3mm} r,q\ge0.
\end{equation}

\begin{lemma}[Fr\'echet topology II]\label{lem:Frechet-topology-twisted}
The seminorms $P_{r,q}$ make $C^\infty_{\mathrm{FL}}(\mathbb R;\Psi^{-\infty}(Y;E_Y))$ a Fr\'echet $*$-algebra under pointwise composition and adjoint.
\end{lemma}
\begin{proof}
Conjugation by $C$ is a continuous $*$-automorphism on $\Psi^{-\infty}(Y;E_{Y})$. Hence the twisting condition $A_{\theta+2\pi}=CA_\theta C^{-1}$ defines a closed subspace of $C^\infty(\mathbb R;\Psi^{-\infty})$, and the seminorms $P_{r,q}$ induce its usual smooth-family topology. Continuity of pointwise composition and adjoint follows from the standard continuous $*$-algebra operations on smoothing kernels over compact $Y$, together with the Leibniz rule in $\theta$. Both operations preserve the twisting condition.
\end{proof}

\begin{thm}\label{FL-transform-thm}
Conjugation by the FL transform defines a $*$-algebra isomorphism
\begin{equation}\label{eqn-sus-FL-algebra-isomorphism}
\Psi^{-\infty}_{\mathbb Z} (\widetilde Y;\widetilde E) \xrightarrow{\ \cong\ } C^\infty_{\mathrm{FL}}(\mathbb{R};\Psi^{-\infty}(Y;E_Y)), \hspace{3mm} \widetilde{A} \mapsto \mathcal{F}\widetilde{A}\mathcal{F}^{-1}.
\end{equation}
Here the image is the family $A_{\theta}$ given by \eqref{eqn-sus-FL-intertwining}. The inverse of \eqref{eqn-sus-FL-algebra-isomorphism} is given on Schwartz kernels by \eqref{eqn-kernel}. Moreover, \eqref{eqn-sus-FL-algebra-isomorphism} and its inverse are continuous, so it is an isomorphism of Fr\'echet spaces.
\end{thm}

\begin{proof}
Choose finitely many open sets $U_i\subset Y$, smooth lifts $\sigma_i:U_i\to\widetilde Y$, and relatively compact subsets $V_i\Subset U_i$ whose closures $F_i=\overline{V_i}$ cover $Y$. For $(z,z')\in U_i\times U_j$, define the relative height
\[
d_{ij}(z,z')=f(\sigma_j(z'))-f(\sigma_i(z)).
\]
All derivatives of $d_{ij}$ are bounded on $F_i\times F_j$, and, for a constant $D$ independent of $i,j$, we have $|d_{ij}|\le D$ there. In particular, $\langle n+d_{ij}\rangle\asymp\langle n\rangle$ uniformly for $n\in\mathbb Z$ and $(z,z')\in F_i\times F_j$, cf. \eqref{asymp-notation}.

\medskip
\noindent\emph{The forward kernel formula.}
Let $\widetilde A\in\Psi_{\mathbb Z}^{-\infty}
(\widetilde Y;\widetilde E)$. For $z,z'\in Y$, choose lifts $x,y\in\widetilde Y$ and define
\begin{equation}\label{eqn-sus-forward-kernel}
K_{A_\theta}(z,z') =e^{-i\theta(f(y)-f(x))}
\sum_{n\in\mathbb Z}e^{-in\theta} K_{\widetilde A}(x,T^ny).
\end{equation}
This formula is independent of the chosen lifts.

\medskip
\noindent\emph{The twisted family.}
Conversely, for a twisted family $A_{(\cdot)}\in C^\infty_{\mathrm{FL}}(\mathbb R; \Psi^{-\infty}(Y;E_Y))$, the inverse kernel is
\eqref{eqn-kernel}. Put
\begin{equation}\label{eqn-sus-Btheta}
B(\theta;x,y) :=e^{i\theta(f(y)-f(x))}K_{A_\theta}(p(x),p(y)).
\end{equation}
The twisting relation $A_{\theta+2\pi}=CA_\theta C^{-1}$ implies that
\begin{align}
B(\theta+2\pi;x,y) &=e^{i(\theta+2\pi)(f(y)-f(x))}e^{2\pi i(f(x)-f(y))} K_{A_\theta}(p(x),p(y))\notag\\
&=B(\theta;x,y). \label{eqn-sus-B-periodicity}
\end{align}
Thus $B$ is $2\pi$-periodic in $\theta$, and the same goes for its derivatives. In particular, integration by parts in $\theta$ produces no boundary terms.

\medskip
\noindent\emph{The forward estimate.}
On $F_i\times F_j$ let
\[
H_{ij,n}(z,z') :=K_{\widetilde A}(\sigma_i(z),T^n\sigma_j(z')).
\]
The relative Schwartz estimates give, for every $M,s\ge0$,
\[
|\nabla^s H_{ij,n}(z,z')| \le \nu_{M,s}(\widetilde A) \langle n+d_{ij}(z,z')\rangle^{-M}.
\]
In these coordinates, \eqref{eqn-sus-forward-kernel} becomes
\[
K_{A_\theta}(z,z') =\sum_{n\in\mathbb Z}e^{-i\theta(n+d_{ij}(z,z'))}H_{ij,n}(z,z').
\]
Fix $r,q\ge0$. Since all derivatives of $d_{ij}$ are bounded on $F_i\times F_j$, the Leibniz rule yields, for $0\le k\le r$ and $0\le\theta\le2\pi$,
\[
\left|\nabla^k\partial_\theta^q
\bigl(e^{-i\theta(n+d_{ij})}H_{ij,n}\bigr)\right|
\le C_{i,j,r,q}\langle n+d_{ij}\rangle^q
\sum_{s=0}^k|\nabla^sH_{ij,n}|.
\]
Taking $M=q+2$ and using $\sup_{u\in\mathbb R}\sum_{n\in\mathbb Z}\langle n+u\rangle^{-2}<\infty$, we obtain
\begin{equation}\label{eq:FL-kernel-forward-estimate}
\max_{0\le k\le r} \sup_{\substack{0\le\theta\le2\pi\\(z,z')\in Y^2}} |\nabla^k\partial_\theta^qK_{A_\theta}(z,z')| \le C_{r,q}\sum_{s=0}^r\nu_{q+2,s}(\widetilde A).
\end{equation}
The same argument applies on any bounded interval of $\theta$. Thus the series and all its derivatives converge uniformly, and $A_\theta$ is a smooth family of smoothing operators. Moreover, \eqref{eqn-sus-forward-kernel} gives
\[
K_{A_{\theta+2\pi}}(z,z') =e^{2\pi i(f(x)-f(y))}K_{A_\theta}(z,z') =K_{CA_\theta C^{-1}}(z,z').
\]
Hence $A_{(\cdot)}\in
C^\infty_{\mathrm{FL}}(\mathbb R;\Psi^{-\infty}(Y;E_Y))$.

\medskip
\noindent\emph{The inverse estimate.}
Conversely, let
$A_{(\cdot)}\in C^\infty_{\mathrm{FL}}(\mathbb R;\Psi^{-\infty}(Y;E_Y))$ and define $K_{\widetilde A}$ by \eqref{eqn-kernel}. For $(z,z')\in F_i\times F_j$, put
\[
B_{ij}(\theta;z,z') :=B(\theta;\sigma_i(z),\sigma_j(z')) =e^{i\theta d_{ij}(z,z')}K_{A_\theta}(z,z').
\]
Writing $x=\sigma_i(z)$ and $y=T^n\sigma_j(z')$ in
\eqref{eqn-kernel}, we obtain
\begin{equation}\label{eq:FL-local-inverse-coefficients}
K_{\widetilde A}(\sigma_i(z),T^n\sigma_j(z'))
=\frac1{2\pi}\int_0^{2\pi} e^{in\theta}B_{ij}(\theta;z,z')\,d\theta.
\end{equation}
By the Leibniz rule and the bounds on $d_{ij}$, for every $L,r\ge0$ and $0\le k\le r$,
\[
\sup_{\substack{0\le\theta\le2\pi\\(z,z')\in F_i\times F_j}} |\nabla^k\partial_\theta^L B_{ij}(\theta;z,z')| \le C_{i,j,r,L}\sum_{q=0}^L \max_{0\le s\le r} \sup_{\substack{0\le\theta\le2\pi\\(z,z')\in Y^2}}
|\nabla^s\partial_\theta^q K_{A_\theta}(z,z')|.
\]
For $n\ne0$, integrate by parts $L$ times in
\eqref{eq:FL-local-inverse-coefficients} and use
\eqref{eqn-sus-B-periodicity}; for $n=0$, use the integral directly. Combining the two cases gives
\begin{align*}
& \sup_{(z,z')\in F_i\times F_j}|\nabla^k K_{\widetilde A} (\sigma_i(z),T^n\sigma_j(z'))| \\
& \hspace{6mm} \le C_{i,j,r,L}(1+|n|)^{-L} \sum_{q=0}^L\max_{0\le s\le r} \sup_{\substack{0\le\theta\le2\pi\\(z,z')\in Y^2}}
|\nabla^s\partial_\theta^qK_{A_\theta}(z,z')|.
\end{align*}
The kernel \eqref{eqn-kernel} is smooth and invariant under deck translation. Every pair of points in $\widetilde Y^2$ is obtained by translating a pair $(\sigma_i(z),T^n\sigma_j(z'))$, and its relative height is $n+d_{ij}(z,z')$. Taking $L=N+1$ and using $\langle n+d_{ij}\rangle\asymp\langle n\rangle$, we obtain
\begin{equation}\label{eq:FL-kernel-inverse-estimate}
\nu_{N,r}(\widetilde A) \le C_{N,r}\sum_{q=0}^{N+1}
\max_{0\le s\le r} \sup_{\substack{0\le\theta\le2\pi\\(z,z')\in Y^2}}|\nabla^s\partial_\theta^qK_{A_\theta}(z,z')|.
\end{equation}
Thus $\widetilde A\in\Psi_{\mathbb Z}^{-\infty}
(\widetilde Y;\widetilde E)$. It is clear that the two kernel constructions are inverse to one another.

\medskip
\noindent \emph{Continuity.} Expressed in terms of the seminorms, the estimates \eqref{eq:FL-kernel-forward-estimate} and \eqref{eq:FL-kernel-inverse-estimate} become
\[
\begin{aligned}
P_{r,q}(\mathcal F\widetilde A\mathcal F^{-1})
&\le C_{r,q}\sum_{s=0}^r\nu_{q+2,s}(\widetilde A),\\
\nu_{N,r}(\mathcal F^{-1}A_{(\cdot)}\mathcal F)
&\le C_{N,r}\sum_{q=0}^{N+1}P_{r,q}(A).
\end{aligned}
\]
Thus the algebraic isomorphism and its inverse are continuous.
\end{proof}

For $R\in\dt{\Psi}_{\mathrm{ep}}^{-\infty}(M;E)$, put
\begin{equation}\label{eq:residual-frechet-seminorms}
\mu_{N,r}(R):=\sum_{a+b\le r}\sup_{x,y\in M}
\langle\rho(x)+\rho(y)\rangle^N |\nabla_x^a\nabla_y^bK_R(x,y)|, \hspace{3mm} N,r\ge0.
\end{equation}
These seminorms make the residual algebra a Fr\'echet space. Indeed, they are the weighted $C^\infty$ seminorms on $M\times M$, and the same local $C^\infty$-limit and uniform weighted-Cauchy argument as in Lemma \ref{lem:Frechet-topology-Z-sus} proves completeness.

Fix an end cutoff $\chi$ as in the definition of
$\Psi_{\mathrm{ep}}^{-\infty}(M;E)$, and define
\begin{equation}\label{eq:Ext-chi-topology}
\operatorname{Ext}_\chi: \Psi_{\mathbb Z}^{-\infty}(\widetilde Y;\widetilde E) \to\Psi_{\mathrm{ep}}^{-\infty}(M;E), \hspace{3mm} \widetilde A\mapsto\chi\widetilde A\chi.
\end{equation}
By Lemma \ref{ses-lemma-1}, $N\circ\operatorname{Ext}_\chi=\mathrm{Id}$
and $\ker N=\dt{\Psi}_{\mathrm{ep}}^{-\infty}(M;E)$.
Consequently, putting
\begin{equation}\label{eq:residual-component}
R_\chi(A):=A-\operatorname{Ext}_\chi(N(A)),
\end{equation}
we obtain the linear bijection
\begin{equation}\label{eq:ep-direct-sum-map}
\begin{aligned}
\Phi_\chi:\Psi_{\mathrm{ep}}^{-\infty}(M;E)
&\to \Psi_{\mathbb Z}^{-\infty}(\widetilde Y;\widetilde E) \oplus\dt{\Psi}_{\mathrm{ep}}^{-\infty}(M;E),\\
A&\mapsto(N(A),R_\chi(A)),
\end{aligned}
\end{equation}
whose inverse is
\begin{equation}\label{eq:ep-direct-sum-inverse}
\Phi_\chi^{-1}(\widetilde A,R) =\operatorname{Ext}_\chi(\widetilde A)+R.
\end{equation}
This leads to a natural Fr\'echet topology on $\Psi_{\mathrm{ep}}^{-\infty}(M;E)$.

\begin{lemma}[Fr\'echet topology III]\label{lem:Frechet-topology-ep}
The seminorms
\begin{equation}\label{eq:ep-full-frechet-seminorms}
p_{N,r}^\chi(A):=\nu_{N,r}(N(A))+\mu_{N,r}(R_\chi(A)), \hspace{3mm} N,r\ge0,
\end{equation}
make $\Psi_{\mathrm{ep}}^{-\infty}(M;E)$ a Fr\'echet $*$-algebra. The map $\Phi_\chi$ is a Fr\'echet-space isomorphism, the residual algebra is a closed two-sided $*$-ideal, and $N$ is a continuous surjective $*$-homomorphism. The topology is independent of the cutoff, metric, and connection.
\end{lemma}
\begin{proof}
Write $\mathcal N=\Psi_{\mathbb Z}^{-\infty}(\widetilde Y;\widetilde E)$,
$\mathcal R=\dt{\Psi}_{\mathrm{ep}}^{-\infty}(M;E)$, and $\mathcal A=\Psi_{\mathrm{ep}}^{-\infty}(M;E)$, and set $\bar\nu_{N,r}=\sum_{s=0}^r\nu_{N,s}$. By construction, the map $\Phi_{\chi}$ transports the product Fr\'echet topology on $\mathcal{N}\oplus\mathcal{R}$ to the topology on $\mathcal{A}$ deteremined by the seminorms \eqref{eq:ep-full-frechet-seminorms}. Thus $\mathcal A$ is Fr\'echet and $\Phi_\chi$ is a topological isomorphism. In these coordinates $N$ is the first projection and $\mathcal R=\ker N$, so $N$ is continuous and $\mathcal R$ is closed.

To check cutoff independence, let $\chi_0,\chi_1$ be admissible, put $q=\chi_1-\chi_0$, and define
\begin{equation}\label{eq:ep-cutoff-difference-proof}
D_{10}(P):=\operatorname{Ext}_{\chi_1}(P)
-\operatorname{Ext}_{\chi_0}(P) = qP\chi_1+\chi_0Pq.
\end{equation}
The function $q$ is compactly supported in the end. Whenever $x\in\operatorname{supp}q$ and
$y\in\operatorname{supp}\chi_i$, one has
\[
\langle\rho(x)+\rho(y)\rangle \le C\langle f(y)-f(x)\rangle.
\]
The same inequality holds with $x,y$ interchanged. The cutoffs have bounded derivatives, so the Leibniz rule and the relative Schwartz estimates give
\begin{equation}\label{eq:ep-cutoff-continuity-proof}
\mu_{N,r}(D_{10}(P))\le C_{N,r}\bar\nu_{N,r}(P).
\end{equation}
In particular, $D_{10}:\mathcal N\to\mathcal R$ is continuous. Hence
\[
\Phi_{\chi_1}\Phi_{\chi_0}^{-1}(P,R)=(P,R-D_{10}(P)),
\]
is a Fr\'echet isomorphism, and the topologies induced by $\chi_{0}$ and $\chi_{1}$ coincide.

For two end-periodic metrics and connections, the fiber norms and covariant derivatives are uniformly comparable, since the data descend to the closed manifold $Y$ outside a compact part of $M$. Thus the seminorms $\mu_{N,r}$ are equivalent for these choices.

Finally, it is straightforward to verify continuity of the algebra operations using the estimates in Lemma \ref{schwartz-convolution-lemma}.
\end{proof}

\subsection{Mapping properties}

Having established the algebraic structure and natural Fr\'echet topologies of the end-periodic smoothing spaces, we now turn to their analytic properties as operators on $L^2(M;E)$. Residual smoothing operators are trace-class, while general end-periodic smoothing operators are bounded (but, in general, not compact) on $L^2(M;E)$. These properties will be used in the Fredholm theory developed below, and in particular will allow us to apply the ordinary operator trace to residual operators.

\begin{lemma}[Residual operators are trace-class]\label{lem:ep-residual-trace-class}
Every $R \in \dt{\Psi}_{\mathrm{ep}}^{-\infty}(M;E)$ is trace-class (in particular, compact) on $L^{2}(M;E)$. Moreover, the inclusion
\[
\dt{\Psi}_{\mathrm{ep}}^{-\infty}(M;E)\hookrightarrow \mathcal L^{1}(L^2(M;E))
\]
is continuous for the residual Fr\'echet topology. Here we write $\mathcal L^1(L^2(M;E))$ for the trace-class ideal with norm $\|\cdot\|_1$.
\end{lemma}


\begin{proof}
Put $d=\dim M$ and fix an integer $L>d/4$. We first recall the estimate underlying the trace-class criterion of \cite[Corollary 4.4]{dr2014} on compact manifolds. For a smooth matrix-valued kernel $F$ on $\mathbb T^d\times\mathbb T^d$ (here $\mathbb{T}^{d} = (\mathbb{R}/\mathbb{Z})^{d}$), let $T_F$ be its integral operator. Factoring
\[
T_F=(1-\Delta)^{-L}\bigl((1-\Delta)^L T_F\bigr)
\]
and using the product estimate for Hilbert-Schmidt operators gives
\begin{equation}\label{eq:ep-residual-local-trace}
\|T_F\|_1 \leq \|(1-\Delta)^{-L}\|_2\,
\|(1-\Delta_u)^L F\|_{L^2} \leq C_L\max_{|\alpha|\leq 2L} \|\partial_u^\alpha F\|_{L^\infty}.
\end{equation}
Here $(1-\Delta)^{-L}$ is Hilbert-Schmidt because
$\sum_{m\in\mathbb Z^d}(1+|m|^2)^{-2L}<\infty$.

Choose a uniformly locally finite coordinate cover
$\{U_i\}_{i\in\mathcal I}$ of $M$, with a subordinate smooth partition of unity $\{\phi_i\}$, by translating finitely many charts and bundle frames along the periodic end and adding finitely many charts on the compact core. Each chart has a level $\ell(i)\in\mathbb Z_{\geq0}$; there are uniformly many charts per level, and \eqref{eq:uniform-weight-comparability-3} gives
\[
c(1+\ell(i)) \leq 1+\rho(x) \leq C(1+\ell(i)),
\hspace{3mm} x\in\operatorname{supp}\phi_i,
\]
where $c,C>0$ are constants independent of $i$. We may take all coordinate images to be $(-\pi,\pi)^d$ and the cutoffs to have support in their interiors.

For $R\in\dt{\Psi}_{\mathrm{ep}}^{-\infty}(M;E)$, put $R_{ij}=\phi_i R\phi_j$. In local unitary frames, transport its kernel to $(-\pi,\pi)^d\times(-\pi,\pi)^d$, incorporating the square roots of the coordinate volume densities, and extend it by zero to a smooth matrix-valued kernel $F_{ij}$ on $\mathbb T^d\times\mathbb T^d$.
The absolute Schwartz estimate \eqref{absolute-schwartz-decay} and the uniform bounds on the local data yield
\[
\max_{|\alpha|\leq2L} \|\partial_u^\alpha F_{ij}\|_{L^\infty} \leq C_L\,(1+\ell(i)+\ell(j))^{-4}\,\mu_{4,2L}(R),
\]
where $\mu_{N,r}$ is defined in
\eqref{eq:residual-frechet-seminorms}. Since $R_{ij}$ is unitarily equivalent to a compression of the torus operator with kernel $F_{ij}$,
\eqref{eq:ep-residual-local-trace} implies
\begin{equation}\label{eq:ep-residual-block-trace}
\|R_{ij}\|_1 \leq C_L\,(1+\ell(i)+\ell(j))^{-4}\,\mu_{4,2L}(R).
\end{equation}
Consequently, the uniform bound on the number of charts at each level gives
\[
\sum_{i,j\in\mathcal I}\|R_{ij}\|_1 \leq C_L\,\mu_{4,2L}(R) \sum_{a,b\geq0}(1+a+b)^{-4}<\infty.
\]
Thus $S:=\sum_{i,j}R_{ij}$ converges in trace norm. Moreover, $S$ is the $L^2$-extension of $R$, and
\[
\|R\|_1\leq C_L\,\mu_{4,2L}(R).
\]
This proves both the trace-class assertion for $R$ and continuity of the inclusion for the residual Fr\'echet topology.
\end{proof}

\begin{lemma}[Boundedness]
\label{ep-theory-lemma-1}
Every operator $A\in\Psi_{\mathrm{ep}}^{-\infty}(M;E)$ extends to a bounded operator on $L^{2}(M;E)$, and the inclusion
\[
\Psi_{\mathrm{ep}}^{-\infty}(M;E)\hookrightarrow
\mathcal B(L^2(M;E))
\]
is continuous for the Fr\'echet topology defined in Lemma \ref{lem:Frechet-topology-ep}.
\end{lemma}

\begin{proof}
Fix an end cutoff $\chi$ as in
\eqref{eq:Ext-chi-topology} and write
\[
A=\operatorname{Ext}_\chi(\widetilde A)+R =\chi\widetilde A\chi+R, \hspace{3mm} \widetilde A\in\Psi_{\mathbb Z}^{-\infty} (\widetilde Y;\widetilde E), \hspace{3mm} R \in\dt{\Psi}_{\mathrm{ep}}^{-\infty}(M;E).
\]
Since $R$ is trace-class on $L^{2}(M;E)$, it is bounded, so it suffices to treat the term $\chi\widetilde{A}\chi$. By Theorem \ref{FL-transform-thm}, $\widetilde A$ has an indicial family
\[
A_{(\cdot)}\in C^\infty_{\mathrm{FL}} (\mathbb R;\Psi^{-\infty}(Y;E_Y)), \hspace{3mm} A_{\theta+2\pi}=CA_\theta C^{-1}.
\]
Set $H_Y=L^2(Y;E_Y)$. Since $Y$ is compact, the kernel of each $A_\theta$ is square-integrable on $Y\times Y$. The twisted-family seminorm \eqref{eq:FL-twisted-seminorms} has
\[
P_{0,0}(A_{(\cdot)}) = \sup_{0\leq\theta\leq2\pi}p_0(A_\theta)<\infty.
\]
Therefore
\begin{equation}\label{eq:ep-bounded-normal-norm}
\begin{aligned}
\|\widetilde A\|_{\mathcal B(L^2(\widetilde Y;\widetilde E))} &=\operatorname*{ess\,sup}_{0\leq\theta<2\pi}\|A_\theta\|_{\mathcal B(H_Y)}\\
&=\sup_{0\leq\theta\leq2\pi} \|A_\theta\|_{\mathcal B(H_Y)}\leq c \operatorname{vol}(Y)P_{0,0}(A_{(\cdot)}).
\end{aligned}
\end{equation}
Here $c$ is a constant. Let $J_\chi:L^2(M;E)\to L^2(\widetilde Y;\widetilde E)$ be the operator obtained by multiplying by $\chi$ on the identified end and extending the result by zero to the full cyclic cover. It is bounded with $\|J_\chi\|\leq\|\chi\|_\infty$. Since $\chi$ is real-valued we have $\operatorname{Ext}_\chi(\widetilde A)=J_\chi^*\widetilde A J_\chi$, hence
\begin{equation}\label{eq:ep-bounded-cutoff-estimate}
\|\operatorname{Ext}_\chi(\widetilde A)\|_{\mathcal B(L^2(M;E))} \leq\|\chi\|_\infty^2 c_E\operatorname{vol}(Y) P_{0,0}(A_{(\cdot)}).
\end{equation}
We conclude that $A$ defines a bounded operator on $L^{2}(M;E)$.
\end{proof}

\subsection{Renormalized trace}

For any $A \in \Psi_{\mathrm{ep}}^{-\infty}(M;E)$, fix an end cutoff $\chi$ as in \eqref{eq:Ext-chi-topology}, and write
\[
A = \chi\widetilde{A}\chi + R, \hspace{3mm} \widetilde{A} \in \Psi_{\mathbb{Z}}^{-\infty}(M;E), \hspace{3mm} R \in \dt{\Psi}_{\mathrm{ep}}(M;E).
\]
By Theorem \ref{FL-transform-thm}, there is a uniquely associated family of operators
\[
A_{(\cdot)} \in C^{\infty}_{\mathrm{FL}}(\mathbb{R};\Psi^{-\infty}(Y;E_{Y})), \hspace{3mm} A_{(\cdot)} := \mc{F}\widetilde{A}\mc{F}^{-1}.
\]
We call $A_{(\cdot)}$ the (Fourier-Laplace) \emph{indicial family} associated with $A$.

We adopt the following convention. Although the indicial family is technically indexed by $\theta \in \mathbb{R}$, satisfying $A_{\theta+2\pi}=CA_{\theta}C^{-1}$, we will often write the indicial family as $A_{\xi}$ with $\xi = e^{i\theta} \in S^{1}$. Composition, adjoints, and inversion are performed fiberwise, whereas traces, Fredholm determinants, and invertibility are independent of the choice of real lift because the transition acts by unitary conjugation. As a result, this convention is usually harmless.

There is a nice formula relating the smooth integrals kernels of $\widetilde{A}$ and $A_{\xi}$:
\[
K_{\widetilde{A}}(y,y') = \frac{1}{2\pi i}\oint_{S^{1}} \xi^{f(y')-f(y)} K_{A_{\xi}}(p(y),p(y'))\frac{d\xi}{\xi}.
\]
See \cite[Proposition 2.5]{mrs1}.

\medskip

Let $W\subset\widetilde Y$ be the compact fundamental segment for the periodic end, and write
\[
M_N = Z\cup W_0\cup\cdots\cup W_N, \hspace{3mm} W_{j} \cong W, \hspace{3mm} N \in \mathbb{Z}_{\geq 0}.
\]
There is a natural renormalized trace in the end-periodic setting, cf. \cite[eqn. (14)]{mrs1}, defined by
\begin{equation}\label{eqn-renormalized-trace}
\rTr(A) := \lim_{N\to\infty}\left[\int_{M_N} \mathrm{tr}(K_{A}(x,x))\,dx- (N+1)\int_W \tr(K_{\widetilde{A}}(x,x))\,dx\right].
\end{equation}
In \cite[Lemma 3.1]{mrs1} it is shown that this limit converges whenever $A$ is of the form $A = D^{m}e^{-tD^{2}}$, where $D$ is an end-periodic Dirac operator. This is a special subclass of end-periodic smoothing operators that suffices for the purposes of \cite{mrs1}.

Below, we extend the renormalized trace and the trace defect formula given in \cite[Theorem 4.6]{mrs1} to the algebra $\Psi_{\mathrm{ep}}^{-\infty}(M;E)$. The point is that $\rTr$ is not a trace in the usual sense because it need not vanish on commutators, and the trace defect formula expresses this defect solely in terms of the \emph{indicial families} of the involved operators.

\begin{prop}\label{renormalized-trace-lemma}
For every $A\in\Psi_{\mathrm{ep}}^{-\infty}(M;E)$, the limit \eqref{eqn-renormalized-trace} exists in $\mathbb C$. It defines a continuous linear functional
\[
\rTr:\Psi_{\mathrm{ep}}^{-\infty}(M;E)\to\mathbb C.
\]
Moreover, for any $A,B\in\Psi_{\mathrm{ep}}^{-\infty}(M;E)$, one has
\begin{align*}
\rTr[A,B] &=-\frac{1}{2\pi i}\oint_{S^1}\Tr \left(\frac{\partial A_\xi}{\partial\xi}B_\xi\right) \,d\xi\\
&\quad+\frac{1}{2\pi i}\oint_{S^1}\int_W f(x)\,\tr (K_{[A_\xi,B_\xi]}(x,x))\,dx\,\frac{d\xi}{\xi}, \numberthis \label{eq:trace-defect}
\end{align*}
with Fourier-Laplace indicial families $A_\xi,B_\xi$.
\end{prop}

\begin{proof}
We first prove existence of the limit by adapting the argument of \cite[Lemma 3.1]{mrs1}, using the absolute Schwartz estimates \eqref{absolute-schwartz-decay} in place of their heat kernel estimates.

Choose $J\geq1$ so that $\chi=1$ on $W_j=T^jW$. By deck-periodicity \eqref{eqn-sus-periodicity}, the integral of $\tr(K_{\widetilde{A}}(x,x))$ gives the same value over each $W_{j}$. On the other hand, \eqref{absolute-schwartz-decay} and the comparison \eqref{eq:uniform-weight-comparability-3} on $W_j$ give, for every integer $s\geq0$,
\[
\int_{W_j}|\tr K_R(x,x)|\,dx \leq C_s\,\mu_{s,0}(R)(1+j)^{-s}.
\]
Thus the diagonal integral of $R$ is absolutely convergent on $M$. Since $R$ is trace-class by Lemma \ref{lem:ep-residual-trace-class}, its integral equals $\Tr(R)$. For $N\geq J$, periodicity and the support of $\chi$ in the end imply
\begin{align*}
&\int_{M_N}\tr K_A(x,x)\,dx-(N+1)\int_W\tr K_{\widetilde{A}}(x,x)\,dx \\
&\hspace{3mm}=\int_{M_N}\tr K_R(x,x)\,dx + \sum_{j=0}^{J-1}\int_{W_j}(\chi(x)^2-1)\tr K_{\widetilde{A}}(x,x)\,dx.
\end{align*}
Passing to the limit yields
\begin{equation}\label{eq:ren-tr-splitting-formula}
\rTr(A)=\Tr(R)+\sum_{j=0}^{J-1}\int_{W_j} (\chi(x)^2-1)\tr K_{\widetilde{A}}(x,x)\,dx.
\end{equation}
This proves existence and linearity. The continuity
assertion follows from continuity of the normal operator and Lemma \ref{lem:ep-residual-trace-class}. 


Now consider the trace defect. The relevant limits will be justified using the relative Schwartz decay \eqref{eqn-sus-schwartz}. Write
\[
A=A^\chi+R_A, \quad B=B^\chi+R_B,\quad
A^\chi=\chi \widetilde{A}\chi,\quad B^\chi=\chi \widetilde{B}\chi.
\]
Consider the commutator $[A,B]$. By Lemma \ref{composition-lemma}, the commutator terms involving one of $R_A,R_B$ are residual, hence trace-class by Lemma \ref{lem:ep-residual-trace-class}. Consequently $\rTr[A,B] = \rTr[A^{\chi},B^{\chi}]$. Using the inverse kernel formula \eqref{eqn-kernel} we have
\begin{equation}\label{eq:ren-tr-normal-comm-zero}
\int_W\tr K_{[\widetilde{A},\widetilde{B}]}(x,x)\,dx =\frac{1}{2\pi}\int_0^{2\pi} \Tr([A_\theta,B_\theta])\,d\theta=0.
\end{equation}
Hence the subtraction term in \eqref{eqn-renormalized-trace} is zero for this commutator. Set
\begin{equation}\label{eq:hn}
H_N:=\int_{M_N}\tr K_{[A^\chi,B^\chi]}(x,x)\,dx, 
\end{equation}
so that
\[
\rTr[A,B] = \rTr[A^{\chi},B^{\chi}] =\lim_{N\to\infty}H_N.
\]
Write $M$ as the union of $M_{N}$ and its complement, and split up the integral \eqref{eq:hn} accordingly. The integrals over $M_N\times M_N$ cancel out (by interchanging $x,y$ and using cyclicity of the fiberwise trace) so we are left with
\begin{align}
H_N=\int_{M_N}\int_{M\setminus M_N}\tr\bigl(&K_{A^\chi}(x,y)K_{B^\chi}(y,x)\notag\\
&-K_{B^\chi}(x,y)K_{A^\chi}(y,x)\bigr)\,dy\,dx. \label{eq:ren-tr-boundary-integral}
\end{align}
These integrals are absolutely convergent because the product of two kernels on $W_j\times W_\ell$ is bounded by $C_s(1+|j-\ell|)^{-s}$ for any $s\geq0$ by \eqref{eqn-sus-schwartz}.

For $m\in\mathbb Z$, define
\begin{equation}\label{eq:ren-tr-Im-definition}
I_m:=\int_{W\times W} \tr\bigl(K_P(T^m x,y)K_Q(y,T^m x)\bigr)\,dx\,dy.
\end{equation}
Because $f$ is bounded on $W$ and $f(T^m x)=f(x)+m$, \eqref{eqn-sus-schwartz} implies the uniform estimate
\begin{equation}\label{eq:ren-tr-Im-bound}
|I_m|\leq C\,\nu_{2,0}(P)\nu_{2,0}(Q)(1+|m|)^{-4}.
\end{equation}
For $N\geq J$, the factor $\chi(y)$ in
\eqref{eq:ren-tr-boundary-integral} is $1$, and $\chi(x)$ is $1$ except possibly on the finite collection $W_0,\ldots,W_{J-1}$. Replacing the remaining cutoff factors by $1$ makes an error $\varepsilon_N$ satisfying
\[
|\varepsilon_N|\leq C\,\nu_{2,0}(P)\nu_{2,0}(Q)
\sum_{j=0}^{J-1}\sum_{\ell=N+1}^{\infty}(1+\ell-j)^{-4} \to 0 \text{ \, as \, } N \to \infty.
\]
Using deck-periodicity, the integral of the first
kernel product in \eqref{eq:ren-tr-boundary-integral} over $W_j\times W_\ell$ is $I_{j-\ell}$. Interchanging the two integration variables, using periodicity, and cyclically permuting the fiberwise trace gives
\[
\int_{W\times W}\tr\bigl(K_Q(T^m x,y)K_P(y,T^m x)\bigr)\,dx\,dy = I_{-m}.
\]
Thus the integral of the second product in \eqref{eq:ren-tr-boundary-integral} is $I_{\ell-j}$.

For each $k\geq1$, there are exactly $\min(k,N+1)$ pairs $(j,\ell)$ with $0\leq j\leq N$, $\ell\geq N+1$ and $\ell-j=k$. Therefore
\begin{equation}\label{eq:ren-tr-crossing-count}
H_N=\sum_{k=1}^{\infty}\min(k,N+1)(I_{-k}-I_k) +\varepsilon_N.
\end{equation}
Since \eqref{eq:ren-tr-Im-bound} implies $\sum_{m\in\mathbb Z}|m I_m|<\infty$, taking the limit as $N \to \infty$ yields
\begin{equation}\label{eq:ren-tr-defect-series}
\rTr[A,B]=-\sum_{m\in\mathbb Z}mI_m.
\end{equation}
This is the same as the expression in \cite[eqn. (21)]{mrs1} (or rather, \eqref{eq:ren-tr-defect-series} is the trace of that expression). Mrowka-Ruberman-Saveliev explain how this leads precisely to the trace defect \eqref{eq:trace-defect}, and the same argument is valid in the present setting. This completes the proof.

\end{proof}

\section{Fredholm theory}\label{sec-2}

\subsection{Fredholmness criterion}

The Fredholm criterion for end-periodic elliptic differential operators was originally established by Taubes: such an operator is Fredholm if and only if its indicial family is invertible (with all operators acting on the appropriate $L^{2}$-Sobolev spaces). See \cite[Lemma 4.3]{tau}. On the other hand, \cite[Proposition 1]{mr2} establishes the same Fredholmness criterion for operators of the form $I+A$ where $A$ is a cusp smoothing operator. 

In this section, we will extend the Fredholm criterion to end-periodic operators of the form $I+A$ where $A\in\Psi_{\mathrm{ep}}^{-\infty}(M;E)$. One direction of this result is particularly explicit: assuming the indicial family is invertible, we construct a two-sided parametrix of the same form, with error terms in $\dt{\Psi}_{\mathrm{ep}}^{-\infty}$.

\begin{prop}[Fredholm criterion]\label{ep-theory-prop-1}
Let $A\in\Psi_{\mathrm{ep}}^{-\infty}(M;E)$. Then
\[
I+A:L^{2}(M;E)\to L^{2}(M;E)
\]
is Fredholm if and only if
\[
I+A_{\theta}: L^{2}(Y;E_Y) \to L^{2}(Y;E_Y)
\]
is invertible for every $\theta \in \mathbb{R}$. In this case, $I+A$ admits a two-sided parametrix of the form
\[
I+B, \hspace{3mm} B\in\Psi_{\mathrm{ep}}^{-\infty}(M;E),
\]
whose errors lie in the residual ideal:
\[
\begin{aligned}
(I+A)(I+B)&=I-R_{+},\\
(I+B)(I+A)&=I-R_{-},
\end{aligned}
\hspace{3mm} R_{\pm}\in\dt{\Psi}_{\mathrm{ep}}^{-\infty}(M;E).
\]
\end{prop}



In order to establish the Fredholm criterion we will use two lemmas below. The first lemma is a ``closure under inverses'' property for indicial families of the form $I+A_{\theta}$.

\begin{lemma}[Closure under inverses]\label{ep-theory-lemma-2}
Let $\widetilde A\in\Psi_{\mathbb Z}^{-\infty}
(\widetilde Y;\widetilde E)$, with indicial family
\[
A_{\theta}\in C^\infty_{\mathrm{FL}}
(\mathbb R;\Psi^{-\infty}(Y;E_Y)),
\]
Suppose that $T_\theta:=I+A_\theta$ is invertible on
$L^2(Y;E_Y)$ for all $\theta\in\mathbb{R}$. Then $I+\widetilde A$ is invertible on $L^2(\widetilde Y;\widetilde E)$, and there is a unique
$\widetilde B\in\Psi_{\mathbb Z}^{-\infty} (\widetilde Y;\widetilde E)$ such that
\[
(I+\widetilde A)^{-1}=I+\widetilde B.
\]
The indicial family of $\widetilde B$ is
\[
B_\theta=(I+A_\theta)^{-1}-I =-T_\theta^{-1}A_\theta,
\]
and it belongs to $C^\infty_{\mathrm{FL}}(\mathbb R;\Psi^{-\infty}(Y;E_Y))$.
\end{lemma}


\begin{proof}
Fix $m\geq 0$. Since $A_\theta$ is smoothing,
$T_\theta=I+A_\theta$ defines a bounded operator on
$H^m(Y;E_Y)$. It is invertible on this space:
indeed, if $v\in H^m$ and $u=T_\theta^{-1}v\in L^2$, then
\[
u=v-A_\theta u\in H^m.
\]
Injectivity follows from invertibility on $L^2$.
The map $\theta\mapsto T_\theta$ is smooth with
values in $\mathcal B(H^m)$. Smooth inversion in
this Banach algebra therefore shows that
$\theta\mapsto T_\theta^{-1}$ is smooth in
$\mathcal B(H^m)$.

Now put
\[
B_\theta=T_\theta^{-1}-I =-T_\theta^{-1}A_\theta.
\]
For every $m$, the family $A_\theta$ is smooth
with values in $\mathcal B(H^{-m},H^m)$.
Consequently $B_\theta$ is also smooth in
$\mathcal B(H^{-m},H^m)$. Since $m$ is arbitrary,
the standard Sobolev characterization of smooth
kernels on the compact manifold $Y$ implies
$B_{(\cdot)}\in
C^\infty(\mathbb R;\Psi^{-\infty}(Y;E_Y))$. Finally, the identity
$T_{\theta+2\pi}=CT_\theta C^{-1}$ gives
$B_{\theta+2\pi}=CB_\theta C^{-1}$.
Thus $B_{(\cdot)} \in C^\infty_{\mathrm{FL}}\left(\mathbb{R}; \Psi^{-\infty}(Y;E_Y)\right)$. The FL isomorphism \eqref{eqn-sus-FL-algebra-isomorphism} gives the asserted inverse on $\widetilde Y$.
\end{proof}

For the next lemma, we fix a smooth end cutoff $\chi\in C^{\infty}(M)$ such that
\[
\chi=0 \text{\, near the compact part, \,}
\hspace{3mm} \chi=1 \text{\, sufficiently far out the end. \,}
\]
Recall that $\Psi_{\mathrm{ep}}^{-\infty}(M;E)$ consists of operators of the form $A=\chi\widetilde A\chi+R$, where $\widetilde A \in \Psi_{\mathbb{Z}}^{-\infty}(\widetilde Y;\widetilde E)$ and $R\in\dt{\Psi}_{\mathrm{ep}}^{-\infty}(M;E)$. The normal operator is $N(A)=\widetilde A$.

\begin{lemma}[Escaping translates]\label{ep-theory-lemma-3}
Let $\widetilde A \in \Psi_{\mathbb{Z}}^{-\infty}
(\widetilde Y;\widetilde E)$ and $u\in C_{c}^{\infty}(\widetilde Y;\widetilde E).$ Then we have
\begin{equation}\label{ep-eqn-2}
T^{k}u\rightharpoonup0 \text{\, weakly in \,} L^{2}(\widetilde Y;\widetilde E),
\end{equation}
as $k \to \infty$, and
\begin{equation}\label{ep-eqn-3}
\|(1-\chi)\widetilde{A}T^{k}u\|_{L^{2}(\widetilde{Y})} \to0.
\end{equation}
Moreover, if $R\in\dt{\Psi}_{\mathrm{ep}}^{-\infty}(M;E)$, then
\begin{equation}\label{ep-eqn-4}
\|RT^{k}u\|_{L^{2}(M)} \to0.
\end{equation}
Here, for $k$ sufficiently large, $T^{k}u$ is regarded as a section supported in the periodic end of $M$.
\end{lemma}

\begin{proof}
The supports of $T^ku$ leave every compact set and
$\|T^ku\|_{L^2}=\|u\|_{L^2}$, proving \eqref{ep-eqn-2}. Extend $\chi$ by zero to the full cyclic cover, so that $1-\chi$ vanishes on $W_j$ for all $j>K$, for some $K$. Choose $J$ with $\operatorname{supp}u\subset\bigcup_{|m|\leq J}W_m$.
For $j\leq K$, $x\in W_j$, and $y\in\operatorname{supp}(T^ku)$, equation \eqref{eq:uniform-weight-comparability-2} gives $\langle f(y)-f(x)\rangle\asymp\langle k-j\rangle$ uniformly for all sufficiently large $k$. Relative Schwartz decay therefore gives
\[
\|(1-\chi)\widetilde A T^ku\|_{L^2}^2
\leq C_{N,u}\sum_{j\leq K}\langle k-j\rangle^{-2N}
\leq C'_{N,u}\langle k\rangle^{1-2N},
\]
for every $N \geq 1$. This proves \eqref{ep-eqn-3}, with rapid decay in $k$. Finally, $R$ is compact by Lemma \ref{lem:ep-residual-trace-class}, so weak convergence in \eqref{ep-eqn-2} implies \eqref{ep-eqn-4}.
\end{proof}

Now we will establish Proposition \ref{ep-theory-prop-1}.
\begin{proof}[Proof of Prop. \ref{ep-theory-prop-1}]
First, suppose that 
\[
I+A_{\theta}: L^{2}(Y;E_Y) \to L^{2}(Y;E_Y)
\]
is invertible for every $\xi\in S^{1}$. We will construct a two-sided parametrix for $I+A$. Write
\[
A=\chi\widetilde A\chi+R, \hspace{3mm} R\in\dt{\Psi}_{\mathrm{ep}}^{-\infty}(M;E).
\]
By Lemma \ref{ep-theory-lemma-2}, there exists $\widetilde B
\in \Psi_{\mathbb Z\mathrm{-sus}}^{-\infty} (\widetilde Y;\widetilde E)$ such that $I+\widetilde{B} = (I+\widetilde{A})^{-1}$. Set $ B:=\chi\widetilde B\chi$. Then
\[
(I+A)(I+B)-I = \chi\widetilde A\chi + \chi\widetilde B\chi + \chi\widetilde A\chi^{2}\widetilde B\chi + R + R\chi\widetilde B\chi.
\]
Using $\widetilde A+\widetilde B+\widetilde A\widetilde B=0$, we obtain
\begin{equation}\label{ep-eqn-5}
(I+A)(I+B)-I = \chi\widetilde A(\chi^{2}-1)\widetilde B\chi + R + R\chi\widetilde B\chi.
\end{equation}
Extend $\chi$ by zero to the full cyclic cover. Then $\chi^2-1$ is bounded and vanishes sufficiently far along the positive end. The third Schwartz convolution estimate of Lemma \ref{schwartz-convolution-lemma} therefore shows that
\[
\chi\widetilde A(\chi^{2}-1)\widetilde B\chi \in \dt{\Psi}_{\mathrm{ep}}^{-\infty}(M;E).
\]
The other terms in \eqref{ep-eqn-5} also belong to $\dt{\Psi}^{-\infty}$ because $\dt{\Psi}^{-\infty}$ is a two-sided ideal. Thus
\[
(I+A)(I+B)=I-R_{+}, \hspace{3mm} R_{+}\in\dt{\Psi}_{\mathrm{ep}}^{-\infty}(M;E).
\]
The same argument shows that
\[
(I+B)(I+A)=I-R_{-}, \hspace{3mm} R_{-}\in\dt{\Psi}_{\mathrm{ep}}^{-\infty}(M;E).
\]
By Lemma \ref{lem:ep-residual-trace-class}, $R_{\pm}$ are compact. Thus $I+A$ is Fredholm with two-sided parametrix given by $I+B$. 

For the converse, suppose that there is some $\theta_{0} \in \mathbb{R}$ such that $I+A_{\theta_{0}}$ is not invertible. We will use Lemma \ref{ep-theory-lemma-3} to construct a sequence $q_{j} \in C_{c}^{\infty}(M;E)$ such that
\begin{equation}\label{ep-eqn-8}
\|q_j\|=1, \hspace{3mm} q_j\rightharpoonup0, \hspace{3mm} (I+A)q_j\to0,
\end{equation}
as $j\to\infty$, which implies that $I+A$ cannot be Fredholm. Since $A_{\theta_{0}}$ is a compact operator, $I+A_{\theta_{0}}$ has index zero. Noninvertibility therefore implies that $\ker(I+A_{\theta_{0}})\neq\{0\}$. Choose any nonzero $v_{0}\in\ker(I+A_{\theta_{0}})$. In fact $v_{0}$ is smooth, because $v_{0}=-A_{\theta_{0}}v_{0}$. Choose a small arc $U$ containing $\theta_{0}$, and extend $v_{0}$ to a smooth local section
\[
\theta \mapsto v(\theta), \hspace{3mm} v(\theta_{0}) = v_{0}, \hspace{3mm} \theta \in U.
\]
By continuity, we have
\[
\|(I+A_{\theta})v(\theta)\|_{L^{2}} \to 0 \text{\, as \,} \theta \to \theta_{0}.
\]
Choose functions $\varphi_j\in C_{c}^{\infty}(U)$
with $\|\varphi_jv\|_{L^{2}}=1$, and with supports shrinking to $\theta_{0}$ as $j\to\infty$. Define $\widehat u_j(\theta)
:= \varphi_j(\theta)v(\theta)$. Then
\[
\|\widehat u_j\|_{L^{2}}=1,\hspace{3mm}
\left\| (I+A_{(\cdot)}) \widehat u_j\right\|_{L^{2}}
\to 0
\]
as $j \to \infty$. Let $u_j:=\mathcal F^{-1}\widehat u_j$,
so that
\[
\|u_j\|_{L^{2}(\widetilde Y)}=1,
\hspace{3mm} \|(I+\widetilde{A})u_j\|_{L^{2}(\widetilde Y)} \to0.
\]
Since $C_{c}^{\infty}(\widetilde Y;\widetilde E)$ is dense in $L^{2}(\widetilde Y;\widetilde E)$ and $I+\widetilde{A}$ is bounded, we may replace $u_j$ by compactly supported smooth sections $w_j\in C_{c}^{\infty}(\widetilde Y;\widetilde E)$
such that
\[
\|w_j\|_{L^{2}}\to1, \hspace{3mm}
\|(I+\widetilde{A}) w_j\|_{L^{2}}\to0.
\]
After normalizing, we may assume $\|w_j\|_{L^{2}}=1$.

\medskip

For each $j$, we choose a nonnegative integer $k_{j}$, and set
$q_j:=T^{k_j}w_j$. We choose the sequence $k_j$ so that $k_{j} \to \infty$, and
\begin{enumerate}
\item $\mathrm{supp}(q_j)$ lies in the region where $\chi=1$;

\item $\mathrm{supp}(q_j)$ leaves every compact subset of $M$;

\item $\|(1-\chi)\widetilde A q_j\|_{L^{2}} <\frac1j$;

\item $\|Rq_j\|_{L^{2}} <\frac1j$.
\end{enumerate}
Here we have fixed a decomposition $A=\chi\widetilde A\chi+R$ with $R\in\dt{\Psi}_{\mathrm{ep}}^{-\infty}(M;E)$. This construction of $q_{j}$ is made possible by Lemma \ref{ep-theory-lemma-3}. In particular, we have
\[
\|q_j\|=1, \hspace{3mm} q_j\rightharpoonup 0 \text{\, weakly in \,} L^{2}(M;E),
\]
and, since $\widetilde A$ commutes with deck translation,
\[
\|(I+\widetilde{A})q_j\| = \|(I+\widetilde{A})w_j\|
\to 0.
\]
Also, because $\chi q_j=q_j$, we have
\[
(I+A)q_j = q_j+\chi\widetilde A q_j+Rq_j = (I+\widetilde A)q_j - (1-\chi)\widetilde A q_j + Rq_j.
\]
We conclude that $\|(I+A)q_j\|_{L^{2}(M)} \to 0$ as $j\to\infty$. We have constructed the sequence $q_{j} \in C_{c}^{\infty}(M;E)$ so that \eqref{ep-eqn-8} holds. No Fredholm operator can admit such a sequence, thus $I+A$ is not Fredholm. 

\end{proof}

\subsection{Winding number formula}

With the Fredholm criterion and the trace defect formula in hand, we can now establish an end-periodic version of \cite[Proposition 1]{mr2}. Just like the cusp case, the Fredholm index of an end-periodic smoothing perturbation of the identity can be expressed in terms of its indicial family. The additional feature in the end-periodic setting is that the renormalized trace defect formula has an extra term, however in this case the extra term simply vanishes.

\begin{prop}[Winding number formula]\label{prop-winding}
Let $A \in \Psi_{\mathrm{ep}}^{-\infty}(M;E)$. If $I+A$ is Fredholm, then its index is computed as the winding number of the Fredholm determinant of the indicial family:
\begin{align*}
\ind(I+A) &= -\mathrm{wn}(\mathrm{det}(I+A_{\xi})) \\
&:= -\frac{1}{2\pi i}\oint_{S^{1}}\Tr\left((I+A_{\xi})^{-1}\frac{\del A_{\xi}}{\del\xi}\right)d\xi.
\end{align*}
\end{prop}
\begin{proof}
The same argument as in Melrose-Rochon can be applied here.
Let $I+B$ denote a two-sided parametrix for $I+A$ from Proposition \ref{ep-theory-prop-1}. Then
\[
[A,B] = [I+A,I+B] = R_{-} - R_{+} \in \dt{\Psi}_{\mathrm{ep}}^{-\infty}(M;E).
\]
In particular, $[A,B]$ is trace-class by Lemma \ref{lem:ep-residual-trace-class}. Therefore Calderon's formula for the index gives
\[
\ind(I+A) = \Tr[A,B] = \rTr[A,B].
\]
The only new point in the end-periodic setting is that we use the trace defect formula from Proposition \ref{renormalized-trace-lemma}, which involves an additional term:
\begin{align*}
\rTr[A,B] &= -\frac{1}{2\pi i}\oint_{S^{1}}\Tr\left(\frac{\partial A_{\xi}}{\partial \xi}B_{\xi} \right)d\xi \\
&\hspace{6mm} + \frac{1}{2\pi i}\oint_{S^{1}}\int_{W}f(x)\cdot\tr(K_{[A_{\xi},B_{\xi}]}(x,x))\,dx\,\frac{d\xi}{\xi}.
\end{align*}
However, the second term vanishes because $(I+B_{\xi})=(I+A_{\xi})^{-1}$ implies $[A_{\xi},B_{\xi}]=0$.
\end{proof}

\begin{remark}
There is a typo in \cite[Proposition 1]{mr2}. There should be a minus sign in the first line of their equation (1.4).
\end{remark}

\subsection{Reduced index map}

We consider the subspace of Fredholm operators which are end-periodic smoothing perturbations of the identity:
\[
\mathcal F_{\mathrm{ep}}(M;E) := \bigl\{ A\in\Psi_{\mathrm{ep}}^{-\infty}(M;E): I+A\text{ is Fredholm on }L^{2}(M;E) \bigr\}.
\]
We consider the Fredholm index map on this space. Our goal in this section to establish the following end-periodic version of \cite[Proposition 2]{mr2}:
\begin{prop}\label{prop:ep-reduced-index}
The index map
\[
\mathrm{ind}: \mathcal F_{\mathrm{ep}}(M;E) \to \mathbb Z, \hspace{3mm} A\mapsto \mathrm{ind}(I+A),
\]
is surjective. Its zero set is
\[
\begin{split}
\ker(\mathrm{ind}) = \bigl\{
A\in\Psi_{\mathrm{ep}}^{-\infty}(M;E):
&\ \text{there exists } S\in\dt\Psi_{\mathrm{ep}}^{-\infty}(M;E) \\
&\ \text{such that } (I+A+S)^{-1}=I+B \\ 
&\ \text{for some } B\in\Psi_{\mathrm{ep}}^{-\infty}(M;E)
\bigr\}.
\end{split} 
\numberthis \label{eqn-index-zero-set}
\]
\end{prop}

In order to prove Proposition \ref{prop:ep-reduced-index}, we will first establish a sequence of lemmas. In Lemma \ref{lem:ep-residual-trace-class}, we showed that any operator $S \in \dt\Psi_{\mathrm{ep}}^{-\infty}(M;E)$ is trace-class on $L^{2}$. The following lemma describes the mapping properties of such an operator in more detail.

\begin{definition}[Rapidly decaying smooth sections]
We let $\dt{C}_{\mathrm{ep}}^{\infty}(M;E)$ denote the space of smooth sections $u\in C^{\infty}(M;E)$ such that, for every $N,k\geq0$,
\[
\sup_{x\in M} \langle\rho(x)\rangle^{N} |\nabla^{k}u(x)| <\infty.
\]
\end{definition}

\begin{lemma}[Mapping properties]\label{lem:dot-ep-mapping}
Let $R\in\dt\Psi_{\mathrm{ep}}^{-\infty}(M;E)$. Then:

\begin{enumerate}
\item
$R$ defines a continuous map
\[
R:L^{2}(M;E) \to \dt{C}_{\mathrm{ep}}^{\infty}(M;E);
\]
\item
for every $r,s,N_{1},N_{2}\geq0$, $R$ extends continuously as an operator
\[
R: \langle\rho\rangle^{N_{2}}H^{-r}(M;E) \to \langle\rho\rangle^{-N_{1}}H^{s}(M;E).
\]
\end{enumerate}
\end{lemma}
\begin{proof}
For $u\in L^{2}(M;E)$, the Cauchy-Schwarz inequality gives
\[
|\nabla_{x}^{\alpha}Ru(x)| \leq \left(
\int_{M} |\nabla_{x}^{\alpha}K_{R}(x,y)|^{2}\,dy \right)^{1/2}
\|u\|_{L^{2}}.
\]
The smoothing kernel estimates from Lemma \ref{ep-theory-lemma-1} imply that, for every $N$,
\[
\sup_{x\in M} \langle\rho(x)\rangle^{N} \left( \int_{M}
|\nabla_{x}^{\alpha}K_{R}(x,y)|^{2}\,dy \right)^{1/2} <\infty.
\]
This shows that $R$ defines a continuous map $R:L^{2}(M;E)\to
\dt C_{\mathrm{ep}}^{\infty}(M;E)$. The weighted Sobolev mapping properties then follow in the same way, using derivatives in both kernel variables and duality.
\end{proof}

The next lemma is a useful tool for constructing finite-rank elements of $\dt\Psi_{\mathrm{ep}}^{-\infty}(M;E)$.

\begin{lemma}[Finite-rank residual operators]
\label{lem:finite-rank-dot-ep}
Suppose that
\[
u_{1},\ldots,u_{N}, v_{1},\ldots,v_{N} \in \dt{C}_{\mathrm{ep}}^{\infty}(M;E).
\]
Then the finite-rank operator $R$ given by
\begin{equation}\label{eqn-10}
u \mapsto Ru = \sum_{j=1}^{N} \langle u,u_{j}\rangle_{L^{2}}v_{j}
\end{equation}
belongs to $\dt\Psi_{\mathrm{ep}}^{-\infty}(M;E)$.
\end{lemma}

\begin{proof}
Its Schwartz kernel is
\[
K_{R}(x,y) = \sum_{j=1}^{N} v_{j}(x)\otimes u_{j}(y)^{*}.
\]
Since each $u_{j}$ and $v_{j}$, together with all derivatives, decays faster than every inverse power of $\rho$, the kernel $K_{R}$ satisfies the defining estimates \eqref{absolute-schwartz-decay} for the algebra $\dt{\Psi}^{-\infty}_{\mathrm{ep}}(M;E)$. 
\end{proof}

For the next lemma, let
\[
T=I+A, \hspace{3mm} A\in\Psi_{\mathrm{ep}}^{-\infty}(M;E),
\]
and suppose that $T$ is Fredholm. By the Fredholm criterion, there is a parametrix
\[
P=I+B, \hspace{3mm} B\in\Psi_{\mathrm{ep}}^{-\infty}(M;E),
\]
such that
\[
\begin{aligned}
PT&=I-R_{-},\\
TP&=I-R_{+},
\end{aligned}
\hspace{3mm} R_{\pm} \in \dt{\Psi}_{\mathrm{ep}}^{-\infty}(M;E).
\]

\begin{lemma}[Rapid decay of nullvectors]\label{lem:nullvectors-dot}
With the above setup, we have
\[
\ker T \subset \dt{C}_{\mathrm{ep}}^{\infty}(M;E), \hspace{3mm} \ker T^{*} \subset \dt{C}_{\mathrm{ep}}^{\infty}(M;E).
\]
\end{lemma}
\begin{proof}
If $u\in\ker T$, then the identity $PT=I-R_{-}$ gives $u=R_{-}u$. Lemma \ref{lem:dot-ep-mapping} therefore implies that $u\in\dt{C}_{\mathrm{ep}}^{\infty}(M;E)$. If $v\in\ker T^{*}$, taking the adjoint of the identity $TP =I-R_{+}$ yields $P^{*}T^{*}=I-R_{+}^{*}$. It follows that $v=R_{+}^{*}v$, and since $R_{+}^{*}
\in \dt\Psi_{\mathrm{ep}}^{-\infty}(M;E)$, the same argument gives
$v\in\dt{C}_{\mathrm{ep}}^{\infty}(M;E)$.
\end{proof}

We shall also use the following fact about sandwiching a bounded operator between two residual operators.
\begin{lemma}
\label{lem:rapid-sandwich}
Let $R_{1},R_{2} \in \dt\Psi_{\mathrm{ep}}^{-\infty}(M;E)$, and let $D$ be a bounded operator on $L^{2}(M;E)$. Then
\[
R_{1}DR_{2} \in \dt\Psi_{\mathrm{ep}}^{-\infty}(M;E).
\]
\end{lemma}

\begin{proof}
Since the periodic end consists of compact segments of constant volume and $\langle\rho(z)\rangle\asymp\langle j\rangle$ on the $j$th segment by \eqref{eq:uniform-weight-comparability-3}, we have
\[
\int_M\langle\rho(z)\rangle^{-2}\,dz<\infty.
\]
For any $N\geq0$ and any covariant derivative $\nabla_x^\alpha$, the absolute Schwartz estimates \eqref{absolute-schwartz-decay} give
\[
|\nabla_x^\alpha K_{R_1}(x,z)|
\leq \mu_{N+1,|\alpha|}(R_1)\langle\rho(x)\rangle^{-N}\langle\rho(z)\rangle^{-1},
\]
because $\rho\geq0$ and $\langle\rho(x)+\rho(z)\rangle^{N+1}\geq \langle\rho(x)\rangle^N\langle\rho(z)\rangle$.
The analogous bound for $R_2$ therefore yields
\begin{equation}\label{eq:rapid-sandwich-L2-bounds}
\begin{aligned}
\|\nabla_x^\alpha K_{R_1}(x,\cdot)\|_{L^2} &\leq c\,\mu_{N+1,|\alpha|}(R_1) \langle\rho(x)\rangle^{-N},\\
\|\nabla_y^\beta K_{R_2}(\cdot,y)\|_{L^2} &\leq c\,\mu_{N+1,|\beta|}(R_2) \langle\rho(y)\rangle^{-N}.
\end{aligned}
\end{equation}
Here $\mu_{N,r}$ denotes the residual Fr\'echet seminorm \eqref{eq:residual-frechet-seminorms}.

The kernel maps $x\mapsto K_{R_1}(x,\cdot)$ and
$y\mapsto K_{R_2}(\cdot,y)$ are smooth as $L^2$-valued maps, since the same bounds hold for all derivatives and justify $L^2$ differentiation. Thus $R_1DR_2$ has smooth kernel
\[
K_{R_1DR_2}(x,y) =\int_M K_{R_1}(x,z)\,
\bigl[D(K_{R_2}(\cdot,y))\bigr](z)\,dz,
\]
where $D$ acts columnwise on the $E_y$-valued kernel. Differentiating this pairing and applying Cauchy-Schwarz gives
\[
\begin{aligned}
|\nabla_x^\alpha\nabla_y^\beta K_{R_1DR_2}(x,y)|
&\leq \|D\|_{\mathcal B(L^2)} \|\nabla_x^\alpha K_{R_1}(x,\cdot)\|_{L^2} \|\nabla_y^\beta K_{R_2}(\cdot,y)\|_{L^2}\\
&\leq c\|D\|_{\mathcal B(L^2)} \mu_{N+1,|\alpha|}(R_1)\mu_{N+1,|\beta|}(R_2) \langle\rho(x)\rangle^{-N}\langle\rho(y)\rangle^{-N},
\end{aligned}
\]
where $c:=\int_M\langle\rho(z)\rangle^{-2}\,dz$.
Since $\langle\rho(x)+\rho(y)\rangle^N\leq
C_N\langle\rho(x)\rangle^N\langle\rho(y)\rangle^N$,
we conclude that, for every $N,r\geq0$,
\begin{equation}\label{eq:rapid-sandwich-seminorm-estimate}
\mu_{N,r}(R_1DR_2) \leq C_{N,r}\|D\|_{\mathcal B(L^2)} \mu_{N+1,r}(R_1)\mu_{N+1,r}(R_2).
\end{equation}
In particular this implies $R_1DR_2\in\dt\Psi_{\mathrm{ep}}^{-\infty}(M;E)$.
\end{proof}


\begin{lemma}[Inverse-closedness of $I+\Psi_{\mathrm{ep}}^{-\infty}$]\label{lem:ep-inverse-closed}
Let $C\in\Psi_{\mathrm{ep}}^{-\infty}(M;E)$, and suppose that $T:=I+C$ is invertible on $L^{2}(M;E)$. Then $T^{-1}=I+D$ for some $D\in\Psi_{\mathrm{ep}}^{-\infty}(M;E)$.
\end{lemma}

\begin{proof}
Since $T$ is invertible, its indicial family is invertible. By the Fredholm criterion (Proposition \ref{ep-theory-prop-1}), $T$ has a parametrix
\[
P=I+D_{0}, \hspace{3mm} D_{0}\in\Psi_{\mathrm{ep}}^{-\infty}(M;E),
\]
such that
\begin{equation}\label{eqn-9}
\begin{aligned}
TP&=I-R_{+},\\
PT&=I-R_{-},
\end{aligned}
\hspace{3mm} R_{\pm} \in \dt\Psi_{\mathrm{ep}}^{-\infty}(M;E).
\end{equation}
Set $Q:=T^{-1}-P$. Multiplying the second equation in \eqref{eqn-9} on the right by $T^{-1}$ gives $Q=R_{-}T^{-1}$. Multiplying the first equation in \eqref{eqn-9} on the left by $T^{-1}$ gives $Q=T^{-1}R_{+}$. Using $T^{-1}=P+Q$, these together yield
\[
Q = R_{-}(P+Q) = R_{-}P+R_{-}T^{-1}R_{+}.
\]
The first term belongs to $\dt\Psi_{\mathrm{ep}}^{-\infty}$ because this space is a two-sided ideal in $\Psi_{\mathrm{ep}}^{-\infty}$. The second belongs to $\dt\Psi_{\mathrm{ep}}^{-\infty}$ by Lemma \ref{lem:rapid-sandwich}. Therefore $Q\in\dt\Psi_{\mathrm{ep}}^{-\infty}(M;E)$. Consequently,
\[
T^{-1} = I+(D_{0}+Q),
\]
with $D_{0}+Q \in \Psi_{\mathrm{ep}}^{-\infty}(M;E)$.
\end{proof}

Now with Lemmas \ref{lem:dot-ep-mapping} -- \ref{lem:ep-inverse-closed} in hand, we will prove Proposition \ref{prop:ep-reduced-index}.

\begin{proof}[Proof of Prop. \ref{prop:ep-reduced-index}]
Suppose that there exists $R \in \dt{\Psi}_{\mathrm{ep}}^{-\infty}(M;E)$ and $B \in \Psi_{\mathrm{ep}}^{-\infty}(M;E)$ such that
\[
(I+A+R)^{-1} = I+B.
\]
Then $I+A+R$ is invertible, and therefore has index zero. Moreover, $R$ is compact on $L^{2}$, so $I+A$ is Fredholm with index
\[
\ind(I+A) = \ind(I+A+R) = 0.
\]
This proves one inclusion in \eqref{eqn-index-zero-set}.

Conversely, suppose that $A \in \mathcal{F}_{\mathrm{ep}}(M;E)$ has $\ind(I+A) = 0$. Set $T:= I+A$ and define
\[
\mathcal{N} := \ker T, \hspace{3mm} \mathcal{C} := \ker T^{*} = \mathrm{ran}(T)^{\perp}.
\]
Since $T$ is Fredholm and has index zero, we have
\[
d := \dim \mathcal{N} = \dim \mathcal{C} < \infty.
\]
By Lemma \ref{lem:nullvectors-dot}, we know that $\mathcal{N},\mathcal{C} \subset \dt{C}_{\mathrm{ep}}^{\infty}(M;E)$. Choose orthonormal bases
\[
u_{1},\ldots,u_{d} \text{ for } \mathcal{N}, \hspace{3mm} v_{1},\ldots,v_{d} \text{ for } \mathcal{C},
\]
and define an operator $R$ by the formula \eqref{eqn-10}:
\[
Rw := \sum_{j=1}^{d}\langle w,u_{j}\rangle_{L^{2}}v_{j}.
\]
By Lemma \ref{lem:finite-rank-dot-ep}, we have $R\in\dt\Psi_{\mathrm{ep}}^{-\infty}(M;E)$. Using the orthogonal decompositions for the domain and codomain
\[
L^{2}(M;E) = \mathcal{N}^{\perp}\oplus\mathcal{N}, \hspace{3mm} L^{2}(M;E) = \mathrm{ran}(T)\oplus\mathcal{C},
\]
we can identify $T+R$ with the block diagonal operator
\[
T+R = 
\begin{pmatrix}
T|_{\mathcal N^{\perp}} & 0\\
0 & J
\end{pmatrix},
\]
where
\[
J: \mathcal{N} \to \mathcal{C}, \hspace{3mm} Ju_{j} = v_{j}
\]
is an isomorphism. The restriction
$T|_{\mathcal{N}^{\perp}}: \mathcal{N}^{\perp} \to \mathrm{ran}(T)$ is also a bounded isomorphism, hence $T+R = I+A+R$ is invertible. 
Since $A+R\in\Psi_{\mathrm{ep}}^{-\infty}(M;E)$, Lemma \ref{lem:ep-inverse-closed} now applies and gives $(I+A+R)^{-1}=I+B$ for some $B\in\Psi_{\mathrm{ep}}^{-\infty}(M;E)$. This establishes the description \eqref{eqn-index-zero-set} for the zero set of the index. 

It remains to show that $\ind: \mathcal{F}_{\mathrm{ep}}(M;E) \to \mathbb{Z}$ is surjective. The point requiring care is that the indicial family of a periodic normal operator is not, in general, an ordinary smooth loop of smoothing operators on $Y$. Rather, it belongs to
\[
C^\infty_{\mathrm{FL}}\bigl(\mathbb R;\Psi^{-\infty}(Y;E_Y)\bigr), \hspace{3mm} A_{\theta+2\pi}=CA_\theta C^{-1},
\]
where $C$ is the unitary multiplication operator of
\eqref{eqn-sus-FL-vector-monodromy}. We will construct a family of rank-one projections satisfying this twisting condition.

Fix $k\in\mathbb Z$. Choose an open subset $U\subset Y$ trivializing $E_{Y}$, and a smooth section $\sigma:U\to\widetilde Y$ of $p$. Set
$f_U:=f\circ\sigma\in C^\infty(U;\mathbb R)$. Recall that $C$ is multiplication by $e^{2\pi if_U}$ on $U$. Choose $e\in C_c^\infty(U;E_Y)$ with $\|e\|_{L^2(Y;E_Y)}=1$,
and for $\theta\in\mathbb R$ define
\begin{equation}\label{eq:ep-surj-e-theta}
e_\theta(z):=
\begin{cases}
e^{i\theta f_U(z)}e(z),&z\in U,\\
0,&z\notin U.
\end{cases}
\end{equation}
Because $e$ has compact support in $U$, the extension by zero is smooth on $Y$, and $(\theta,z)\mapsto e_\theta(z)$ is smooth on $\mathbb R\times Y$. Moreover, $\|e_\theta\|_{L^2}=1$, and we have $e_{\theta+2\pi}=Ce_\theta$. Let $\Pi_\theta$ be the orthogonal projection onto $\mathbb Ce_\theta$:
\begin{equation}\label{eq:ep-surj-Pi-definition}
\Pi_\theta u=\langle u,e_\theta\rangle_{L^2}e_\theta.
\end{equation}
Its kernel is
\begin{equation}\label{eq:ep-surj-Pi-kernel}
K_{\Pi_\theta}(z,z')=e_\theta(z)\otimes e_\theta(z')^*.
\end{equation}
This is a smooth function of $(\theta,z,z')\in\mathbb R\times Y^2$, with support in the fixed compact set $\operatorname{supp}(e)\times\operatorname{supp}(e)$ in $(z,z')$. Therefore $\theta\mapsto\Pi_\theta$ is smooth in the
usual Fr\'echet topology on $\Psi^{-\infty}(Y;E_Y)$. Furthermore, since $C$ is unitary, we have $\Pi_{\theta+2\pi}=C\Pi_\theta C^{-1}$.
Define
\begin{equation}\label{eq:ep-surj-a-definition}
a_{k,\theta}:=(e^{-ik\theta}-1)\Pi_\theta, \hspace{3mm} g_{k,\theta}:=I+a_{k,\theta} =(I-\Pi_\theta)+e^{-ik\theta}\Pi_\theta.
\end{equation}
These satisfy the twisting condition
\[
a_{k,\theta+2\pi}=Ca_{k,\theta}C^{-1},
\hspace{3mm} g_{k,\theta+2\pi}=Cg_{k,\theta}C^{-1}.
\]
The smoothness of $\Pi_\theta$ implies that
\begin{equation}\label{eq:ep-surj-a-twisted-membership}
a_{k,(\cdot)}\in C^\infty_{\mathrm{FL}}\bigl(\mathbb R;\Psi^{-\infty}(Y;E_Y)\bigr).
\end{equation}
Since $\Pi_\theta^2=\Pi_\theta$, the decomposition
$L^2(Y;E_Y)=\mathbb Ce_\theta\oplus(\mathbb Ce_\theta)^\perp$ shows that $g_{k,\theta}$ is invertible for every $\theta$, with
\begin{equation}\label{eq:ep-surj-g-inverse}
g_{k,\theta}^{-1} =(I-\Pi_\theta)+e^{ik\theta}\Pi_\theta
=I+(e^{ik\theta}-1)\Pi_\theta.
\end{equation}
It has Fredholm determinant $\det g_{k,\theta}=e^{-ik\theta}$. Consequently,
\begin{equation}\label{eq:ep-surj-winding}
\operatorname{wind}\bigl(\det g_{k,(\cdot)}\bigr)
=\frac{1}{2\pi i}\int_0^{2\pi} \frac{\partial_\theta\det g_{k,\theta}}{\det g_{k,\theta}}\,d\theta =-k.
\end{equation}
By Theorem \ref{FL-transform-thm}, the twisted family $a_{k,(\cdot)}$ determines 
\begin{equation}\label{eq:ep-surj-normal-operator}
\widetilde A_k:=\mathcal F^{-1}a_{k,(\cdot)}\mathcal F \in\Psi_{\mathbb Z}^{-\infty}(\widetilde Y;\widetilde E).
\end{equation}
Using the cutoff extension defined in \eqref{eq:Ext-chi-topology}, set
\begin{equation}\label{eq:ep-surj-Ak-extension}
A_k:=\operatorname{Ext}_\chi(\widetilde A_k)
=\chi\widetilde A_k\chi \in\Psi_{\mathrm{ep}}^{-\infty}(M;E).
\end{equation}
The indicial family of $I+A_k$ is $g_{k,(\cdot)}$, which is invertible for all $\theta\in\mb{R}$. Thus the Fredholm criterion (Proposition \ref{ep-theory-prop-1}) implies that $A_k\in\mathcal F_{\mathrm{ep}}(M;E)$. Finally, the winding number formula (Proposition \ref{prop-winding}) and \eqref{eq:ep-surj-winding} give
\[
\operatorname{ind}(I+A_k) =-\operatorname{wind}\bigl(\det g_{k,(\cdot)}\bigr) =k.
\]
Thus the index map is surjective.
\end{proof}

\section{End-periodic smoothing groups}\label{sec-3}

We now study the topology and characteristic forms of the smoothing group associated with the algebra $\Psi_{\mathrm{ep}}^{-\infty}(M;E)$. Compared with the cusp smoothing group, the new feature is that FL indicial families satisfy the twisting \eqref{eqn- sus-FL-operator-twist}, rather than ordinary periodicity. We shall incorporate this additional structure to analyze the normal and end-periodic smoothing groups and their transgression forms.

We first recall the ordinary smoothing group for the closed manifold $Y$. We write
\[
\mathfrak g:=\Psi^{-\infty}(Y;E_Y),\hspace{3mm}
\mathcal G:=G^{-\infty}(Y;E_Y) :=\{I+A:A\in\mathfrak g,\ (I+A)^{-1}\in I+\mathfrak g\}.
\]
The group $\mathcal G$ is a Fr\'echet Lie group. By finite-rank approximation and Bott periodicity (cf. \cite[Proposition 3]{mr2}), one has
\[
\pi_{2k}(\mathcal G)=0,\hspace{3mm} \pi_{2k-1}(\mathcal G)\cong\mathbb Z,\hspace{3mm} k\geq1.
\]
The odd generators are detected by the normalized closed forms
\begin{equation}\label{beta-form-1}
\beta_k^{\mathrm{odd}}:=c_k^{\mathrm{odd}}\Tr
\bigl((g^{-1}dg)^{2k-1}\bigr),
\end{equation}
with coefficient
\begin{equation}\label{c-coeff-1}
c_k^{\mathrm{odd}} :=\frac{1}{(2\pi i)^k}\frac{(k-1)!}{(2k-1)!}.
\end{equation}
Thus we have an isomorphism $\pi_{2k-1}(\mathcal G)\to\mathbb Z$ given by $[F]\mapsto\int_{S^{2k-1}}F^*\beta_k^{\mathrm{odd}}$.

Let $C$ be the unitary multiplier $(Cv)(p(x))=e^{2\pi i f(x)}v(p(x))$ from \eqref{eqn-sus-FL-vector-monodromy}. In view of Theorem \ref{FL-transform-thm}, the suspended or normal smoothing group is the \emph{twisted loop group}
\begin{equation}\label{eq:ep-twisted-loop-group}
G_{\mathbb Z}^{-\infty}(Y;E_Y)=L_C\mathcal G :=\bigl\{g\in C^\infty(\mathbb R;\mathcal G): g_{\theta+2\pi}=Cg_\theta C^{-1}\bigr\},
\end{equation}
We also define the end-periodic and residual smoothing groups
\[
\begin{aligned}
\mathcal E &:=G_{\mathrm{ep}}^{-\infty}(M;E):=\bigl\{I+A:A\in\Psi_{\mathrm{ep}}^{-\infty}(M;E),\,
(I+A)^{-1}\in I+\Psi_{\mathrm{ep}}^{-\infty}(M;E)\bigr\},\\
\mathcal K &:=\dt G_{\mathrm{ep}}^{-\infty}(M;E) :=\bigl\{I+R:R\in\dt\Psi_{\mathrm{ep}}^{-\infty}(M;E),\,
(I+R)^{-1}\in I+\dt\Psi_{\mathrm{ep}}^{-\infty}(M;E)\bigr\}.
\end{aligned}
\]

\begin{remark}[Smooth inversion]\label{remark-smooth-inversion}
These are Fr\'echet Lie groups: their defining invertibility conditions specify open subsets of the corresponding affine Fr\'echet spaces, and multiplication and inversion are smooth. Indeed, openness follows from the continuous inclusions into $\mathcal B(L^2)$ and inverse-closedness (Lemmas \ref{ep-theory-lemma-2} and \ref{lem:ep-inverse-closed}); the same conclusion for $\mathcal K$ follows from $(I+R)^{-1}-I=-R+R(I+R)^{-1}R$ and Lemma \ref{lem:rapid-sandwich}.

Continuity of multiplication follows from the Schwartz-convolution estimates and the Fr\'echet algebra topology established above. For $L_C\mathcal G$, smooth inversion follows from smooth Banach-space inversion on Sobolev spaces and the smoothing-kernel estimates in Lemma \ref{ep-theory-lemma-2}.

For $\mathcal E$, a smoothly varying inverse indicial family yields (using the FL isomorphism) a smooth parametrix $P$ with $gP=I-R_+$ and $Pg=I-R_-$, where $R_\pm$ are residual. Then
\[
g^{-1}=P+R_-P+R_-g^{-1}R_+
\]
and the residual sandwich estimate of Lemma \ref{lem:rapid-sandwich} reduces smooth inversion to smooth inversion in $\mathcal B(L^2)$.
\end{remark}

The based normal group is
\[
G_{\mathbb Z,*}^{-\infty}(Y;E_Y) =\Omega_C\mathcal G:=\{g\in L_C\mathcal G:g_0=I\}.
\]
Evaluation at $\xi=1$ gives the exact sequence
\begin{equation}\label{eq:ep-twisted-evaluation-sequence}
1\to\Omega_C\mathcal G \to L_C\mathcal G \xrightarrow{\ \operatorname{ev}_1\ }\mathcal G \to1,\hspace{3mm} \operatorname{ev}_1(g)=g_0.
\end{equation}
It is surjective because $\mathcal G$ is path connected: any
$h\in\mathcal G$ can be joined by a smooth path, constant near its endpoints, to $ChC^{-1}$. The fibration structure of this sequence will be considered in the twisted case below.

\subsection{Untwisted case}

For comparison, if $C=I$, the normal group is the ordinary smooth free-loop group $L\mathcal G=C^\infty(S^1;\mathcal G)$. In this case, the constant loops provide a splitting for evaluation $\mathrm{ev}_{1}: L\mathcal{G} \to \mathcal{G}$, hence one obtains a isomorphism
\begin{equation}\label{eq-untwisted-splitting}
L\mathcal G\cong\mathcal G\times\Omega\mathcal G,
\hspace{3mm} g\mapsto\bigl(g(1),g(1)^{-1}g(\cdot)\bigr).
\end{equation}
Thus, by the homotopy groups of $\mathcal G$ above,
\[
\pi_{2k-1}(L\mathcal G)\cong\mathbb Z, \hspace{3mm}
\pi_{2k}(L\mathcal G)\cong\mathbb Z, \hspace{3mm} k\geq1.
\]
Let $\operatorname{ev}:S^1\times L\mathcal G\to\mathcal G$
be evaluation, and let $\operatorname{ev}_1$ denote evaluation at
$1\in S^1$. The corresponding odd and even forms are
\[
\alpha_k^{\mathrm{odd}} :=\operatorname{ev}_1^*\beta_k^{\mathrm{odd}}, \hspace{3mm}
\beta_k^{\mathrm{even}} :=\int_{S^1}\operatorname{ev}^*\beta_{k+1}^{\mathrm{odd}}.
\]
For $k\geq1$, their pairings give isomorphisms
\begin{equation}\label{eq:untwisted-odd-pairing}
\pi_{2k-1}(L\mathcal G)\to\mathbb Z, \hspace{3mm} [F]\mapsto
\int_{S^{2k-1}}F^*\alpha_k^{\mathrm{odd}},
\end{equation}
\begin{equation}\label{eq:untwisted-even-pairing}
\pi_{2k}(L\mathcal G)\to\mathbb Z, \hspace{3mm} [F]\mapsto
\int_{S^{2k}}F^*\beta_k^{\mathrm{even}}.
\end{equation}
Indeed, under \eqref{eq-untwisted-splitting} the odd pairing comes from $\mathcal G$, whereas the even pairing comes from $\Omega\mathcal G$. The latter is modeled by the cusp suspended group on $\mathbb R$: one flattens based loops near their basepoint and identifies the punctured circle with $\mathbb R$. In this way, \cite[Proposition 4]{mr2} computes the homotopy groups of the cusp suspended group. In the twisted setting, the constant-loop splitting \eqref{eq-untwisted-splitting} is no longer available; we incorporate the twisting and study the corresponding forms and homotopy groups below.


\subsection{Twisted case}

We now turn to the twisted normal smoothing group
$L_C\mathcal G$ arising from the FL transform. Unlike the ordinary free-loop group $L\mathcal{G}$ considered above, the twisting condition prevents a splitting by constant loops. Nevertheless, we will show that the homotopy groups of $L_C\mathcal G$ can be detected by natural differential forms, analogous to the odd and even characteristic forms used by \cite[Propositions 3--4]{mr2}. 

The odd forms arise by evaluating the indicial family at a fixed point of $S^1$, while the even forms are obtained by integrating the ordinary trace forms along the FL circle. Although the indicial families satisfy a twisting condition, the trace is invariant under conjugation, so these forms are well-defined. 

Let $M=Z\cup\widetilde Y_+$ be an end-periodic manifold with end modeled on $p:\widetilde Y\to Y$, and let $E\to M$ be an end-periodic Hermitian vector bundle. We write
\[
a_{\mathcal G}:\mathcal G\to I+\mathfrak g, \hspace{3mm} a_{\mathcal G}(g)=g
\]
for the tautological map and
\[
\Theta_{\mathcal G}:=a_{\mathcal G}^{-1}da_{\mathcal G} \in\Omega^1(\mathcal G;\mathfrak g)
\]
for its left Maurer-Cartan form.

For a real angle $\vartheta$, evaluation defines
$\operatorname{ev}_\vartheta:L_C\mathcal G\to\mathcal G$,
$g\mapsto g_\vartheta$, and we put
\[
a_\vartheta :=a_{\mathcal G}\circ\operatorname{ev}_\vartheta.
\]
Although $\operatorname{ev}_\vartheta$ is not a single-valued $\mathcal G$-valued map on the circle, its values at $\vartheta$ and $\vartheta+2\pi$ differ by conjugation with $C$. Since the ordinary trace is invariant under this conjugation, the forms below are single-valued on $S^1$.

Evaluation at $\xi=1$ gives odd forms
\[
\alpha_k^{\mathrm{odd}} :=\operatorname{ev}_1^*\beta_k^{\mathrm{odd}} \in\Omega^{2k-1}(L_C\mathcal G),
\]
while fiber integration gives even forms
\begin{equation}\label{eq:ep-even-fiber-integration}
\beta_k^{\mathrm{even}} :=\int_0^{2\pi}\iota_{\partial_\vartheta}\, \Tr\!\left[ c_{k+1}^{\mathrm{odd}} (a_\vartheta^{-1}d a_\vartheta)^{2k+1}\right]d\vartheta \in \in\Omega^{2k}(L_C\mathcal G).
\end{equation}
Here contraction and integration apply to the differential form on $L_C\mathcal G\times\mathbb R$ obtained by regarding $a_\vartheta$ as a function of both variables, and $d$ denotes the exterior differential on $L_{C}\mathcal G \times \mathbb{R}$.

\begin{lemma}\label{lemma-alpha-beta-forms}
Consider the $\mathfrak g$-valued forms
\begin{equation}\label{eq:omega-evaluation-pullback}
\omega_\vartheta :=a_\vartheta^{-1}d a_\vartheta
=\operatorname{ev}_\vartheta^*\Theta_{\mathcal G}
\in\Omega^1(L_C\mathcal G;\mathfrak g),
\end{equation}
\begin{equation}\label{eq:v-evaluation-pullback}
v_\vartheta :=a_\vartheta^{-1}\frac{\partial a_\vartheta}{\partial\vartheta} \in\Omega^0(L_C\mathcal G;\mathfrak g).
\end{equation}
The odd forms $\alpha_k^{\mathrm{odd}}$ are given by
\begin{equation}\label{eq:alpha-explicit}
\alpha_k^{\mathrm{odd}} =c_k^{\mathrm{odd}}\Tr(\omega_0^{2k-1}) =c_k^{\mathrm{odd}}\Tr\bigl[(a_0^{-1}da_0)^{2k-1}\bigr].
\end{equation}
The even forms $\beta_k^{\mathrm{even}}$ are given by
\begin{equation}\label{eq:beta-even-explicit}
\begin{aligned}
\beta_k^{\mathrm{even}} &=c_k^{\mathrm{even}} \int_0^{2\pi}\Tr(\omega_\vartheta^{2k} v_\vartheta)\,d\vartheta\\
&=c_k^{\mathrm{even}}\int_0^{2\pi}\Tr\!\left[(a_\vartheta^{-1}da_\vartheta)^{2k} a_\vartheta^{-1}\frac{\partial a_\vartheta}{\partial\vartheta} \right]d\vartheta,
\end{aligned}
\end{equation}
where
\[
c_k^{\mathrm{even}} :=(2k+1)c_{k+1}^{\mathrm{odd}}
=\frac{1}{(2\pi i)^{k+1}}\frac{k!}{(2k)!} =\frac{c_k^{\mathrm{odd}}}{4\pi i}.
\]
The forms $\alpha_k^{\mathrm{odd}}$ and $\beta_k^{\mathrm{even}}$ are closed.
\end{lemma}

\begin{proof}
The formula for $\alpha_k^{\mathrm{odd}}$ follows from evaluation at $\vartheta=0$. On the cylinder
$L_C\mathcal G\times\mathbb R$ the pulled-back Maurer-Cartan form has the decomposition
\[
a_\vartheta^{-1}d a_\vartheta =\omega_\vartheta+v_\vartheta\,d\vartheta.
\]
Thus only terms containing one factor $d\vartheta$ contribute to \eqref{eq:ep-even-fiber-integration}. Graded cyclicity of the trace
makes all $2k+1$ contributions equal and gives
\[
\beta_k^{\mathrm{even}}=(2k+1)c_{k+1}^{\mathrm{odd}} \int_0^{2\pi}\Tr(\omega_\vartheta^{2k}v_\vartheta)\,d\vartheta,
\]
which is \eqref{eq:beta-even-explicit}.

Both forms are closed: $\alpha_k^{\mathrm{odd}}$ is the pullback of the closed form $\beta_k^{\mathrm{odd}}$, while $\beta_k^{\mathrm{even}}$ is obtained by integrating the closed form $\mathrm{ev}^*\beta_{k+1}^{\mathrm{odd}}$.
\end{proof}

In the end-periodic setting we must calculate the homotopy groups of the normal smoothing group $L_C\mathcal G$, but here the splitting \eqref{eq-untwisted-splitting} is unavailable. To get around this obstacle we will construct an evaluation fibration and its exact sequence, and use an endpoint-flattening argument as a bridge to the ordinary based-loop homotopy theory (i.e., the untwisted case).

The following proposition describes the homotopy groups of the normal smoothing group $L_C\mathcal G$, and also that the forms $\alpha_k^{\mathrm{odd}}$ and $\beta_k^{\mathrm{even}}$ detect their generators.

\begin{prop}\label{prop:ep-free-loop-basic-classes}
For every $k\geq1$, the maps
\begin{equation}\label{eqn-13}
\pi_{2k-1}(L_C\mathcal G)\to\mathbb Z, \hspace{3mm} [F]\mapsto\int_{S^{2k-1}}F^*\alpha_k^{\mathrm{odd}},
\end{equation}
\begin{equation}\label{eqn-14}
\pi_{2k}(L_C\mathcal G)\to\mathbb Z, \hspace{3mm} [F]\mapsto\int_{S^{2k}}F^*\beta_k^{\mathrm{even}},
\end{equation}
are isomorphisms.
\end{prop}


Write
\begin{equation}\label{eq:twisted-group-automorphism}
\sigma:=\operatorname{Ad}_C:\mathcal G\to\mathcal G, \hspace{3mm} \sigma(h)=ChC^{-1},
\end{equation}
where $C$ is the unitary operator from \eqref{eqn-sus-FL-vector-monodromy}. Since conjugation by $C$ preserves $\Psi^{-\infty}(Y;E_Y)$ and its Fr\'echet topology, $\sigma$ is a smooth Lie group automorphism. Its differential on the Lie algebra $\mathfrak g=\Psi^{-\infty}(Y;E_Y)$ is $X\mapsto CXC^{-1}$. With this notation,
\[
L_C\mathcal G =\{g\in C^\infty(\mathbb R;\mathcal G):
g_{\vartheta+2\pi}=\sigma(g_\vartheta)\},\hspace{3mm}\Omega_C\mathcal G=\{g\in L_C\mathcal G:g_0=I\}.
\]
In the following, all homotopy groups will be based at the identity.

\begin{lemma}[Conjugation acts trivially on homotopy groups]
\label{lem:ep-conjugation-homotopy}
For every $j\geq1$, the induced map
\[
\sigma_*:\pi_j(\mathcal G)\to\pi_j(\mathcal G)
\]
is the identity. More precisely, for each $k\geq1$ one has
$\sigma^*\beta_k^{\mathrm{odd}}=\beta_k^{\mathrm{odd}}$; this identity implies $\sigma_*=\mathrm{id}$ on $\pi_{2k-1}(\mathcal G)\cong\mathbb Z$, while $\pi_{2k}(\mathcal G)=0$.
\end{lemma}

\begin{proof}
Let $h\in\mathcal G$ and $V\in T_h\mathcal G$. Since $C$ is independent of $h$,
\[
(\sigma(h))^{-1}d\sigma_h(V) =Ch^{-1}VC^{-1}.
\]
Consequently
\[
\sigma^*\Tr\bigl((h^{-1}dh)^{2k-1}\bigr)
=\Tr\bigl(C(h^{-1}dh)^{2k-1}C^{-1}\bigr)
=\Tr\bigl((h^{-1}dh)^{2k-1}\bigr).
\]
This proves $\sigma^*\beta_k^{\mathrm{odd}}=\beta_k^{\mathrm{odd}}$.

Let $[F]\in\pi_{2k-1}(\mathcal G,I)$ be represented by a smooth based map $F:S^{2k-1}\to\mathcal G$. We have
\[
\int_{S^{2k-1}}(\sigma\circ F)^*\beta_k^{\mathrm{odd}}
=\int_{S^{2k-1}}F^*\sigma^*\beta_k^{\mathrm{odd}}
=\int_{S^{2k-1}}F^*\beta_k^{\mathrm{odd}}.
\]
By the normalization in \eqref{beta-form-1}--\eqref{c-coeff-1}, this pairing is an isomorphism
$\pi_{2k-1}(\mathcal G)\to\mathbb Z$. It follows that
$\sigma_*[F]=[F]$ in every odd degree. The even homotopy groups vanish, so the claim follows in all positive degrees.
\end{proof}


\begin{lemma}[The twisted evaluation fibration]
\label{lem:ep-twisted-path-fibration}
Evaluation at $\vartheta=0$ defines a locally trivial fibration
\begin{equation}\label{eq:ep-twisted-evaluation-fibration}
\Omega_C\mathcal G\to L_C\mathcal G
\xrightarrow{\operatorname{ev}_1}\mathcal G,
\hspace{3mm} \operatorname{ev}_1(g)=g_0.
\end{equation}
Under the natural identification $\pi_j(\Omega_C\mathcal G)\cong\pi_{j+1}(\mathcal G)$, its connecting homomorphism is
\begin{equation}\label{eq:ep-evaluation-boundary}
\partial:\pi_j(\mathcal G)\to\pi_j(\mathcal G),
\hspace{3mm} \partial=\sigma_*-\mathrm{id},\hspace{3mm} j\geq1,
\end{equation}
up to the overall sign determined by the orientation of the interval. In particular, the associated exact sequence can be written as
\begin{equation}\label{eq:ep-twisted-wang-sequence}
\begin{aligned}
\cdots\to\pi_{j+1}(\mathcal G) &\xrightarrow{\ \sigma_*-\mathrm{id}\ } \pi_{j+1}(\mathcal G) \to\pi_j(L_C\mathcal G)\\
&\xrightarrow{\ (\operatorname{ev}_1)_*\ }
\pi_j(\mathcal G) \xrightarrow{\ \sigma_*-\mathrm{id}\ } \pi_j(\mathcal G)\to\cdots.
\end{aligned}
\end{equation}
\end{lemma}

\begin{proof}
Consider first the space of smooth twisted paths
\[
P_\sigma\mathcal G :=\{p\in C^\infty([0,2\pi];\mathcal G): p(2\pi)=\sigma(p(0))\}.
\]
It has no requirement on the higher endpoint derivatives. Restriction to $[0,2\pi]$ gives a map $L_C\mathcal G\to P_\sigma\mathcal G$.

Choose a smooth nondecreasing function $\rho:[0,2\pi]\to[0,2\pi]$ which is identically $0$ near $0$,
identically $2\pi$ near $2\pi$, and satisfies
$\rho(0)=0$, $\rho(2\pi)=2\pi$. For any $p\in P_\sigma\mathcal G$, the path $p\circ\rho$ is constant near each endpoint and its endpoint values satisfy the twisting condition. Extending it to $\mathbb R$ by
\[
(p\circ\rho)(\vartheta+m2\pi) :=\sigma^m\bigl(p(\rho(\vartheta))\bigr), \hspace{3mm} 0\leq\vartheta\leq 2\pi, \quad m\in\mathbb Z,
\]
therefore produces a smooth twisted loop. The interpolations
\[
\rho_s(\vartheta):=(1-s)\vartheta+s\rho(\vartheta), \hspace{3mm} 0\leq s\leq1,
\]
give homotopies from the identity to endpoint flattening on the smooth twisted-path space. For paths already satisfying the twisting condition, these homotopies remain in $L_C\mathcal G$. It follows that restriction $L_C\mathcal G\to P_\sigma\mathcal G$ is a homotopy equivalence, compatible with evaluation.

The evaluation $P_\sigma\mathcal G\to\mathcal G$ is locally trivial. To see this explicitly, fix $h_0\in\mathcal G$. Because $\mathcal G$ is path connected, choose a smooth path $p_0$ from $h_0$ to $\sigma(h_0)$, constant near both endpoints. For $h$ sufficiently close to $h_0$, write $h=h_0\exp X$ with $X\in\mathfrak g$ small. Let $\chi:[0,2\pi]\to[0,1]$ be smooth, equal to $0$ near $0$ and to $1$ near $2\pi$. Then
\[
p_h(\vartheta) :=p_0(\vartheta) \exp\bigl((1-\chi(\vartheta))X +\chi(\vartheta)CXC^{-1}\bigr)
\]
depends smoothly on $h$ and has endpoint values
$p_h(0)=h$ and $p_h(2\pi)=\sigma(h)$. It supplies a smooth local section of evaluation. Pointwise group multiplication by these local sections gives local trivializations, with fiber consisting of paths from $I$ to $I$; the smooth twisted-loop fiber is precisely $\Omega_C\mathcal G$. Thus \eqref{eq:ep-twisted-evaluation-fibration} is a fibration.

The connecting map can be identified from the endpoint condition. The ordinary smooth path space $P\mathcal G=C^\infty([0,2\pi];\mathcal G)$ has an endpoint fibration
\[
(\operatorname{ev}_0,\operatorname{ev}_{2\pi}):
P\mathcal G\to\mathcal G\times\mathcal G
\]
with fiber $\Omega\mathcal G$ over $(I,I)$. The space $P_\sigma\mathcal G$ is its pullback by the graph map
\[
\mathcal G\to\mathcal G\times\mathcal G, \hspace{3mm} h\mapsto(h,\sigma(h)).
\]
Under the adjunction $\pi_{j-1}(\Omega\mathcal G)\cong\pi_j(\mathcal G)$, the endpoint-fibration boundary sends a pair of homotopy classes $([a],[b])$ to $[b]-[a]$. Pulling back by the graph gives $\sigma_*[a]-[a]$, which is \eqref{eq:ep-evaluation-boundary}. Identifying the homotopy groups of the based fiber with those of $\Omega\mathcal G$ yields \eqref{eq:ep-twisted-wang-sequence}.
\end{proof}

By Lemma \ref{lem:ep-conjugation-homotopy}, every boundary map $\sigma_*-\mathrm{id}$ in the exact sequence \eqref{eq:ep-twisted-wang-sequence} vanishes. In odd degrees, this implies that the evaluation map induces an isomorphism
\[
(\operatorname{ev}_1)_*:\pi_{2k-1}(L_C\mathcal G) \xrightarrow{\ \cong\ }\pi_{2k-1}(\mathcal G) \cong\mathbb Z, \hspace{3mm} k\geq1.
\]
In even degrees, because $\pi_{2k}(\mathcal G)=0$, the inclusion induces an isomorphism
\[
0\to\pi_{2k}(\Omega_C\mathcal G) \xrightarrow{\ \iota_*\ }
\pi_{2k}(L_C\mathcal G)\to0, \hspace{3mm} k\geq1.
\]
So in order to compute the homotopy group in degree $2k$, we next look at the based subgroup $\Omega_{C}\mathcal{G}$.

\begin{lemma}[Flattening based twisted loops]
\label{lem:ep-based-loop-flattening}
Let $\Omega^{\mathrm{flat}}\mathcal G$ denote the space of smooth paths which equal $I$ on neighborhoods of $\theta=0,2\pi$. We have smooth embeddings
\[
\Omega^{\mathrm{flat}}\mathcal G \hookrightarrow \Omega\mathcal G, \hspace{3mm} \Omega^{\mathrm{flat}}\mathcal G \hookrightarrow \Omega_C\mathcal G,
\]
and both inclusions induce homotopy equivalences. Consequently, for every $j \geq 0$, we have
\begin{equation}\label{eq:ep-based-loop-adjunction}
\pi_j(\Omega_C\mathcal G) \cong \pi_{j}(\Omega^{\mathrm{flat}}\mathcal{G}) \cong\pi_j(\Omega\mathcal G) \cong\pi_{j+1}(\mathcal G).
\end{equation}
\end{lemma}

\begin{proof}
For $g\in\Omega_C\mathcal G$, one has $g_0=I$ and hence $g_{2\pi}=\sigma(I)=I$. Define
\[
R_sg(\vartheta):=g(\rho_s(\vartheta)), \hspace{3mm} \rho_s(\vartheta)=(1-s)\vartheta+s\rho(\vartheta),
\]
using the flattening function from Lemma \ref{lem:ep-twisted-path-fibration}. Near $\vartheta=0$, the function $\rho_s$ is $(1-s)\vartheta$; near $2\pi$ it is $2\pi+(1-s)(\vartheta-2\pi)$. Since
\[
\partial_\vartheta^r g_{2\pi} =C(\partial_\vartheta^r g_0)C^{-1}, \hspace{3mm} r\geq0,
\]
the chain rule shows that $R_sg$ satisfies the same twisting condition. Thus $R_s$ is a homotopy through based twisted loops, with $R_0g=g$ and $R_1g=g\circ\rho$.

The loop $R_1g$ equals $I$ near both endpoints, and hence belongs to $\Omega^{\mathrm{flat}}\mathcal G$. Such a path has all endpoint derivatives equal to zero; it therefore extends both as an ordinary smooth based loop and as a smooth twisted based loop. The same homotopy $R_s$ works for ordinary smooth based loops, since matching endpoint derivatives are likewise preserved. Moreover, if $g$ is already flat, $R_sg$ remains flat for all $s$, possibly on smaller endpoint neighborhoods. Thus flattening provides a homotopy inverse to each inclusion of the flat-loop space. This argument yields
\[
\pi_j(\Omega_C\mathcal G) \cong \pi_{j}(\Omega^{\mathrm{flat}}\mathcal{G})
\cong\pi_j(\Omega\mathcal G)
\]
The final isomorphism in \eqref{eq:ep-based-loop-adjunction} is the usual loop-space adjunction (with basepoint):
\[
\pi_j(\Omega\mathcal G) \cong\pi_{j+1}(\mathcal G).
\]
\end{proof}

Combining Lemma \ref{lem:ep-twisted-path-fibration} and Lemma \ref{lem:ep-based-loop-flattening}, we can now identify the homotopy groups of $L_{C}\mathcal{G}$. They are the same as those of $L\mathcal{G}$.

\begin{lemma}[Homotopy groups of $L_C\mathcal G$]
\label{lem:ep-twisted-homotopy-groups}
In odd degrees, the evaluation map induces an isomorphism
\[
(\operatorname{ev}_1)_*:\pi_{2k-1}(L_C\mathcal G)
\xrightarrow{\ \cong\ }\pi_{2k-1}(\mathcal G)
\cong\mathbb Z,
\hspace{3mm} k\geq1.
\]
In even degrees, inclusion of the based subgroup induces an isomorphism
\[
\iota_{*}: \pi_{2k}(\Omega_C\mathcal G)
\xrightarrow{\ \cong\ }\pi_{2k}(L_C\mathcal G)
\cong\pi_{2k+1}(\mathcal G)
\cong\mathbb Z,
\hspace{3mm} k\geq1.
\]
Moreover $\pi_0(L_C\mathcal G)\cong\pi_1(\mathcal G)\cong\mathbb Z$.
\end{lemma}

\begin{proof}
By Lemma \ref{lem:ep-conjugation-homotopy}, every map $\sigma_*-\mathrm{id}$ in the exact sequence
\eqref{eq:ep-twisted-wang-sequence} vanishes. In degree $2k-1$, exactness and $\pi_{2k}(\mathcal G)=0$ give
\[
0\to\pi_{2k-1}(L_C\mathcal G)
\xrightarrow{\ (\operatorname{ev}_1)_*\ }
\pi_{2k-1}(\mathcal G)\to0.
\]
In degree $2k$, exactness and $\pi_{2k}(\mathcal G)=0$ give
\[
0\to\pi_{2k}(\Omega_C\mathcal G)
\xrightarrow{\ \iota_*\ }
\pi_{2k}(L_C\mathcal G)\to0.
\]
Lemma \ref{lem:ep-based-loop-flattening} identifies the first group with $\pi_{2k+1}(\mathcal G)\cong\mathbb Z$. Finally, because $\mathcal G$ is connected, the degree-zero part of the same exact sequence identifies
$\pi_0(L_C\mathcal G)$ with $\pi_0(\Omega_C\mathcal G)\cong\pi_1(\mathcal G)\cong\mathbb Z$.
\end{proof}

The following lemma describes the isomorphism $\pi_{0}(L_{C}\mathcal{G}) \cong \mb{Z}$ in more detail.

\begin{lemma}[Components of $L_{C}\mathcal{G}$]
\label{lem:ep-twisted-determinant-winding}
For $g\in L_C\mathcal G$, the function $\vartheta\mapsto\det(g_\vartheta)$ descends to a smooth
map $S^1\to\mathbb C^\times$. The assignment
\begin{equation}\label{eq:ep-twisted-winding}
\operatorname{wind}(g)
:=\frac{1}{2\pi i}\int_{0}^{2\pi} \Tr\bigl(g_\vartheta^{-1}
\partial_\vartheta g_\vartheta\bigr) \,d\vartheta =\operatorname{wind}\bigl(\det g_\vartheta\bigr)
\end{equation}
induces a bijection $\pi_0(L_C\mathcal G)\xrightarrow{\ \cong\ }\mathbb Z$.
\end{lemma}

\begin{proof}
Each $g_\vartheta-I$ is smoothing and therefore trace class, so the Fredholm determinant is defined and nonzero. The conjugation-invariance of this determinant gives
\[
\det(g_{\vartheta+2\pi}) =\det(Cg_\vartheta C^{-1}) =\det(g_\vartheta).
\]
The logarithmic derivative identity for the Fredholm determinant then gives the integral in \eqref{eq:ep-twisted-winding}; its integrand is periodic by the same conjugation argument.

To show that winding detects every component, note that $\det:\mathcal G\to\mathbb C^\times$ induces an isomorphism on $\pi_1$. Indeed,
\[
\frac{1}{2\pi i}\Tr(h^{-1}dh) =\frac{1}{2\pi i}\,d\log\det h,
\]
and the left-hand side is $\beta_1^{\mathrm{odd}}$, whose pairing generates $\pi_1(\mathcal G)\cong\mathbb Z$. Since $\det\circ\sigma=\det$, the determinant maps the twisted path-space evaluation fibration to the ordinary loop-space evaluation fibration of $\mathbb C^\times$. By naturality of the low-degree exact sequences and the isomorphism on $\pi_1$, it induces an isomorphism
\[
\pi_0(L_C\mathcal G) \xrightarrow{\ \cong\ } \pi_0(L\mathbb C^\times)\cong\mathbb Z.
\]
The last isomorphism is precisely the winding number.
\end{proof}

We now combine the preceding lemmas to prove Proposition \ref{prop:ep-free-loop-basic-classes}. We show that the isomorphisms in Lemma \ref{lem:ep-twisted-homotopy-groups} are given explicitly by pairing the basic forms from Lemma \ref{lemma-alpha-beta-forms}. In particular, the isomorphisms are given by the same pairings as in the untwisted case.

\begin{proof}[Proof of Prop. \ref{prop:ep-free-loop-basic-classes}]
In odd degrees, Lemma \ref{lem:ep-twisted-homotopy-groups} gives an isomorphism
\[
(\operatorname{ev}_1)_*:\pi_{2k-1}(L_C\mathcal G)
\xrightarrow{\ \cong\ }\pi_{2k-1}(\mathcal G).
\]
Since $\alpha_k^{\mathrm{odd}}
=\operatorname{ev}_1^*\beta_k^{\mathrm{odd}}$,
for a smooth based representative $F:S^{2k-1}\to L_C\mathcal G$ one has
\[
\int_{S^{2k-1}}F^*\alpha_k^{\mathrm{odd}} =\int_{S^{2k-1}}(\operatorname{ev}_1\circ F)^*\beta_k^{\mathrm{odd}}.
\]
The right-hand side detects the generator of
$\pi_{2k-1}(\mathcal G)\cong\mathbb Z$, establishing \eqref{eqn-13}.

For even degrees, let $\iota:\Omega_C\mathcal G\hookrightarrow L_C\mathcal G$ be the embedding from Lemma \ref{lem:ep-twisted-homotopy-groups}. It induces an isomorphism $\iota_*:\pi_{2k}(\Omega_C\mathcal G)\to \pi_{2k}(L_C\mathcal G)$. It therefore suffices to evaluate $\iota^*\beta_k^{\mathrm{even}}$ on a smooth based map
\[
F:S^{2k}\to\Omega_C\mathcal G, \hspace{3mm} F(x_0)=I,
\]
where $x_0\in S^{2k}$ is the basepoint and $I$ on the right denotes the constant identity loop. Apply the reparametrization from Lemma \ref{lem:ep-based-loop-flattening} to every loop in the image of
$F$. The homotopy $F_{s}(x)(\vartheta) :=F(x)(\rho_s(\vartheta))$ consists of smooth based maps into $\Omega_C\mathcal G$; hence it preserves the homotopy class of $F$. After replacing $F$ by $F_1$, we may assume that every loop is the identity on fixed neighborhoods of $0$ and $2\pi$. The fixed-neighborhood assertion follows from the choice of the single function $\rho$, independent of $x$.

These flattened loops also belong to the ordinary smooth based loop group. Consequently the adjoint formula
\[
\widetilde F:S^{2k}\times S^1\to\mathcal G, \hspace{3mm} \widetilde F(x,e^{i\vartheta})=F(x)(\vartheta),
\]
defines a smooth map. It satisfies
\[
\widetilde F(x,1)=I\quad\text{for all }x, \hspace{3mm}
\widetilde F(x_0,\xi)=I\quad\text{for all }\xi.
\]
If desired, precompose $F$ with a smooth based self-map of $S^{2k}$ homotopic to the identity and constant on a small neighborhood of $x_0$. We may then arrange that $\widetilde F$ is the identity on a neighborhood of the wedge
\[
(S^{2k}\times\{1\})\cup(\{x_0\}\times S^1).
\]
Thus $\widetilde F$ descends to a based map
\[
\widehat F:S^{2k}\wedge S^1 \xrightarrow{\ \cong\ }S^{2k+1}
\to\mathcal G.
\]
The map $[F]\mapsto[\widehat F]$ is exactly the usual loop-space adjunction after the homotopy equivalence $\Omega_C\mathcal G\simeq\Omega\mathcal G$ of Lemma \ref{lem:ep-based-loop-flattening}. Consequently it is an isomorphism
\[
\pi_{2k}(\Omega_C\mathcal G) \xrightarrow{\ \cong\ }\pi_{2k+1}(\mathcal G).
\]
Finally, the definition of $\beta_k^{\mathrm{even}}$ as the fiber integral of the descended odd trace form yields
\[
\int_{S^{2k}}F^*\beta_k^{\mathrm{even}} =\int_{S^{2k}\times S^1} \widetilde F^*\beta_{k+1}^{\mathrm{odd}}
=\int_{S^{2k+1}} \widehat F^*\beta_{k+1}^{\mathrm{odd}}.
\]
Here the first equality is fiber integration; in the flattened realization the descended odd trace form is literally the pullback of $\beta_{k+1}^{\mathrm{odd}}$ by $\widetilde F$. The second equality uses the quotient $S^{2k}\times S^1\to S^{2k}\wedge S^1$ and the orientation for which its image of the product fundamental class is the fundamental class of $S^{2k+1}$. The right-hand side is the usual odd pairing, which induces an isomorphism $\pi_{2k+1}(\mathcal G)\to\mathbb Z$. This proves \eqref{eqn-14} and completes the proof.
\end{proof}


\subsection{Transgression identity}

We will use the end-periodic renormalized trace to define odd forms on the full smoothing group and establish a transgression identity relating them to the even forms on $L_C\mathcal G$, analogous to \cite[Proposition 5]{mr2}.

For the end-periodic smoothing group we will use the shorthand $\mathcal{E} := G_{\mathrm{ep}}^{-\infty}(M;E)$. Let $W\subset\widetilde Y$ be the compact fundamental segment for the periodic end, and write
\[
M_N = Z\cup W_0\cup\cdots\cup W_N, \hspace{3mm} W_{j} \cong W.
\]
We recall the renormalized trace \eqref{eqn-renormalized-trace}. By Proposition \ref{renormalized-trace-lemma}, for $C\in\Psi_{\mathrm{ep}}^{-\infty}(M;E)$ with normal operator $N(C)=\widetilde C$, the renormalized trace $\rTr(C)$ converges. We extend $\rTr$ to act on $\Psi_{\mathrm{ep}}^{-\infty}(M;E)$-valued differential forms on $\mc{E}$ in the natural way. If $\eta\in\Omega^p\bigl(\mathcal E;\Psi_{\mathrm{ep}}^{-\infty}(M;E)\bigr)$, then
\[
(\rTr\eta)_g(V_1,\ldots,V_p) := \rTr \left(\eta_g(V_1,\ldots,V_p)\right).
\]
Accordingly, we obtain a map
\[
\rTr:\Omega^\bullet\bigl(\mathcal E;\Psi_{\mathrm{ep}}^{-\infty}(M;E)\bigr)\to\Omega^\bullet(\mathcal E;\mathbb C).
\]
We parametrize the FL circle by $\xi=e^{i\vartheta}$, $0\leq\vartheta\leq2\pi,$ and denote the indicial family by $C_{\vartheta} = C_{e^{i\vartheta}}$. We write
\[
\dt{C}_{\vartheta} = \frac{\partial C_{e^{i\vartheta}}}{\partial\vartheta}.
\]
For any such smooth loop $C_{\vartheta}$, define the $f$-weighted functional
\[
\mathcal L_f(C_{\vartheta}) := \frac{1}{2\pi} \int_0^{2\pi}\left( \int_W f(x)\, \tr(K_{C_{\vartheta}}(x,x))\,dx \right)d\vartheta
\]
so that the trace defect is given by
\[
\rTr[P,Q] = \mathcal L_f\bigl([P_{\vartheta},Q_{\vartheta}]
\bigr) - \frac{1}{2\pi i} \int_0^{2\pi}
\mathrm{Tr}\left(\dt P_{\vartheta}Q_{\vartheta}\right)d\vartheta.
\]
for $P,Q \in \Omega^{\bullet}\bigl(\mathcal{E};
\Psi_{\mathrm{ep}}^{-\infty}(M;E)\bigr)$. 
Here we use the graded commutator convention
\[
[P,Q] := P Q-(-1)^{\deg P\,\deg Q}Q P.
\]
Now we consider the tautological map on the end-periodic smoothing group: 
\[
a_{\mathcal{E}}: \mathcal{E} \to I + \Psi_{\mathrm{ep}}^{-\infty}(M;E), \hspace{3mm} a_{\mathcal{E}}(g) = g.
\]
The associated left Maurer-Cartan form is
\[
\Theta_{\mathcal{E}} := a_{\mathcal{E}}^{-1}da_{\mathcal{E}} \in\Omega^1(\mathcal{E};\Psi_{\mathrm{ep}}^{-\infty}(M;E)).
\]
Pointwise, one has
\[
(\Theta_{\mathcal{E}})_g(V)=g^{-1}V, \hspace{3mm}
g\in\mathcal E,\quad V\in T_g\mathcal E.
\]
This form is indeed $\Psi_{\mathrm{ep}}^{-\infty}(M;E)$-valued: if $g^{-1}=I+B$ with $B\in \Psi_{\mathrm{ep}}^{-\infty}(M;E)$, then
\[
g^{-1}V=V+BV\in \Psi_{\mathrm{ep}}^{-\infty}(M;E),
\]
since $\Psi_{\mathrm{ep}}^{-\infty}(M;E)$ is an algebra under composition.

Using the end-periodic renormalized trace, we define a new collection of odd forms:
\begin{equation}\label{gamma-forms}
\gamma_k := c_k^{\mathrm{odd}} \rTr\left(\Theta_{\mathcal{E}}^{2k-1}\right) \in \Omega^{2k-1}(\mathcal E;\mathbb C),
\end{equation}
for $k \geq 1$, where
\[
c_k^{\mathrm{odd}} = \frac{1}{(2\pi i)^k} \frac{(k-1)!}{(2k-1)!}
\]
as in \eqref{c-coeff-1}. This is the end-periodic analogue of the regularized odd Maurer-Cartan form introduced by Melrose-Rochon.

The indicial family induces the multiplicative map
\[
\mathcal I:I+\Psi_{\mathrm{ep}}^{-\infty}(M;E)
\to I+C^\infty_{\mathrm{FL}}(\mathbb{R};\mathfrak g),
\hspace{3mm} I+D\mapsto(\vartheta\mapsto I+D_\vartheta).
\]
Its restriction to $\mathcal E$ is a smooth group homomorphism
\[
\mathcal I:\mathcal E\to L_C\mathcal G.
\]
Because $g\in\mathcal E$ is invertible on $L^2(M;E)$,
its Fredholm index is zero; by the winding-number formula the
image lies in the zero-winding component $L_C^0\mathcal G$.
We also write
\[
\mathcal I_2:\Psi_{\mathrm{ep}}^{-\infty}(M;E)
\to C^\infty_{\mathrm{FL}}(\mathbb{R};\mathfrak g),
\hspace{3mm} D\mapsto(\vartheta\mapsto D_\vartheta)
\]
for the linear indicial map on smoothing operators. Note that the tautological maps and indicial maps are related by
\begin{equation}\label{eq:tautological-indicial-compatibility}
a_{\mathcal{G}} \circ \mathrm{ev}_{\vartheta} \circ \mathcal{I}
=
\mathrm{ev}_{\vartheta} \circ \mathcal{I} \circ a_{\mathcal{E}}.
\end{equation}
Recall $\omega_{\vartheta}$ and $v_{\vartheta}$ from Lemma \ref{lemma-alpha-beta-forms}. The Maurer-Cartan identities are
\[
d\omega_\vartheta+\omega_\vartheta^2=0 \text{ \, in \, }\Omega^2(L_{C}\mathcal{G};\mathfrak g),
\]
\[
\partial_\vartheta\omega_\vartheta = d v_\vartheta +\omega_\vartheta v_\vartheta-v_\vartheta\omega_\vartheta \text{ \, in \, } \Omega^1(L_{C}\mathcal{G};\mathfrak g),
\]
where $d$ here denotes the exterior differential on $L\mathcal{G}$. Consider the form $\Theta_{\mathcal{E}} = a_{\mathcal{E}}^{-1}da_{\mathcal{E}} \in\Omega^1(\mathcal{E};\Psi_{\mathrm{ep}}^{-\infty}(M;E))$. Applying $\mathcal{I}_{2}$ on the $\Psi_{\mathrm{ep}}^{-\infty}$ coefficient and evaluating at $e^{i\vartheta} \in S^{1}$ yields
\[
\Theta_{\vartheta} := \mathcal{I}_{2}(\Theta_{\mathcal{E}})(e^{i\vartheta}) \in \Omega^1(\mathcal E;\mathfrak g).
\]

\begin{lemma}\label{lem:indicial-maurer-cartan}
For every $\vartheta\in\mathbb R/2\pi\mathbb Z$, one has
\[
\mathcal{I}^{*}\omega_\vartheta = \Theta_\vartheta \text{ \, in \, } \Omega^1(\mathcal E;\mathfrak g).
\]
Here $\mathcal{I}^{*}$ is the pullback by $\mathcal{I}: \mathcal{E} \to L_{C}\mathcal{G}$.
\end{lemma}

\begin{proof}
Let $g\in\mathcal E$ and $V\in T_g\mathcal E\cong\Psi_{\mathrm{ep}}^{-\infty}(M;E)$. Multiplicativity of the indicial map gives
\[
(\Theta_\vartheta)_g(V)=\bigl(\mathcal I_2(g^{-1}V)\bigr)_\vartheta=g_\vartheta^{-1}V_\vartheta.
\]
Since the differential of $\mathcal I$ is the linear indicial map $\mathcal{I}_{2}$,
\[
(\mathcal I^*\omega_\vartheta)_g(V)=(\omega_\vartheta)_{\mathcal I(g)}(d\mathcal I_g(V))=g_\vartheta^{-1}V_\vartheta.
\]
This proves the identity.
\end{proof}

We will use a few elementary facts in the computation below. Since the trace defines a linear functional $\Tr: \mathfrak{g} \to \mathbb{C}$, naturality implies that
\[
\Tr\bigl(\Theta_\vartheta\bigr)=\Tr(\mathcal I^{*}\omega_{\vartheta})=\mathcal I^{*}\Tr(\omega_\vartheta).
\]
The same principal applies to the linear functional $\mathcal{L}_{f}$ in the trace defect formula. With this setup, we establish a transgression identity for the forms $\gamma_k$. 

\begin{prop}[Transgression identity]\label{prop:ep-regularized-odd-transgression}
For every $k\geq1$, we have
\begin{equation}\label{ep-trans-formula}
d\gamma_k = \mathcal I^{*} \left( -\beta_{k}^{\mathrm{even}}
+ dT_{k} \right),
\end{equation}
where the transgression form is defined as
\[
T_{k} := T_{k}^{(1)} + T_{k}^{(2)} \in \Omega^{2k-1}
\left( L_{C}\mathcal{G} \right),
\]
\[
T_{k}^{(1)} := \frac{c_k^{\mathrm{odd}}}{4\pi i}
\int_0^{2\pi}\mathrm{Tr}\left(\omega_{\vartheta}^{2k-1}v_{\vartheta}\right)d\vartheta,
\]
\[
T_{k}^{(2)} := c_k^{\mathrm{odd}}
\mathcal L_f \left( \omega_{\vartheta}^{2k-1} \right) = \frac{c_k^{\mathrm{odd}}}{2\pi} \int_0^{2\pi} \int_{W} f(x)\,\tr
(K_{\omega_{\vartheta}^{2k-1}}(x,x)) \,dx\,d\vartheta.
\]
In particular,
\[
\mathcal I^*[\beta_k^{\mathrm{even}}]=0 \, \text{ in } \, H_{\mathrm{dR}}^{2k}(\mathcal E).
\]
\end{prop}

\begin{proof}
Throughout the proof we will write $\Theta_{\mathcal{E}} = \Theta$ for simplicity. The Maurer-Cartan identity gives $d\Theta+\Theta^{2}=0$. For an odd power, one has
\begin{align*}
d(\Theta^{2k-1}) &= \sum_{j=0}^{2k-2} (-1)^j \Theta^j(d\Theta)\Theta^{2k-2-j}\\
&= -\left( \sum_{j=0}^{2k-2}(-1)^j \right)\Theta^{2k} \\
&= -\Theta^{2k}.
\end{align*}
Therefore
\begin{equation}\label{trans-eqn-1}
d\gamma_k = -c_k^{\mathrm{odd}}\rTr\left(\Theta^{2k}\right)
= -\frac{c_k^{\mathrm{odd}}}{2}\rTr\left([\Theta^{2k-1},\Theta]\right)
\end{equation}
because both $\Theta^{2k-1}$ and $\Theta$ have odd degree and hence
\[
[\Theta^{2k-1},\Theta]=\Theta^{2k}+ \Theta^{2k} = 2\Theta^{2k}.
\]
Now apply the trace defect formula with $P=\Theta^{2k-1}, Q=\Theta$. By Lemma \ref{lem:indicial-maurer-cartan}, the corresponding indicial forms are 
\[
P_{\vartheta} = \mathcal{I}^{*}\omega_{\vartheta}^{2k-1}, \hspace{3mm} Q_{\vartheta} = \mathcal{I}^{*}\omega_{\vartheta}.
\]
Hence $[P_{\vartheta},Q_{\vartheta}] = 2\omega_{\vartheta}^{2k}$ and therefore
\begin{equation}\label{trans-eqn-0}
\rTr\left([\Theta^{2k-1},\Theta]\right) = 2\,\mathcal I^{*} \mathcal L_f\left(\omega_{\vartheta}^{2k}\right)
- \frac{1}{2\pi i} \mathcal I^{*} \int_0^{2\pi} \mathrm{Tr}\left(\frac{\partial}{\partial\vartheta}(\omega_{\vartheta}^{2k-1})
\omega_{\vartheta} \right)d\vartheta.
\end{equation}
Let's focus on the second term on the right-hand side in \eqref{trans-eqn-0}. Integration by parts in $\vartheta$ gives
\[
\int_0^{2\pi} \mathrm{Tr}\left( \frac{\partial}{\partial\vartheta} (\omega_{\vartheta}^{2k-1}) \omega_{\vartheta} \right)d\vartheta
= -\int_0^{2\pi} \mathrm{Tr}\left(\omega_{\vartheta}^{2k-1}
\dt\omega_{\vartheta}\right)d\vartheta.
\]
Looking at the integrand on the RHS, note that
\[
\dt\omega_{\vartheta} = d v_{\vartheta} + \omega_{\vartheta}v_{\vartheta} - v_{\vartheta}\omega_{\vartheta},
\]
hence
\begin{align*}
\mathrm{Tr} \left( \omega_{\vartheta}^{2k-1} \dt\omega_{\vartheta}
\right) &= \mathrm{Tr}\left(\omega_{\vartheta}^{2k-1}d v_{\vartheta}\right)+\mathrm{Tr}\left(\omega_{\vartheta}^{2k}v_{\vartheta} \right) \\
&\quad - \mathrm{Tr}\left(\omega_{\vartheta}^{2k-1}v_{\vartheta}\omega_{\vartheta}\right). \numberthis \label{trans-eqn-2}
\end{align*}
Since $d(\omega_{\vartheta}^{2k-1})=-\omega_{\vartheta}^{2k}$, we have
\[
d\mathrm{Tr}\left(\omega_{\vartheta}^{2k-1}v_{\vartheta}\right)
= -\mathrm{Tr}\left(\omega_{\vartheta}^{2k}v_{\vartheta}
\right)-\mathrm{Tr}\left(\omega_{\vartheta}^{2k-1}d v_{\vartheta}
\right).
\]
Using this in \eqref{trans-eqn-2} yields
\begin{align*}
\mathrm{Tr} \left( \omega_{\vartheta}^{2k-1} \dt\omega_{\vartheta}
\right) &= -\mathrm{Tr}\left(\omega_{\vartheta}^{2k-1}v_{\vartheta}\omega_{\vartheta}\right) - d\Tr(\omega_{\vartheta}^{2k-1}v_{\vartheta}) \\
&= \mathrm{Tr}\left(\omega_{\vartheta}^{2k}v_{\vartheta} \right) - d\Tr(\omega_{\vartheta}^{2k-1}v_{\vartheta}),
\end{align*}
where, in the second line, swapping odd forms under the graded-cyclic trace acquires a sign. Therefore, the integral in the second term in \eqref{trans-eqn-0} is
\begin{align*}
\int_0^{2\pi} \mathrm{Tr}\left(
\frac{\partial}{\partial\vartheta}(\omega_{\vartheta}^{2k-1})
\omega_{\vartheta} \right)d\vartheta &= -\int_{0}^{2\pi}\Tr(\omega_{\vartheta}^{2k}v_{\vartheta})\,d\vartheta + d\int_{0}^{2\pi}\Tr(\omega_{\vartheta}^{2k-1}v_{\vartheta})\,d\vartheta.
\end{align*}
Pulling this back by $\mathcal{I}$ and reintroducing the coefficients then yields
\begin{align*}
&-\frac{c_k^{\mathrm{odd}}}{4\pi i}\mathcal I^{*} \int_0^{2\pi}
\mathrm{Tr}\left(\omega_{\vartheta}^{2k}v_{\vartheta}\right)d\vartheta \\
&\hspace{3mm} + d\,\mathcal I^{*}\left[\frac{c_k^{\mathrm{odd}}}{4\pi i}
\int_0^{2\pi}\mathrm{Tr}\left(\omega_{\vartheta}^{2k-1}v_{\vartheta}
\right)d\vartheta\right]
=\mathcal I^{*} \left(-\beta_{k}^{\mathrm{even}}+
dT_k^{(1)} \right). \numberthis \label{trans-eqn-4}
\end{align*}
This is the contribution to $d\gamma_{k}$ from the second term in \eqref{trans-eqn-0}.
It remains to analyze the first term in \eqref{trans-eqn-0}. The linear functional $\mathcal{L}_{f}$ commutes with exterior differentiation, so
\begin{equation}\label{trans-eqn-5}
dT_{k}^{(2)} = c_k^{\mathrm{odd}}d\mathcal L_f
\left( \omega_{\vartheta}^{2k-1} \right)
= c_k^{\mathrm{odd}} \mathcal L_f \left(d\omega_{\vartheta}^{2k-1}
\right)
= - c_k^{\mathrm{odd}} \mathcal L_f\left(\omega_{\vartheta}^{2k}\right).
\end{equation}
Combining \eqref{trans-eqn-1}, \eqref{trans-eqn-0}, \eqref{trans-eqn-4} and \eqref{trans-eqn-5} yields
\[
d\gamma_k =\mathcal I^{*}\left(-\beta_{k}^{\mathrm{even}}+dT_k^{(1)}+dT_{k}^{(2)}\right)=\mathcal I^{*}\left(
-\beta_{k}^{\mathrm{even}}+dT_{k}\right).
\]
\end{proof}

\begin{remark}
There is a sign error in \cite[(2.13)]{mr2}, where they write $+\beta_{k}^{\mathrm{even}}$ instead of minus. 
\end{remark}


In the next section we will study the long exact sequence associated with the normal map. In particular, we will show that the connecting homomorphism in even degrees
\[
\partial_{2k}:\pi_{2k}(L_C^0\mathcal G)
\to \pi_{2k-1} \bigl(\dt{G}_{\mathrm{ep}}^{-\infty}(M;E)\bigr)
\]
is an isomorphism. The transgression identity in Proposition \ref{prop:ep-regularized-odd-transgression} will be a key tool for that step.

\section{Contractibility theorem}\label{sec-4}

The weak contractibility theorem for the cusp smoothing group in \cite[Proposition 6]{mr2} has a different structure in the end-periodic setting, because the FL transform of the normal operator is not an ordinary smooth loop of smoothing operators. It is a smooth family satisfying the twisting condition \eqref{eqn- sus-FL-operator-twist}. Consequently, the correct normal smoothing group is the twisted loop group $L_C\mathcal G$ introduced in \eqref{eq:ep-twisted-loop-group}, not the ordinary free loop group $C^\infty(S^1;\mathcal G)$. They have the same homotopy groups, as proved in the preceding section, but their evaluation sequences are not identified by the ordinary constant-loop splitting.

To prove the theorem, we first show that the normal homomorphism on $G_{\mathrm{ep}}^{-\infty}$ is a locally trivial principal bundle with structure group $\dt{G}_{\mathrm{ep}}^{-\infty}$. The transgression identity of Proposition \ref{prop:ep-regularized-odd-transgression} then allows us to compute the even-degree connecting maps in the homotopy long exact sequence. We deduce that evaluation of the indicial family at $\theta=0$ is a weak homotopy equivalence and that its kernel is weakly contractible.

The principal bundle structure for the normal map $N$ reflects the fact that our end-periodic smoothing algebra is somewhat simplified in comparison with the cusp smoothing algebra in \cite{mr2}. The local splitting for $N$ is possible because choosing any end cutoff $\chi$ determines an extension map $\mathrm{Ext}_{\chi}$ which is a continuous linear right inverse $N$, cf. \eqref{eq:Ext-chi-topology}. That splitting is not possible for the cusp normal map, because the latter encodes more information (a full Taylor series at the cusp front face). Hence the statement in \cite[Lemma 1]{mr2} asserts that the cusp normal map is merely a Serre fibration rather than a principal bundle.

\subsection{The normal map as a principal bundle}
\label{subsec:ep-normal-fibration}

We adopt the shorthand:
\[
\begin{aligned}
\mathcal G&:=G^{-\infty}(Y;E_Y),
&\mathcal E&:=G_{\mathrm{ep}}^{-\infty}(M;E),\\
\mathcal K&:=\dt G_{\mathrm{ep}}^{-\infty}(M;E),
&L_C\mathcal G&:=G_{\mathbb Z}^{-\infty}(Y;E_Y).
\end{aligned}
\]
Thus $L_C\mathcal G$ consists of families $a_\theta=I+A_\theta$ with
\[
A_{(\cdot)}\in C^\infty_{\mathrm{FL}}
(\mathbb R;\Psi^{-\infty}(Y;E_Y)),\hspace{3mm}
a_{\theta+2\pi}=Ca_\theta C^{-1},
\]
and each $a_\theta$ invertible. Write
\begin{equation}\label{eq:ep-contractibility-base}
L_{C}^{0}\mathcal{G} := \{a\in L_C\mathcal G:\operatorname{wind}(a)=0\},
\end{equation}
where the winding homomorphism is the one in
\eqref{eq:ep-twisted-winding}. By Lemma \ref{lem:ep-twisted-determinant-winding}, $L_{C}^{0}\mathcal{G}$ is precisely the identity component of $L_C\mathcal G$.

For $g=I+A\in\mathcal E$, let $\widetilde A=N(A)$ be its periodic normal operator and put
\begin{equation}\label{eq:ep-contractibility-indicial}
\mathcal I(g):=I+\mathcal F\widetilde A\mathcal F^{-1} =(I+A_\theta)_{\theta\in\mathbb R},
\end{equation}
acting fiberwise on the FL Hilbert bundle
$\mathcal H^{\mathrm{FL}}\to S^1$.

\begin{prop}[Exact normal sequence]
\label{prop:ep-exact-normal-sequence}
The indicial homomorphism fits into the short exact sequence of topological groups
\begin{equation}\label{zero-index-ses}
1\to\mathcal K\to\mathcal E \xrightarrow{\ \mathcal I\ }L_{C}^{0}\mathcal{G}\to1.
\end{equation}
\end{prop}

\begin{proof}
The normal map is multiplicative by Lemma \ref{lem:normal-op-homom}, and conjugation by $\mathcal{F}$ is a $*$-algebra isomorphism by Theorem \ref{FL-transform-thm}. Therefore $\mathcal I$ is a group homomorphism. Its continuity
follows from the continuity of the normal map in the direct-sum Fr\'echet topology and from Fr\'echet continuity, cf. Theorem \ref{FL-transform-thm}.

If $g=I+A\in\mathcal E$, then $g$ is invertible on $L^2(M;E)$, so its Fredholm index is zero. The winding-number formula (Proposition \ref{prop-winding}) gives
\[
0=\operatorname{ind}(g)=-\operatorname{wind}(\mathcal I(g)).
\]
Thus $\mathcal I(g)\in L_{C}^{0}\mathcal{G}$. Moreover, $\mathcal I(g)=I$ if and only if $N(A)=0$: this is the injectivity of the FL transform. By the exact sequence in Lemma \ref{ses-lemma-1}, the latter condition is equivalent to $A\in\dt\Psi_{\mathrm{ep}}^{-\infty}(M;E)$. Since $g$ is invertible, this is precisely the defining condition for $g\in\mathcal K$; hence $\ker\mathcal I=\mathcal K$.

For surjectivity, let $a=(I+A_\theta)\in L_{C}^{0}\mathcal{G}$. By definition, $A_{(\cdot)}$ belongs to $C_{\mathrm{FL}}^{\infty}(\mathbb{R};\Psi^{-\infty}(Y;E_{Y}))$. Theorem \ref{FL-transform-thm} therefore supplies a unique
\[
\widetilde A=\mathcal F^{-1}A_{(\cdot)}\mathcal F
\in\Psi_{\mathbb Z}^{-\infty}(\widetilde Y;\widetilde E).
\]
Choose an end cutoff $\chi$ as in \eqref{eq:Ext-chi-topology} and set $A=\operatorname{Ext}_\chi(\widetilde A)=\chi\widetilde A\chi$. Since each $I+A_\theta$ is invertible, the Fredholm criterion (Proposition \ref{ep-theory-prop-1}) makes $I+A$ Fredholm. Furthermore,
\[
\operatorname{ind}(I+A) =-\operatorname{wind}(a)=0.
\]
Hence Proposition \ref{prop:ep-reduced-index} gives $R\in\dt\Psi_{\mathrm{ep}}^{-\infty}(M;E)$ such that $I+A+R$ is invertible. The inverse-closedness lemma \ref{lem:ep-inverse-closed} shows that its inverse is again an end-periodic smoothing perturbation of $I$, so $g:=I+A+R\in\mathcal E$. Finally $N(R)=0$ and $N\circ\operatorname{Ext}_\chi=\mathrm{Id}$; therefore $\mathcal I(g)=a$. This proves surjectivity.
\end{proof}

The linear right inverse $\operatorname{Ext}_\chi$ also gives local sections at the level of invertible operators. If
$\widetilde A\in\Psi_{\mathbb Z}^{-\infty}$ and
$A_{(\cdot)}=\mathcal F\widetilde A\mathcal F^{-1}$, then the operator norms satisfy
\begin{equation}\label{eq:Ext-operator-bound}
\begin{aligned}
\|\operatorname{Ext}_\chi(\widetilde A)\|_{\mathcal B(L^2(M;E))}
&\leq\|\chi\|_\infty^2 \|\widetilde A\|_{\mathcal{B}(L^2(\widetilde Y;\widetilde E))}\\
&=\|\chi\|_\infty^2 \sup_{0\leq\theta\leq2\pi}
\|A_\theta\|_{\mathcal B(L^2(Y;E_Y))}.
\end{aligned}
\end{equation}
The identity in the second line is the norm formula for a decomposable operator on $L^2(S^1;\mathcal H^{\mathrm{FL}})$, using the unitary FL transform \eqref{eqn-sus-FL-unitary}. The supremum is finite and is continuous in the twisted-family Fr\'echet topology:
\[
\sup_\theta\|A_\theta\|_{\mathcal B(L^2(Y;E_Y))}
\leq\operatorname{vol}(Y)\,P_{0,0}(A).
\]

\begin{lemma}[Local lifting near the identity]
\label{lem:local-section-identity}
There are an open neighborhood $\mathcal V_I\subset L_{C}^{0}\mathcal{G}$ of the identity family and a smooth map $s_I:\mathcal V_I\to\mathcal E$ satisfying $\mathcal I\circ s_I=\mathrm{Id}_{\mathcal V_I}$.
\end{lemma}

\begin{proof}
Choose the cutoff $\chi$ with $0\leq\chi\leq1$, and put
\[
\mathcal V_I:=\left\{a=(I+A_\theta)\in L_{C}^{0}\mathcal{G}: \sup_{0\leq\theta\leq2\pi}
\|A_\theta\|_{\mathcal B(L^2(Y;E_Y))}<1\right\}.
\]
The supremum is continuous by the estimate preceding the lemma; hence $\mathcal V_I$ is an open neighborhood of the identity in the twisted-loop topology. For $a\in\mathcal V_I$ set
\begin{equation}\label{eq:ep-contractibility-local-identity}
s_I(a):=I+\operatorname{Ext}_\chi
\bigl(\mathcal F^{-1}(a-I)\mathcal F\bigr).
\end{equation}
The parenthesized inverse FL transform exists in the periodic smoothing algebra because $a-I$ satisfies the twisting \eqref{eqn- sus-FL-operator-twist}. Estimate \eqref{eq:Ext-operator-bound} yields $\|s_I(a)-I\|_{\mathcal B(L^2(M;E))}<1$, so $s_I(a)$ is invertible by the Neumann lemma. Inverse-closedness, Lemma \ref{lem:ep-inverse-closed}, shows that $s_I(a)\in\mathcal E$. Because
$N\circ\operatorname{Ext}_\chi=\mathrm{Id}$, its indicial family is $a$ and thus $\mathcal I(s_I(a))=a$. Finally the expression \eqref{eq:ep-contractibility-local-identity} is the restriction of an affine continuous linear map between the defining Fr\'echet spaces. In particular it is smooth.
\end{proof}

\begin{prop}[Normal fibration]\label{prop:ep-serre-fibration}
The homomorphism $\mathcal I:\mathcal E\to L_{C}^{0}\mathcal{G}$ is a locally trivial principal $\mathcal K$-bundle. Its local sections are smooth in the natural Fr\'echet group charts; in particular it is a Serre fibration.
\end{prop}

\begin{proof}
Fix $a_0\in L_{C}^{0}\mathcal{G}$. By
Proposition \ref{prop:ep-exact-normal-sequence} choose $g_0\in\mathcal E$ with $\mathcal I(g_0)=a_0$. Let $\mathcal V_I$ and $s_I$ be given by Lemma \ref{lem:local-section-identity}. The neighborhood $U_{a_0}:=a_0\mathcal V_I$ is open in $L_{C}^{0}\mathcal{G}$,
and 
\begin{equation}\label{eq:ep-contractibility-local-lift}
s_{a_0}(a):=g_0\,s_I(a_0^{-1}a),\hspace{3mm} a\in U_{a_0},
\end{equation}
is a local section, since $\mathcal I(s_{a_0}(a))=a_0(a_0^{-1}a)=a$. Notice that multiplication $a_0^{-1}a$ is pointwise multiplication of \emph{twisted} families: it again satisfies
$(a_0^{-1}a)_{\theta+2\pi} =C(a_0^{-1}a)_\theta C^{-1}$.

For any such pair $(U,s)$, define \begin{equation}\label{eq:ep-contractibility-trivialization}
\begin{aligned}
\Phi_s:\mathcal I^{-1}(U)&\to U\times\mathcal K,\\
g&\mapsto \bigl(\mathcal I(g),s(\mathcal I(g))^{-1}g\bigr),\\
\Phi_s^{-1}:U\times\mathcal K&\to\mathcal I^{-1}(U),\\
(a,r)&\mapsto s(a)r.
\end{aligned}
\end{equation}
The second component of $\Phi_s(g)$ lies in $\mathcal K$, since
\[
\mathcal I\bigl(s(\mathcal I(g))^{-1}g\bigr)
=\mathcal I(g)^{-1}\mathcal I(g)=I.
\]
The formulas show directly that the two maps are mutually inverse, and that the right action of $\mathcal K$ on $\mathcal E$ by multiplication is free and transitive on each fiber. They are continuous, and smooth in the Fr\'echet group structures provided by the smoothing-algebra operations and inverse-closedness. Consequently $\mathcal I$ is a locally trivial principal bundle. Every locally trivial fiber bundle is a Serre fibration.
\end{proof}

\begin{remark}[Comparison with the cusp normal sequence]
It is insightful to compare with the cusp normal sequence in \cite[(2.14)]{mr2}. We have used a special feature of the end-periodic normal map which is the splitting $N\circ\operatorname{Ext}_\chi=\mathrm{Id}$. This permits the explicit local section \eqref{eq:ep-contractibility-local-identity}. By contrast, the normal map \cite[(2.14)]{mr2} retains higher-order Taylor coefficients at the cusp face. Hence the cusp normal map does not admit this sort of splitting, and \cite[Lemma 1]{mr2} shows that it defines a Serre fibration rather than a locally trivial principal bundle.
\end{remark}

\subsection{Homotopy groups and connecting maps}
\label{subsec:ep-boundary-map}

Since \eqref{zero-index-ses} is a Serre fibration, it gives the homotopy long exact sequence, based at the identity,
\begin{equation}\label{eqn-homotopy-les}
\cdots\to\pi_j(\mathcal K) \to\pi_j(\mathcal E)
\xrightarrow{\mathcal I_*}\pi_j(L_{C}^{0}\mathcal{G}) \xrightarrow{\partial_j}\pi_{j-1}(\mathcal K) \to\cdots .
\end{equation}
We first identify its outer terms and the connecting maps.

For $g=I+R\in\mathcal K$, inversion in the residual group and the ideal property of $\dt{\Psi}_{\mathrm{ep}}^{-\infty}(M;E)$ imply that its left Maurer-Cartan
form takes values in the residual algebra:
\[
\Theta_{\mathcal K}:=g^{-1}dg \in\Omega^1(\mathcal K;\dt{\Psi}_{\mathrm{ep}}^{-\infty}(M;E)).
\]
By Lemma \ref{lem:ep-residual-trace-class}, residual operators are trace-class, so we can take the ordinary trace. The smoothness of inversion and multiplication in $\mathcal{K}$ (cf. Remark \ref{remark-smooth-inversion}) shows that
\begin{equation}\label{eq:ep-residual-odd-form}
\beta_k^{\mathrm{odd}} :=c_k^{\mathrm{odd}}\Tr \bigl(\Theta_{\mathcal K}^{2k-1}\bigr), \hspace{3mm}
c_k^{\mathrm{odd}}=\frac{(k-1)!}{(2\pi i)^k(2k-1)!},
\end{equation}
is a smooth differential form on $\mathcal K$ for every $k\geq1$.

\begin{lemma}[Homotopy type for the residual group]
\label{lem:ep-residual-bott}
The inclusion of the stable complex general linear group, realized by finite-rank perturbations with compactly supported smooth kernels, into $\mathcal K$ is a weak homotopy equivalence. Consequently
\begin{equation}\label{eq:ep-H-bott}
\pi_{2k}(\mathcal K)=0,\hspace{3mm} \pi_{2k-1}(\mathcal K)\cong\mathbb Z, \hspace{3mm} k\geq 1.
\end{equation}
Moreover, we have an isomorphism given by
\begin{equation}\label{eq:ep-residual-bott-iso}
\pi_{2k-1}(\mathcal K) \xrightarrow{\ \cong\ }\mathbb Z, \hspace{3mm} [u] \mapsto\int_{S^{2k-1}}u^*\beta_k^{\mathrm{odd}},
\end{equation}
where $\beta_{k}^{\mathrm{odd}}$ is given by \eqref{eq:ep-residual-odd-form}.
\end{lemma}

\begin{proof}
This follows by the same finite-rank approximation argument as outlined in the proof of \cite[Proposition 3]{mr2}.
\end{proof}

The corresponding information about $L_{C}^{0}\mathcal{G}$ has already been proved in the preceding section. Lemma \ref{lem:ep-twisted-homotopy-groups} gives
\begin{equation}\label{eq:ep-twisted-normal-bott}
\begin{aligned}
&\pi_{2k-1}(L_{C}^{0}\mathcal{G})\cong\mathbb Z,
&&(\operatorname{ev}_1)_*:
\pi_{2k-1}(L_{C}^{0}\mathcal{G})
\xrightarrow{\ \cong\ }\pi_{2k-1}(\mathcal G),\\
&\pi_{2k}(L_{C}^{0}\mathcal{G})\cong\mathbb Z,
&&(\iota_C)_*:
\pi_{2k}(\Omega_C^0\mathcal G)
\xrightarrow{\ \cong\ }\pi_{2k}(L_{C}^{0}\mathcal{G}),
\end{aligned}\hspace{3mm} k\geq1,
\end{equation}
where
\begin{equation}\label{eq:ep-contractibility-based-normal}
\Omega_C^0\mathcal G:=\Omega_C\mathcal G\cap L_{C}^{0}\mathcal{G} =\{a\in L_{C}^{0}\mathcal{G}:a_0=I\}
\end{equation}
is the identity component of the \emph{twisted} based-loop group. Here $\iota_C$ denotes its inclusion in $L_{C}^{0}\mathcal{G}$. The identification $\pi_j(\Omega_C^0\mathcal G) \cong\pi_{j+1}(\mathcal G)$ for $j\geq1$ follows from the endpoint-flattening lemma \ref{lem:ep-based-loop-flattening}.

We use the odd and even forms of the preceding section to specify these isomorphisms. Recall the even form $\beta_k^{\mathrm{even}}$ on $L_C\mathcal G$ from \eqref{eq:beta-even-explicit}. Proposition \ref{prop:ep-free-loop-basic-classes} restricted to the zero-winding component gives an isomorphism
\begin{equation}\label{eq:ep-normalized-pairings}
\pi_{2k}(L_{C}^{0}\mathcal{G}) \xrightarrow{\ \cong\ }\mathbb Z, \hspace{3mm} [F]\mapsto \int_{S^{2k}}F^*\beta_k^{\mathrm{even}}.
\end{equation}

\noindent Finally, write $\Theta_{\mathcal E}=g^{-1}dg$ for the Maurer-Cartan form and set
\[
\gamma_k:=c_k^{\mathrm{odd}}\rTr(\Theta_{\mathcal E}^{2k-1}) \in\Omega^{2k-1}(\mathcal E;\mathbb C).
\]
The renormalized trace agrees with the ordinary trace on residual operators, because their normal operator vanishes. Hence
\begin{equation}\label{eq:ep-gamma-restriction}
\left.\gamma_k\right|_{\mathcal K} =\beta_k^{\mathrm{odd}}.
\end{equation}
By Proposition \ref{prop:ep-regularized-odd-transgression},
\begin{equation}\label{eq:ep-transgression-recalled} 
d\gamma_k =\mathcal I^*\bigl(-\beta_k^{\mathrm{even}}+dT_k\bigr),
\end{equation}
where $T_k\in\Omega^{2k-1}(L_C\mathcal G)$ is the explicitly defined correction form in \eqref{ep-trans-formula}. Its restriction to $L_{C}^{0}\mathcal{G}$ will be understood in what follows.

\begin{prop}[Even-degree connecting maps]
\label{prop:ep-even-boundary-map}
For every $k\geq1$ the connecting homomorphism
\[
\partial_{2k}:\pi_{2k}(L_{C}^{0}\mathcal{G}) \to\pi_{2k-1}(\mathcal K)
\]
of \eqref{eqn-homotopy-les} is an isomorphism. Using the boundary orientations chosen below, it is multiplication by $-1$.
\end{prop}

\begin{proof}
Let $F:(S^{2k},q)\to(L_{C}^{0}\mathcal{G},I)$ be a smooth based representative of a generator, oriented so that
\begin{equation}\label{eq:ep-contractibility-generator}
\int_{S^{2k}}F^*\beta_k^{\mathrm{even}}=1.
\end{equation}
The real oriented blow-up $D=[S^{2k};\{q\}]$ of the point $q$ is diffeomorphic to a closed $2k$-ball. Its blow-down $b:D\to S^{2k}$ is smooth, restricts to an orientation-preserving diffeomorphism of interiors, and maps $\partial D$ to $q$. In particular, the pullback of a smooth top-degree form on $S^{2k}$ has the same integral over $D$ as the original form over $S^{2k}$.

Pull back the principal bundle of Proposition \ref{prop:ep-serre-fibration} along $F\circ b:D\to L_{C}^{0}\mathcal{G}$. Since $D$ is contractible, this principal $\mathcal K$-bundle is trivial. A smooth trivialization provides a smooth lift $\widetilde F:D\to\mathcal E$ with
\begin{equation}\label{eq:ep-lift-blowdown}
\mathcal I\circ\widetilde F=F\circ b.
\end{equation}
The restriction of $\mathcal{I}\circ\widetilde F$ to $\partial D$ is constant. Choose $s_0\in\partial D$ and set
\[
u(s):=\widetilde F(s_0)^{-1}\widetilde F(s),
\hspace{3mm} s\in\partial D\cong S^{2k-1}.
\]
Then $\mathcal{I}(u(s))=I$ and $u(s_0)=I$, so $u$ is a based map into $\mathcal K$. By the standard construction of the connecting homomorphism for a principal bundle (see, for example, \cite[Theorem 4.41]{hatcher}), we have $[u]=\partial_{2k}[F]$.

On the boundary, $\widetilde F(s)=\widetilde F(s_0)u(s)$ with $\widetilde F(s_0)$ constant. Therefore the left Maurer-Cartan forms agree on the boundary:
\[
\left.\widetilde F^*\Theta_{\mathcal E}\right|_{\partial D} =u^{-1}du.
\]
Using \eqref{eq:ep-gamma-restriction} then gives
\begin{equation}\label{eq:ep-boundary-form-identity}
\int_{\partial D}\widetilde F^*\gamma_k =\int_{S^{2k-1}}u^*\beta_k^{\mathrm{odd}}.
\end{equation}
Stokes' theorem, \eqref{eq:ep-transgression-recalled}, and \eqref{eq:ep-lift-blowdown} now imply
\[
\begin{aligned}
\int_{S^{2k-1}}u^*\beta_k^{\mathrm{odd}}
&=\int_D d(\widetilde F^*\gamma_k)\\
&=\int_D(F\circ b)^*\bigl(-\beta_k^{\mathrm{even}}+dT_k\bigr).
\end{aligned}
\]
The exact term contributes zero because $F\circ b$ is constant on $\partial D$, hence
\[
\int_D(F\circ b)^*dT_k =\int_{\partial D}(F\circ b)^*T_k=0.
\]
Then by the degree-one property of the blow-down and
\eqref{eq:ep-contractibility-generator},
\[
\int_{S^{2k-1}}u^*\beta_k^{\mathrm{odd}}
=-\int_D(F\circ b)^*\beta_k^{\mathrm{even}}
=-\int_{S^{2k}}F^*\beta_k^{\mathrm{even}}=-1.
\]
Thus $[u]$ is a generator of $\pi_{2k-1}(\mathcal K)\cong\mathbb Z$ (cf. Lemma \ref{lem:ep-residual-bott}). Since $\pi_{2k}(L_{C}^{0}\mathcal G)\cong\mathbb Z$, the connecting map is an isomorphism.
\end{proof}

\subsection{The homotopy type of the end-periodic smoothing group}
\label{subsec:ep-contractibility-main}

\begin{thm}[Homotopy type of the full end-periodic smoothing group]\label{thm:homotopy-type-Gep}
Evaluation at $\theta=0$ of the \emph{twisted} indicial family,
\begin{equation}\label{eq:ep-evaluation}
\begin{aligned}
\varepsilon:\mathcal E&\to\mathcal G,\\
g&\mapsto (\mathrm{ev}_{1}\circ \mathcal{I})(g)
\end{aligned}
\end{equation}
is a weak homotopy equivalence. In particular,
$\mathcal E$ is path connected and
\[
\pi_{2k}(\mathcal E)=0,\hspace{3mm} \pi_{2k-1}(\mathcal E)\cong\mathbb Z, \hspace{3mm} k\geq1.
\]
Thus the full group $G_{\mathrm{ep}}^{-\infty}(M;E)$ is not weakly contractible.
\end{thm}

\begin{proof}
We first calculate the homotopy groups of $\mathcal E$ using \eqref{eqn-homotopy-les}. For each $k\geq1$, its even-degree segment includes
\[
0=\pi_{2k}(\mathcal K)\to \pi_{2k}(\mathcal E)
\xrightarrow{\mathcal I_*}\pi_{2k}(L_{C}^{0}\mathcal{G}) \xrightarrow{\partial_{2k}}\pi_{2k-1}(\mathcal K).
\]
The last arrow is an isomorphism by Proposition \ref{prop:ep-even-boundary-map}. Exactness makes $\mathcal I_*$ injective with zero image, so $\pi_{2k}(\mathcal E)=0$.

In odd degree the relevant part is
\[
\pi_{2k}(L_{C}^{0}\mathcal{G}) \xrightarrow{\partial_{2k}}\pi_{2k-1}(\mathcal K)\to\pi_{2k-1}(\mathcal E)\xrightarrow{\mathcal I_*}\pi_{2k-1}(L_{C}^{0}\mathcal{G})\to\pi_{2k-2}(\mathcal K).
\]
The first map is surjective and the last group vanishes. For $k=1$, the last group means $\pi_0(\mathcal K)$, which is trivial by Lemma \ref{lem:ep-residual-bott}. Exactness therefore shows that
\begin{equation}\label{eq:ep-contractibility-odd-normal-iso}
\mathcal I_*:\pi_{2k-1}(\mathcal E)\xrightarrow{\ \cong\ }\pi_{2k-1}(L_{C}^{0}\mathcal{G}).
\end{equation}
The path-component part of the long exact sequence, together with $\pi_0(\mathcal K)=\pi_0(L_{C}^{0}\mathcal{G})=0$, gives $\pi_0(\mathcal E)=0$.

It remains to prove that the \emph{evaluation map itself} induces the stated isomorphisms. On the normal group, the preceding section's twisted evaluation lemma
\ref{lem:ep-twisted-homotopy-groups} gives
\begin{equation}\label{eq:ep-contractibility-odd-evaluation}
(\operatorname{ev}_1)_*: \pi_{2k-1}(L_{C}^{0}\mathcal{G}) \xrightarrow{\ \cong\ }\pi_{2k-1}(\mathcal G).
\end{equation}
Here $\operatorname{ev}_1(a)=a_0$; the subscript refers to the point $\xi=1=e^{i0}$. Because $\varepsilon=\operatorname{ev}_1\circ\mathcal I$, \eqref{eq:ep-contractibility-odd-normal-iso} and \eqref{eq:ep-contractibility-odd-evaluation} imply that $\varepsilon_*$ is an isomorphism in every odd degree. In even positive degrees both $\pi_{2k}(\mathcal E)$ and $\pi_{2k}(\mathcal G)$ vanish. Finally both groups are connected. Therefore $\varepsilon$ induces an isomorphism on every homotopy group, as claimed.
\end{proof}

Define the subgroup of elements whose indicial family is the identity at the chosen basepoint $\xi=1$ by
\begin{equation}\label{eq:ep-based-group}
\mathcal E_*:=G_{\mathrm{ep},*}^{-\infty}(M;E)
:=\{g\in\mathcal E:(\mathcal I(g))_0=I\}=\ker\varepsilon.
\end{equation}
Its normal homomorphism has as its image the group
$\Omega_C^0\mathcal G$ of \eqref{eq:ep-contractibility-based-normal}, not the ordinary based-loop group. Restricting the exact normal sequence therefore gives
\begin{equation}\label{eq:ep-based-fibration}
1\to\mathcal K\to\mathcal E_*
\xrightarrow{\ \mathcal I\ }\Omega_C^0\mathcal G
\to1.
\end{equation}
Indeed, if $a\in\Omega_C^0\mathcal G$, surjectivity of $\mathcal I:\mathcal E\to L_{C}^{0}\mathcal{G}$ supplies $g$ with $\mathcal I(g)=a$, and the condition $a_0=I$ automatically gives $g\in\mathcal E_*$. Conversely, the normal family of an element of $\mathcal E_*$ satisfies precisely these two conditions. The sequence \eqref{eq:ep-based-fibration} is the pullback of the principal bundle \eqref{zero-index-ses} under the inclusion $\Omega_C^0\mathcal G\hookrightarrow L_{C}^{0}\mathcal{G}$, hence is itself a principal $\mathcal K$-bundle and a Serre fibration.

\begin{cor}[Based end-periodic weak contractibility]
\label{cor:based-ep-weak-contractibility}
The based end-periodic smoothing group $G_{\mathrm{ep},*}^{-\infty}(M;E)$ is weakly contractible:
\[
\pi_j\bigl(G_{\mathrm{ep},*}^{-\infty}(M;E)\bigr)=0
\hspace{3mm}\text{for every }j\geq0.
\]
\end{cor}

\begin{proof}
By the flattening lemma \ref{lem:ep-based-loop-flattening} and Bott periodicity,
\[
\begin{aligned}
\pi_{2k}(\Omega_C^0\mathcal G)&\cong
\pi_{2k+1}(\mathcal G)\cong\mathbb Z, && k\geq1,\\
\pi_{2k-1}(\Omega_C^0\mathcal G)&\cong
\pi_{2k}(\mathcal G)=0, &&k\geq1,\\ \pi_0(\Omega_C^0\mathcal G)&=0.
\end{aligned}
\]
The inclusion $\iota_C:\Omega_C^0\mathcal G\hookrightarrow L_{C}^{0}\mathcal{G}$ induces an isomorphism on every positive even homotopy group by Lemma \ref{lem:ep-twisted-homotopy-groups}. Naturality of the long exact sequence for the pullback fibration therefore identifies the even connecting map for \eqref{eq:ep-based-fibration} with
\[
\pi_{2k}(\Omega_C^0\mathcal G) \xrightarrow{(\iota_C)_*}\pi_{2k}(L_{C}^{0}\mathcal{G}) \xrightarrow{\partial_{2k}}\pi_{2k-1}(\mathcal K).
\]
Both arrows are isomorphisms, the second by Proposition \ref{prop:ep-even-boundary-map}. Thus the based connecting map is an isomorphism for every $k\geq1$.

Now the long exact sequence of \eqref{eq:ep-based-fibration}, combined with $\pi_{2k}(\mathcal K)=0$ and $\pi_{2k-1}(\Omega_C^0\mathcal G)=0$, gives $\pi_j(\mathcal E_*)=0$ for every $j\geq1$. In degree one, exactness uses $\pi_0(\mathcal K)=0$, which was proved in Lemma \ref{lem:ep-residual-bott}. Both $\mathcal K$ and $\Omega_C^0\mathcal G$ are path connected, so the sequence also gives $\pi_0(\mathcal E_*)=0$. This proves weak contractibility.
\end{proof}

In forthcoming work we will use Corollary \ref{cor:based-ep-weak-contractibility} to construct invertible perturbations for families of elliptic end-periodic differential operators. If $B$ is a compact parameter space, and $P = (P_{b})_{b \in B}$ is a smooth family of elliptic end-periodic differential operators, one can construct a family of end-periodic smoothing operators $Q = (Q_{b})_{b \in B}$ such that $P+Q$ is invertible for every $b \in B$, provided that the necessary $K$-theoretic obstruction vanishes. This mirrors \cite[Theorem 1]{mr2}. Here the weakly contractible group $G_{\mathrm{ep},*}^{-\infty}$ plays a role analogous to the weakly contractible cusp smoothing group in \cite{mr2}.

\appendix

\section{Convolution estimates and composition properties}\label{appendix-1}

The basic composition properties for end-periodic smoothing operators follow from some elementary Schwartz convolution estimates which we summarize in the following lemma.

\begin{lemma}[Schwartz convolution estimates]
\label{schwartz-convolution-lemma}
Let $A=(a_{jk})$, $B=(b_{jk})$, and $R=(r_{jk})$ be nonnegative matrices. Assume that $A$ and $B$ have relative Schwartz decay, in the sense that, for every $N\geq 0$, there is a constant $C_N>0$ such that
\[
a_{jk}\leq C_N\langle j-k\rangle^{-N}, \hspace{3mm}
b_{jk}\leq C_N\langle j-k\rangle^{-N}.
\]
Assume that $R$ has absolute Schwartz decay, in the sense that, for every $N_1,N_2\geq 0$, there is a constant $C_{N_1,N_2}>0$ such that
\[
r_{jk} \leq C_{N_1,N_2} \langle j\rangle^{-N_1}
\langle k\rangle^{-N_2}.
\]
Then the following estimates hold.

\begin{enumerate}
\item The matrix
\[
(AR)_{j\ell}:=\sum_{k\geq 0}a_{jk}r_{k\ell},
\hspace{3mm} j,\ell\geq 0,
\]
has absolute Schwartz decay. Equivalently, for every $N_1,N_2\geq 0$,
\[
\sup_{j,\ell\geq 0} \langle j\rangle^{N_1} \langle \ell\rangle^{N_2} |(AR)_{j\ell}| <\infty.
\]

\item The matrix
\[
(RB)_{j\ell}:=\sum_{k\geq 0}r_{jk}b_{k\ell},\hspace{3mm} j,\ell\geq 0,
\]
has absolute Schwartz decay. Equivalently, for every $N_1,N_2\geq 0$,
\[
\sup_{j,\ell\geq 0} \langle j\rangle^{N_1}
\langle \ell\rangle^{N_2} |(RB)_{j\ell}| <\infty.
\]

\item Let $q=(q_k)_{k\in\mathbb Z}$ be a bounded sequence such that $q_k=0$ for all sufficiently large positive $k$. Then
\[
c_{j\ell}:= \sum_{k\in\mathbb Z} a_{jk}q_kb_{k\ell},
\hspace{3mm} j,\ell\geq 0,
\]
has absolute Schwartz decay in $(j,\ell)$. That is, for every $N_1,N_2\geq 0$,
\[
\sup_{j,\ell\geq 0} \langle j\rangle^{N_1} \langle \ell\rangle^{N_2} |c_{j\ell}| <\infty.
\]
\end{enumerate}
\end{lemma}

\begin{proof}
We prove the first assertion. Fix $N_1,N_2\geq 0$, and choose integers $(L,P,Q$ with $L>N_1+2$, $P>N_1+2$, and $Q>N_2$. Using the assumed estimates for $A$ and $R$, we have
\begin{align*}
\langle j\rangle^{N_1} \langle \ell\rangle^{N_2} (AR)_{j\ell} &\leq C\sum_{k\geq 0} \langle j \rangle^{N_{1}} \langle j-k\rangle^{-L} \langle\ell\rangle^{N_{2}-Q} \\
&\leq C\sum_{k\geq 0}\langle j-k\rangle^{N_1-L}\langle k\rangle^{N_1-P}\langle \ell\rangle^{N_2-Q}.
\end{align*}
where in the second line we used $\langle j\rangle \leq C\langle j-k\rangle\langle k\rangle$ for some $C \geq 0$. Since $Q>N_2$, the last factor is uniformly bounded by $1$. Moreover, since $L>N_1+2$ and $P>N_1+2$, the sum $\sum_{k\geq 0}\langle j-k\rangle^{N_1-L}\langle k\rangle^{N_1-P}$ is bounded independently of $j$. Thus
\[
\sup_{j,\ell\geq 0} \langle j\rangle^{N_1} \langle \ell\rangle^{N_2}(AR)_{j\ell} <\infty,
\]
as claimed. The second assertion is proved in the same way, but with the roles of the left and right variables interchanged.

It remains to prove the third assertion. Let $K\in\mathbb Z$ be such that $q_k=0$ for all $k>K$. Since $q$ is bounded, it suffices to estimate
\[
\sum_{k\leq K} \langle j-k\rangle^{-L} \langle k-\ell\rangle^{-L'}
\]
with $L,L'$ chosen arbitrarily large. For $j,\ell\geq 0$ and $k\leq K$, we have
\[
\langle j\rangle \leq C_K\langle j-k\rangle,\hspace{3mm} \langle \ell\rangle \leq C_K\langle \ell-k\rangle,
\]
where we use $\langle j\rangle \leq C\langle j-k\rangle\langle k\rangle$. Hence, 
\[
\langle j\rangle^{N_1} \langle \ell\rangle^{N_2} |c_{j\ell}| \leq C \sum_{k\leq K} \langle j-k\rangle^{N_1-L} \langle \ell-k\rangle^{N_2-L'}.
\]
Choose $L>N_1+2$ and $L'>N_2+2$. Then the right-hand side is uniformly bounded in $j,\ell\geq 0$. Therefore $c_{j\ell}$ has absolute Schwartz decay in $(j,\ell)$.
\end{proof}

\begin{proof}[Proof of Lemma \ref{composition-lemma}]
We explain how the estimates required for these composition properties reduce to the matrix estimates of Lemma \ref{schwartz-convolution-lemma}. Choose a compact fundamental segment $W\subset \widetilde Y$ for the positive periodic end, and write $W_j:=T^jW$ for $j\geq 0$. Since $\rho$ agrees with $f$ sufficiently far along the end, and since $f(T^j x)=f(x)+j$, we can find constants $c,C>0$ such that, for every $j\geq 0$ and every $x\in W_j$,
\[
c(1+j)\leq 1+\rho(x)\leq C(1+j).
\]
Thus the absolute Schwartz condition
\[
\sup_{x,y\in M} \langle \rho(x)+\rho(y)\rangle^N
|\nabla_x^\alpha\nabla_y^\beta K_R(x,y)| <\infty
\]
is equivalent, on the end, to the following family of segmentwise estimates: for every $\alpha,\beta,N_1,N_2$,
\[
\sup_{\substack{x\in W_j\\ y\in W_\ell}}
|\nabla_x^\alpha\nabla_y^\beta K_R(x,y)| \leq C_{\alpha,\beta,N_1,N_2} \langle j\rangle^{-N_1}
\langle \ell\rangle^{-N_2}, \hspace{3mm} j,\ell\geq 0.
\]
Similarly, if $\widetilde A \in \Psi_{\mathbb Z}^{-\infty}(\widetilde Y;\widetilde E)$, then the relative Schwartz estimates defining $\Psi_{\mathbb Z}^{-\infty}$ imply that, for every $\alpha,\beta,N$,
\[
\sup_{\substack{x\in W_j\\ y\in W_\ell}}
|\nabla_x^\alpha\nabla_y^\beta K_{\widetilde A}(x,y)| \leq C_{\alpha,\beta,N} \langle j-\ell\rangle^{-N}.
\]
As a result, the estimates required for the composition properties reduce to the matrix estimates of Lemma \ref{schwartz-convolution-lemma}.

To see how the third composition property follows from the third statement in Lemma \ref{schwartz-convolution-lemma}, note that the kernel of the difference is
\[
K_{(\chi\widetilde A\chi)(\chi\widetilde B\chi) - \chi(\widetilde A\widetilde B)\chi}(x,y) = \chi(x)\chi(y) \int_{\widetilde Y} K_{\widetilde A}(x,z) q(z) K_{\widetilde B}(z,y)\,dz,
\]
where $q$ is a bounded function such that $q(z) = 0$ for $z$ sufficiently far along the positive periodic end. Thus, the result follows by applying Lemma \ref{schwartz-convolution-lemma} with
\[
a_{jk} := \sup_{\substack{x\in W_j\\ z\in W_k}} |\nabla_x^\alpha K_{\widetilde A}(x,z)|, \hspace{3mm} 
b_{k\ell} := \sup_{\substack{z\in W_k\\ y\in W_\ell}}
|\nabla_y^\beta K_{\widetilde B}(z,y)|, \hspace{3mm} q_k:=\sup_{z\in W_k}|q(z)|.
\]

\end{proof}

\printbibliography

\end{document}